\documentclass{siamltex}
\usepackage{fullpage,enumerate,setspace,changepage,multirow,rotating,xcolor}
\usepackage{graphicx,amsmath,amssymb,amsfonts,url, algorithm, algorithmic}
\usepackage[small]{caption}
\def\grad{\nabla}

\def\bo{\mathbf{0}}

\def\cA{\mathcal{A}}

\def\cC{\mathcal{C}}

\def\cE{\mathcal{E}}
\def\cF{\mathcal{F}}

\def\cK{\mathcal{K}}

\def\cN{\mathcal{N}}
\def\cO{\mathcal{O}}

\def\cQ{\mathcal{Q}}

\def\cS{\mathcal{S}}

\def\cX{\mathcal{X}}

\def\mR{\mathbb{R}}

\def\eq{\[}
\def\en{\]}

\def\smskip{\smallskip}

\def\texitem#1{\par\smskip\noindent\hangindent 25pt
               \hbox to 25pt {\hss #1 ~}\ignorespaces}

\def\abs#1{\left|#1\right|}
\def\norm#1{\|#1\|}

\newcommand{\BEAS}{\begin{eqnarray*}}
\newcommand{\EEAS}{\end{eqnarray*}}
\newcommand{\BEA}{\begin{eqnarray}}
\newcommand{\EEA}{\end{eqnarray}}
\newcommand{\BEQ}{\begin{eqnarray}}
\newcommand{\EEQ}{\end{eqnarray}}
\newcommand{\BIT}{\begin{itemize}}
\newcommand{\EIT}{\end{itemize}}
\newcommand{\BNUM}{\begin{enumerate}}
\newcommand{\ENUM}{\end{enumerate}}

\newcommand{\BA}{\begin{array}}
\newcommand{\EA}{\end{array}}

\newcommand{\ones}{\mathbf 1}

\newcommand{\reals}{\mathbb{R}}
\newcommand{\integers}{\mathbb{Z}}

\newcommand{\Rank}{\mathop{\bf rank}}
\newcommand{\Tr}{\mathop{\bf Tr}}

\newcommand{\argmin}{\mathop{\rm argmin}}

\newcommand{\dom}{\mathop{\bf dom}}

\newcommand{\intr}{\mathop{\bf int}}

\newcommand{\bdry}{\mathop{\bf bd}}
\newcommand{\vect}{\mathop{\bf vec}}

\newif\ifpagenumbering
\pagenumberingtrue

\pagenumberingfalse

\newsavebox{\remarkbox}
\savebox{\remarkbox}{\noindent\bf Remark}
\def\fprod#1{\left\langle#1 \right\rangle}
\def\etak{\eta_1^{(k)}}
\def\etakii{\eta_2^{(k)}}
\def\xbp{x_{\ast}}

\def\lk{\lambda^{(k)}}
\def\tk{\tau^{(k)}}
\def\epsk{\epsilon^{(k)}}

\def\Xk{X^{(k)}}
\def\Xkopt{X^{(k)}_{\ast}}
\def\Xopt{X_{\ast}}
\def\Nk{N^{(k)}}
\def\Sk{\cS^{(k)}}

\def\Pk{P^{(k)}}

\def\sopt{s_{\ast}}
\def\yopt{y_{\ast}}
\def\sk{s^{(k)}}
\def\yk{y^{(k)}}
\def\skopt{s^{(k)}_{\ast}}
\def\ykopt{y^{(k)}_{\ast}}
\def\tetki{\theta_1^{(k)}}
\def\tetkii{\theta_2^{(k)}}
\def\xik{\xi^{(k)}}
\def\alphas{\alpha^\ast}\def\betas{\beta^\ast}

\def\fprod#1{\left\langle#1\right\rangle}

\newtheorem{assumption}[theorem]{Assumption}

\title{A Unified Approach for Minimizing Composite Norms}
\author{N. S. Aybat\thanks{IEOR Department, Columbia University.
    Email: {\tt nsa2106@columbia.edu}} \and G. Iyengar\thanks{IEOR
    Department, Columbia University.
    Email: {\tt gi10@columbia.edu}}}
\begin{document}
\maketitle
\begin{abstract}
We propose a first-order augmented Lagrangian algorithm (FALC) to solve
the composite norm minimization problem
\eq
\begin{array}{rl}
  \min_{X\in\reals^{m\times n}} &
  \mu_1\norm{\sigma(\cF(X)-G)}_\alpha+\mu_2\norm{\mathcal{C}(X)-d}_\beta,\\
  \mbox{subject to} & \mathcal{A}(X)-b\in\cQ,
\end{array}
\en
where $\sigma(X)$ % :\reals^{m\times n}\rightarrow\reals_+^{\min\{m,n\}}$
% is a function returning
denotes the vector of singular values of $X \in \reals^{m\times
  n}$, the matrix norm
$\norm{\sigma(X)}_{\alpha}$ denotes either the Frobenius, the nuclear, or
the $\ell_2$-operator norm of $X$, the vector norm $\norm{.}_{\beta}$
denotes either the $\ell_1$-norm, $\ell_2$-norm or the
$\ell_{\infty}$-norm; $\cQ$ is a closed convex set and $\cA(.)$, $\cC(.)$,
$\cF(.)$ are linear operators from
$\reals^{m\times n}$ to vector spaces of appropriate dimensions.
Basis pursuit, matrix
completion, robust principal component pursuit~(PCP), and stable PCP
problems are all special cases
of the composite norm minimization problem. Thus, FALC is able to solve
all these problems in a unified manner.
We show that any limit point of FALC iterate sequence is an optimal solution
of the composite norm minimization problem.
% \comment{GARUD: for many special cases we don't need to compute svd,
% hence I changed SVD to ``constrained shrinkage"
%whenever the optimal solution is unique.
We also show that
%there exists and we can construct an a-priori fixed sequence
%$\{\lk\}_{k\in\integers_+}$ such that
for {\em all} $\epsilon>0$, the % iterates computed by
FALC iterates are $\epsilon$-feasible
% i.e. $\norm{\mathcal{A}(X^{(k)}) - b}_2 \leq \epsilon$,
and $\epsilon$-optimal
% ,
% $\left|~(\mu_1\norm{X^{(k)}}_*+\mu_2\norm{\mathcal{C}(X^{(k)})-d}_1) -
%   (\mu_1\norm{X_*}_*+\mu_2\norm{\mathcal{C}(X_*)-d}_1)\right| \leq
% \epsilon$,
after $\cO(\log(\epsilon^{-1}))$ iterations, which require
$\mathcal{O}(\epsilon^{-1})$ % opera
% tions in the worst case. The
% complexity of each operation is
% dominated by computing one
constrained shrinkage operations and Euclidean projection onto the
set $\cQ$. Surprisingly, on the problem sets we tested, FALC required
only $\cO(\log(\epsilon^{-1}))$ constrained shrinkage,  instead of the
$\cO(\epsilon^{-1})$ worst case bound, to
compute an $\epsilon$-feasible and $\epsilon$-optimal solution.
To best of our knowledge, FALC is the first algorithm with a known
complexity bound that solves the stable PCP problem.
\end{abstract}
\section{\normalsize Introduction}
In this paper we
consider the class of composite norm minimization problems defined in
\eqref{ch3_eq:composite_norm_minimization}.
\begin{align}
\min_{X\in\reals^{m\times
    n}}\mu_1\norm{\sigma(\cF(X)-G)}_\alpha+\mu_2\norm{\mathcal{C}(X)-d}_\beta
\hbox{ subject to }
\mathcal{A}(X)-b\in\cQ, \label{ch3_eq:composite_norm_minimization}
\end{align}
where $\mu_1,\mu_2 \geq 0$, $b\in \mR^{q}$, $\mathcal{A}: \reals^{m\times
  n}\rightarrow \mR^{q}$ denotes a linear map, $\cQ \subset
\mR^{q}$ is a nonempty, closed convex set,
 $d\in \mR^{p}$, $\mathcal{C}: \reals^{m\times n}\rightarrow
\mR^{p}$ is a linear map, $G\in \mR^{r_1\times r_2}$,
$\mathcal{F}: \reals^{m\times n}\rightarrow
\mR^{r_1\times r_2}$ is a linear map, and the function $\sigma(Z)$
denotes the singular values of the matrix $Z$.  The parameter
$\beta\in\{1,2,\infty\}$, and $\norm{.}_{\beta}$ denotes  the
$\ell_{\beta}$ vector norm. The parameter $\alpha \in \{1,2,\infty\}$,
and for $\alpha = 1, 2, \infty$, the norm $\norm{\sigma(Z)}_{\alpha}$ denotes,
respectively, the nuclear norm $\norm{X}_{\ast} = \norm{\sigma(X)}_1$, the Frobenius
norm $\norm{X}_{F} = \norm{\sigma(X)}_2$, and the $\ell_2$-operator norm
$\norm{X}_2 = \norm{\sigma(X)}_\infty$.
%   and
% $\cQ\subset\cE_1$
% where $\sigma(X)\in\reals_+^{\min\{m,n\}}$ denotes the vector of
% singular values of the matrix $X
% \in \reals^{m \times n}$. Here the affine mapping $\cF(.)-G$ maps
% $X$ to a symmetric matrix.
Except for Section~\ref{ch3_sec:extensions} we assume the following.
\begin{assumption}
\label{ass:injective}
The linear map $\cA$ is surjective,   and,
at least, one of the linear maps
$\cC$ and $\cF$ in the objective function is injective.
\end{assumption}

Since at least one of the linear maps $\cC$ and $\cF$ are injective, it
follows that the objective function
$\mu_1\norm{\sigma(\cF(X)-G)}_\alpha+\mu_2\norm{\mathcal{C}(X)-d}_\beta$
of \eqref{ch3_eq:composite_norm_minimization} goes to  $\infty$ as
$\norm{X}_F\rightarrow \infty$, i.e. Assumption~\ref{ass:injective}
ensures that the objective function of
\eqref{ch3_eq:composite_norm_minimization} is coercive, and hence, an
optimal solution to \eqref{ch3_eq:composite_norm_minimization}
exists.

The composite norm minimization
problem~\eqref{ch3_eq:composite_norm_minimization} appears in the context
of ``structured" or ``sparse" optimization  where desired solution is
``structured'' in some form -- the solution matrix may be sparse, i.e.
has very few non-zero components, or it may be low rank, or the indices of
its non-zero coefficients may
all belong to a union of few given index sets, i.e. ``groups''.
We show in Section~\ref{sec:special-cases} that many well-known ``structured"
optimization problems such as basis pursuit, matrix
completion, robust principal component pursuit~(PCP), and stable PCP, are
all special cases
of~\eqref{ch3_eq:composite_norm_minimization}. Moreover,
Assumption~\ref{ass:injective} is satisfied in all these special cases.

Composite norm minimization problem
\eqref{ch3_eq:composite_norm_minimization} can be reformulated as a
semidefinite programming problem~(SDP); hence, it can be solved
efficiently in theory. However, instances of
\eqref{ch3_eq:composite_norm_minimization} that arise in practice are very
large and typically dense. Therefore, interior point based SDP solvers
perform very poorly on these instances.
% Under these assumptions, in Section~\ref{ch_sec:theory}, we propose
% a first-order augmented Lagrangian algorithm~(FALC) that can solve
% \eqref{ch3_eq:composite_norm_minimization}.
% Please note that

\subsection{New Results}
We propose a first-order augmented Lagrangian algorithm~(FALC) to solve
\eqref{ch3_eq:composite_norm_minimization}. The main results of this paper
are as follows:
\begin{enumerate}[(a)]
\item We establish that every  limit
  point $\bar{X}$ of the sequence of FALC iterates
  $\{\Xk\}_{k\in\integers_+}$ is an optimal solution of
  \eqref{ch3_eq:composite_norm_minimization}, i.e.,
  \begin{align}
    \bar{X} \in \argmin_{X\in\reals^{m\times
        n}}\Big\{\mu_1\norm{\sigma(\cF(X)-G)}_\alpha+\mu_2\norm{\mathcal{C}(X)-d}_\beta:\
    \mathcal{A}(X)-b\in\cQ\Big\}. \nonumber
  \end{align}
\item Let $P^\ast$ denote the optimal value of
    \eqref{ch3_eq:composite_norm_minimization}.
  %We explicitly construct a parameter sequence such that
  For {\em all} $\epsilon>0$, the FALC iterates $X^{(k)}$ are
  $\epsilon$-feasible, i.e. there exists $y^{(k)} \in \cQ$ such that
  \eq
   \norm{\mathcal{A}(X^{(k)}) - y^{(k)}- b}_2 \leq \epsilon,
   \en
   and $\epsilon$-optimal, i.e.
   \eq
   \big|~\big(\mu_1\norm{\sigma(\cF(X^{(k)})-G)}_\alpha
       +\mu_2\norm{\mathcal{C}(X^{(k)})-d}_\beta\big)
     -P^\ast\big|
   \leq \epsilon,
   \en
   after $\cO(\log(\epsilon^{-1}))$ FALC iterations that requires
   $\cO(\epsilon^{-1})$ projections on to $\cQ$ and
   $\cO(\epsilon^{-1})$ ``constrained shrinkage" operations
   % \textbf{Algorithm~APG} iterations
  in the worst case - see \eqref{ch3_eq:matrix_constrained_shrinkage} and
  \eqref{ch3_eq:vector_constrained_shrinkage} for the definition and
  complexity of each ``constrained
  shrinkage" operation.
\item FALC can be extended to solve the following more general optimization
  problem
  \begin{equation}
    \label{ch3_eq:general_problem}
    \begin{array}{rl}
      \mbox{min}_{X\in \reals^{m\times n}} &\mu_1 \norm{\sigma(\cF(X)-G)}_{\alpha}  +
      \mu_2\norm{\cC(X)-d}_{\beta}+\mu_3 H(X)\; \mbox{ s.t.}\; \mathcal{A}(X)-b\in\cQ,
    \end{array}
  \end{equation}
  where $H:\reals^{m\times n}\rightarrow \reals$ is a strongly convex
  function, with the same complexity guarantees.
\item In our numerical tests we observed that FALC required only
    $\cO(\log(\epsilon^{-1}))$  projection and
    shrinkage operations
    % \textbf{Algorithm~APG} iterations
  as opposed to $\cO(\epsilon^{-1})$ the worst case theoretical bound
  proven in the paper to obtain an $\epsilon$-feasible and
  $\epsilon$-optimal solution.
\item We also observed that, although,  FALC is a general-purpose  algorithm
  for the  composite norm minimization
  problem~\eqref{ch3_eq:composite_norm_minimization}, the numerical
  results show that FALC is competitive with the state-of-the-art special purpose
  algorithms designed for all  special cases that we tested.
\end{enumerate}

%  and show that FALC can be extended
% easily to solve problem~\eqref{ch3_eq:general_problem} with the same
% complexity guarantees.  for the optimization
% problem~\eqref{ch3_eq:composite_norm_minimization} and later in
% Section~\ref{ch3_sec:extensions},

% In this paper we establish the following properties of the FALC algorithm.
%  % %\textbf{Algorithm~APG} iteration
%   % is determined by %$\cO(\min\{nm^2,n^2m\})$
%   % the complexity of a singular value decomposition~(SVD) when
%   % $\mu_1>0$ or by complexity of matrix-vector multiplications when
%   % $\mu_1=0$ and $\mu_2>0$. For the definition of ``constrained
%   % shrinkage" operations,

% \end{enumerate}

\subsection{Special Cases}
\label{sec:special-cases}
We show below that many well studied ``structured'' optimization problems are
special cases of~\eqref{ch3_eq:composite_norm_minimization}.
\paragraph{\textbf{Nuclear norm-minimization}}
The {\em nuclear norm minimization problem}
\begin{align}
\min_{X\in\reals^{m\times n}}\norm{X}_* \hbox{ subject to~(s.t.) }
\cA(X) = b. \label{ch3_eq:nuclear_minimization}
\end{align}
is a special case of (\ref{ch3_eq:composite_norm_minimization})
with $\cQ=\{\bo\}$, $\cF(X)=X$, $G=\bo$, $\alpha=1$, i.e. $\norm{\sigma(X)}_1 =
\norm{X}_{\ast}$, $\mu_1 = 1$, and $\mu_2 = 0$.
The nuclear norm minimization problem is a convex approximation for the
NP-hard rank minimization problem $\min_{X\in\reals^{m\times
    n}}\{\rank(X):\ \cA(X) = b\}$,
where $\rank(X)$ denotes the rank of $X\in\reals^{m\times n}$.
Rank
minimization arises in many different contexts, e.g. system
identification \cite{Liu08_1J}, optimal
control~\cite{Faz03_1P,Faz04_1P,Gha93_1P}, low-dimensional
embedding in Euclidean space~\cite{Lin95_1J}, and matrix
completion~\cite{Can08_1J}.

Let $X_0\in\reals^{m\times n}$ be the unknown low-rank matrix such that $\cA(X_0) = b$. Let $r=\rank(X_0)$ and $\bar{n}=\max\{m,n\}$. When
the linear operator $\mathcal{A}:\reals^{m\times n}\rightarrow\reals^q$ satisfies some
regularity properties, and the number of measurements $q = \cO(r(m+n)\log(mn))$,  Recht et al.~\cite{Recht07_1J}
show that, with very high probability, \eqref{ch3_eq:nuclear_minimization} has a unique
optimal solution  and  this solution is also optimal for the rank minimization problem.

Another related special case is the \emph{matrix completion problem}
where the operator $\cA$ picks a subset of the matrix elements,
i.e., the linear constraints are of the form: $X_{ij} = (X_0)_{ij}$ for
$(i,j) \in \Omega$, where $\Omega$ is a given index set of observable entries of an unknown low rank matrix $X_0\in\reals^{m\times n}$.
When indices $(i,j)$ are sampled uniformly at random, and
$|\Omega|=\cO(\bar{n}^{1.2}r\log(\bar{n}))$ and the unknown matrix $X_0$
satisfies some regularity conditions, Cand{\'e}s et al.
\cite{Can08_1J} show that, with high probability, $X_0$ is the unique
solution of the matrix completion problem. The Netflix prize
problem~\cite{Netflix} is an
example of the matrix completion problem.

\eqref{ch3_eq:nuclear_minimization} can be reformulated as an SDP; however, instances of
\eqref{ch3_eq:nuclear_minimization} that arise in practice are so large that standard SDP solvers are unable to
solve them.
For existing algorithmic methodologies for solving the nuclear norm
minimization problem, see
\cite{Can08_2J,Gold10_1J,Ma09_1J,Ma09_1R,Ma08_1J,Toh09_1J} and references
therein.
\paragraph{\textbf{Basis-pursuit problem}}
The {\em basis pursuit problem}
\begin{align}
\min_{x\in\reals^n} \norm{x}_1 \hbox{ s.t. } Ax=b, \label{ch3_eq:l1_minimization}
\end{align}
where $A\in\reals^{q\times n}$ and $b\in\reals^q$,
is a  special case of \eqref{ch3_eq:composite_norm_minimization} with
$\cQ=\{\bo\}$, $C(x)=x \in\reals^{n\times 1}$,
$d=0$, $\beta = 1$, $\mu_1 = 0$ and $\mu_2 = 1$. The basis pursuit problem has
attracted a lot of attention
recently, since it appears in the context of {\em compressed sensing}~(CS)~\cite{Can06_1J, Can06_2J, Can06_3J,
  Don06_1J}. The goal in CS is to recover a sparse signal $x_0\in\reals^n$ from a small
set of linear measurements or transform values $b = Ax_0$, or equivalently,
to solve the NP-hard $\ell_0$-minimization problem
\begin{align}
\min_{x\in\reals^n}\|x\|_0 \hbox{ s.t. } Ax=b, \label{ch3_eq:l0_minimization}
\end{align}
where the $\ell_0$-norm $\norm{x}_0 = \sum_{i=1}^n \ones(x_i \neq 0)$
and $\ones(\cdot)$ is equal to 1 if the argument is true, 0
otherwise.
Recently, Cand\'{e}s, Romberg and Tao~\cite{Can06_1J, Can06_2J, Can06_3J} and
Donoho~\cite{Don06_1J} have shown that,
if the target signal $x_0$ is $s$-sparse, i.e. %only $s$ of the $n$
                                %components are non-zeros
$\norm{x_0}_0=s$, the matrix $A \in \reals^{q \times n}$ has $q =
\cO(s\log(n))$ and is chosen randomly according to a specific
set of distributions, then, with very high probability, the  sparse target signal $x_0$
is the unique optimal solution of the basis pursuit
problem~(\ref{ch3_eq:l1_minimization}).
Thus, $x_0$ can be recovered by solving
a linear program~(LP), and therefore, in theory, signal recovery is very efficient.
In practice, however, simplex and interior point based
general purpose LP solvers are unable to solve such LPs efficiently because the matrix $A$
in \eqref{ch3_eq:l1_minimization} is large, dense, and %in addition, these LPs are
often ill-conditioned.
The measurement matrix $A$ in many CS applications has a lot of
structure, in particular, the
matrix-vector multiplication $Ax$ and $A^Ty$ can be computed efficiently in $\cO(n\log(n))$ time.
Recently, a number of different algorithms have been proposed to
exploit this
structural fact to efficiently
solve~\eqref{ch3_eq:l1_minimization}~\cite{Aybat12,AybatI09:SPA,Dau08_1J,Wri07_1J,Yin07_1R,Yin08_1J,Boy07_1R,Ber08_1J,Wen09_1R,Yin08_2J}
.
Note that the ``noisy'' basis pursuit problem
\[
\min_{x\in\reals^n}\|x\|_1 \hbox{ s.t. } \norm{Ax-b}_2 \leq \delta,
\]
is a special case of \eqref{ch3_eq:composite_norm_minimization} with
$\cQ = \{y\in\reals^q: \norm{y}_2 \leq \delta\}$.
%Please refer to Section~\ref{ch1_sec:intro} for a brief survey of the previous work on $\ell_1$-minimization to obtain sparse solution.
\paragraph{\textbf{Principal component pursuit}}
The
\emph{principal component pursuit problem}
\begin{align}
\min_{X\in\reals^{m\times n}} \norm{X}_* + \mu_2 \norm{\vect(X)-d}_1, \label{ch3_eq:component_pursuit}
\end{align}
is a special case of \eqref{ch3_eq:composite_norm_minimization} with
$\cA = \bo$, $b=0$, $\cQ=\{\bo\}$,
$\cF(X)=X$, $G=\bo$, $\alpha=1$, i.e. $\norm{\sigma(X)}_\alpha =
\norm{X}_{\ast}$,
$\beta=1$, %$\mu_1=1$, $\mu_2>0$,
and %  the
%  $\mathcal{C}:\reals^{m\times n}\rightarrow \reals^{mn}$ is
% given by
$\mathcal{C}(X)= \vect(X)$, where $\vect(X)$ is the vector obtained by
stacking the columns of $X\in\reals^{m\times n}$ in order, and $\beta = 1$.
Suppose the data matrix
$D\in\reals^{m\times n}$ is of the form $D = X_0 + S_0$, where $X_0$ is a low
rank matrix, $S_0$ is a sparse matrix, and both satisfy some
regularity conditions given in \cite{Can09_1J,Ma09_1J}. Then the low
rank and sparse components of $D$ can be recovered by solving
(\ref{ch3_eq:component_pursuit}) with $d=\vect(D)$ and
$\mu_2=1/\sqrt{\bar{n}}$, where $\bar{n}=\max\{m,n\}$~\cite{Can09_1J,Ma09_1J}.
For existing algorithmic approaches for solving principal component pursuit, see ~\cite{Can09_1J,Gold10_1J,Ma09_1J,Ma09_1R,sRPCA} and references therein.

In \cite{sRPCA}, it is shown that recovery is still possible
even when the data matrix $D$ is corrupted by a dense error matrix.
Suppose  the data matrix $D$ is of the form
$D = X_0 + S_0 + \zeta_0$, where % $X_0$ is a low rank
% matrix, $S_0$ is a sparse matrix and
$(\zeta_0)_{ij}$ is independent and
identically distributed~(i.i.d.) for all $i,j$ such that
$\norm{\zeta_0}_F\leq\delta$.
Then the optimal solution $(X_*,S_*)$ of the
\emph{stable principal component pursuit problem}
\begin{equation}
\label{ch3_eq:stable_component_pursuit}
  \begin{array}{rl}
    \mbox{min}_{X,S \in \reals^{m\times n}} &\norm{X}_* + \frac{1}{\sqrt{\bar{n}}}\norm{\vect(S)}_1,\\
    \mbox{s.t.} &  \norm{X+S-D}_F\leq\delta,
  \end{array}
\end{equation}
satisfies $\norm{X_*-X_0}_F^2+\norm{S_*-S_0}_F^2\leq
c~mn\delta^2$ for some constant $c$ with high probability. % , provided that $X_0$ and $S_0$ satisfy some regularity conditions.
The stable principle component pursuit is a special case of
\eqref{ch3_eq:composite_norm_minimization} with $\cA(X,S) =
\vect(X+S)$, $b = \vect(D)$, $\cQ = \{y\in\reals^{mn}: \norm{y}_{2} \leq \delta\}$.

\paragraph{\textbf{Composite norm minimization with conic constraints}}
The goal in the {\em minimal system realization problem} is to design
the lowest order discrete-time, linear time-invariant~(LTI) dynamical system
that is consistent with the observed
data. Let $x_i$ be the true (unknown) impulse
response of the system at time $i$ for $ i = 1,
\ldots, n$. Suppose that we observe noisy data $\tilde{x}_i=x_i+\varepsilon_i$, $ i = 1,
\ldots, n$, where $\{\varepsilon_i\}$ are i.i.d. uniform over $[-\varrho,\varrho]$. It is well-known~\cite{Faz03_1J,Faz11_1J} that
the minimum order system consistent with the observations can be computed by
solving
\eq
\min_{x \in \reals^{2n-1}} \rank(H_n(x)) \mbox{ s.t. } \norm{x_1^n -
  \tilde{x}_1^n}_{\infty} \leq \delta,
\en
where $H_n(x)$ be the Hankel matrix formed by
$x\in\reals^{2n-1}$ and $x_1^n$ (resp. $\tilde{x}_1^n$) denotes the vector formed by the first $n$ components
of $x$ (resp. $\tilde{x}$).
% From \cite{Faz03_1J,Faz11_1J}, it is known that
% and with
% order $r$ if and only if
% $r=\min_{\{x_i\}_{i=n+1}^{2n-1}}\rank(H_n(x))$.
This rank-minimization problem can be approximated by the nuclear norm
minimization problem
\begin{align}
\min_{x \in \reals^{2n-1}} \norm{H_n(x)}_{\ast} \mbox{ s.t. } \norm{x_1^n -
  \tilde{x}_1^n}_{\infty} \leq \delta. \label{ch3_eq:hankel}
\end{align}
% \begin{eqnarray}
% \min_{x\in\reals^{2n-1}} &\norm{H_n(x)}_* \mbox{ subject to }
% |x_i-\tilde{x}_i|\leq \varrho,\; i=1,\ldots,n. \label{ch3_eq:hankel}
% \end{eqnarray}
Since $H_n(.)$ is  injective, $\cA(x)=[x_1,\ldots,x_n]^T$ is
surjective, and $\cQ=\{y\in\reals^n:\ \norm{y}_\infty\leq \delta\}$
is a closed convex set, \eqref{ch3_eq:hankel} is a special case of
\eqref{ch3_eq:composite_norm_minimization}.

In the {\em sparse} PCA problem
\begin{eqnarray}
\min_{x\in\reals^n} &\norm{\Sigma-xx^T}_F \mbox{ s.t. }
\norm{x}_0\leq s, \label{ch3_eq:sparse_pca}
\end{eqnarray}
% is interested in computing a
the goal is to compute an $s$-sparse vector $x$ that
is ``close'' to the eigenvector corresponding to the largest eigenvalue
of the positive semidefinite matrix $\Sigma$.
% where $\mu_2 > 0$ controls the trade-off between sparsity and
% approximation error.
% The optimization problem
% (\ref{ch3_eq:sparse_pca}) is not convex, and
% is, therefore, hard to solve.
Let $X = xx^T$. Then \eqref{ch3_eq:sparse_pca}
is equivalent to
\eq
\min_{X\in\reals^{m\times n}} \norm{X-\Sigma}_F \mbox{ s.t. }
\norm{\vect(X)}_0\leq s^2, \rank(X)=1, X \succeq 0.
\en
% in the sparse PCA
% problem we require that $\rank(X) = 1$
% and $\vect(X)$ be sparse.
Since  $\norm{X}_{\ast}$ is the tightest convex upper bound for $\rank(X)$, and
$\norm{X}_{\ast} = \Tr(X)$ for positive semidefinite~(psd) matrices,
the convex relaxation
\begin{equation}
  \begin{array}{rl}
    \min_{X\in\reals^{n \times n}} &\norm{X-\Sigma}_F +
    \mu \norm{\vect(X)}_1+\nu\fprod{I, X},\label{ch3_eq:sPCA}\\
    \mbox{s.t.} & X\succeq 0,
  \end{array}
\end{equation}
for~\eqref{ch3_eq:sparse_pca}, where $\mu$ and $\nu$
control the sparsity on the entries and the singular values of $X$,
respectively, is a special case of
\eqref{ch3_eq:composite_norm_minimization} with $\cQ$ set as the cone
of psd matrices~(the linear term in the objective function over psd
cone can be handled easily).
See \cite{Mont08_1J, Mont07_1J,Nest10_1J} for existing %convex relaxation
approaches for solving the sparse PCA problem.

In Section~\ref{ch3_sec:extensions}, we show that FALC can be easily
extended to solve problems of the form given in \eqref{ch3_eq:general_problem}.
When $F \not \equiv 0$, we do not require that $\cC$ or
$\cF$ be injective. % or that $\max\{\mu_1,\mu_2\} > 0$
Thus, FALC is
able to solve regularized conic optimization problems of the form
\[
\min_{X \in \reals^{m\times n}}
\fprod{R,X}+\rho\norm{X}_F^2 \hbox{ s.t. } B \preceq \cA(X),
\]
where $\cA$ is a  surjective map.

This paper is organized as follows. In Section~\ref{ch3_sec:theory} we prove
the main convergence results for FALC and in
Section~\ref{ch3_sec:implementation} we discuss all the implementation details. In Section~\ref{ch3_sec:computations} we report the results from our
numerical experiments comparing FALC with other algorithms to solve principle component pursuit problems. Finally, in Section~\ref{ch3_sec:extensions}, we briefly discuss the general problem \eqref{ch3_eq:general_problem} and conclude.

\section{Preliminaries}
\label{sec:preliminary}
\begin{figure}[t]
    \rule[0in]{6.5in}{1pt}\\
    \textbf{Algorithm APG}$(p, f, \cS, x^{(0)}, \textsc{ITERstop}, \textsc{GRADstop})$\\
    \rule[0.125in]{6.5in}{0.1mm}
    \vspace{-0.25in}
    {\footnotesize
    \begin{algorithmic}[1]
    %\STATE \textbf{input:} a prox function $h(.)$, $X^{(0)} \in \dom p$ and $F\subset\reals^{m\times n}$ \textbf{such that}\\
    %\hspace*{0.5in}$X_s\in F$ for some $X_s\in\argmin_{x\in~\reals^{m\times n}}p(X)+f(X)$
    \STATE $x_1^{(0)}\gets x^{(0)}$, $x_2^{(1)} \gets x^{(0)}$, $t^{(1)}\gets 1$, $\ell\gets 0$
    \WHILE{\textsc{ITERstop}$(\ell)$ \textbf{and} \textsc{GRADstop}$(x_1^{(\ell)})$ are \FALSE} \label{ch1_algeq:stop}
    \STATE $\ell \gets \ell + 1$
    \STATE $x_1^{(\ell)} \gets \argmin\left\{p(x)+\fprod{\nabla f(x_2^{(\ell)}),~x-x_2^{(\ell)}} + \frac{L}{2}\norm{x-x_2^{(\ell)}}: x\in \cS\right\}$ \label{ch1_algeq:tseng_z}
    \STATE $t^{(\ell+1)}\gets \left(1+\sqrt{1+4~\left(t^{(\ell)}\right)^2}\right)/2$
    \STATE $x_2^{(\ell+1)} \gets x_1^{(\ell)} + \left(\frac{t^{(\ell)}-1}{t^{(\ell+1)}}\right)\left(x_1^{(\ell)}-x_1^{(\ell-1)}\right)$
    \ENDWHILE
    \RETURN $ x_1^{(\ell)}$
    \end{algorithmic}
    \rule[0.25in]{6.5in}{0.1mm}
    }
    \vspace{-0.5in}
    \caption{Accelerated Proximal Gradient Algorithm}\label{ch1_alg:pga}
    \vspace{-0.25in}
\end{figure}
In this section we state and briefly discuss the details of a particular
implementation of Fast Iterative Shrinkage-Thresholding
Algorithm~(FISTA)~\cite{Beck09_1J} that we use as a subroutine in
FALC. Let $\cE$ be a Hilbert space and
$\norm{.}=\sqrt{\fprod{.,.}}$.  FISTA computes an
$\epsilon$-optimal solution of
\begin{align}
\label{ch1_eq:tseng_problem}
\min_{x\in\cE} p(x)+f(x),
\end{align}
in $\cO(1/\epsilon)$ iterations, where $p:\cE \rightarrow\reals$ and $f:\cE \rightarrow\reals$ are
continuous convex functions such that $\grad f$ is Lipschitz
continuous on $\cE$ with constant $L$. % \eqref{ch1_eq:tseng_problem} is
% obtained.
%, i.e. $p(\bar{x})+f(\bar{x})\leq \inf_{x\in\cE} \{p(x)+f(x)\} +\epsilon$.
Later, Tseng~\cite{Tseng08} showed that this rate result for FISTA also holds
when $p:\cE \rightarrow (-\infty, +\infty]$ and $f:\cE \rightarrow
(-\infty, +\infty]$ that are proper, lower semicontinuous, and convex
functions such that $\dom p$ is closed and $\grad f$ is Lipschitz
continuous on $\cE$.
Moreover, for any given convex set $\cS\subset\cE$ such that $\cS \cap
\argmin_{x \in \cE}\{p(x)+f(x)\} \neq \emptyset$, one can ensure that
the iterate sequence lies in the convex set $\cS$ and all the iterates
are $\epsilon$-optimal after $\cO(1/\epsilon)$ iterations.  We use
this property later
in the paper to uniformly bound the FALC iterates.
%Let $h:\cE \rightarrow\reals$ be a differentiable and strongly convex
%function with convexity parameter $c>0$ with respect to norm $\norm{.}$,
%i.e. $h(y)\geq h(x)+\fprod{\grad h(x), y-x}+\frac{c}{2}\norm{y-x}^2$ for
%all $x,y\in\dom h$. In this paper, we will refer to a function $h(.)$
%with these properties as a {\em prox} function.

\textbf{Algorithm~APG} displayed in Figure~\ref{ch1_alg:pga} takes as
input the functions $f$ and  $p$, the convex set $\cS\subset
\cE$ such that $\cS \cap
\argmin_{x \in \cE}\{p(x)+f(x)\} \neq \emptyset$, an initial iterate $x^{(0)}\in \cE$ and two
stopping criteria \textsc{ITERstop} and \textsc{GRADstop}.
\textbf{Algorithm~APG} is the same as Algorithm~2 in \cite{Tseng08} where we have set
$\cX_\ell \equiv \cS$ (see~\cite{Tseng08} for details). Indeed,
Algorithm~2 in~\cite{Tseng08} is a
modification of FISTA~\cite{Beck09_1J} and reduces to FISTA when
$\cS=\cE$. \noindent FISTA and Algorithm~2 in
\cite{Tseng08} do not use  $\textsc{ITERstop}$ and
$\textsc{GRADstop}$ -- we include them in the definition of
\textbf{Algorithm~APG} because we terminate the algorithm early when we
call it
% depending on $\textsc{ITERstop}$ and $\textsc{GRADstop}$
as a subroutine in FALC. $\textsc{ITERstop}(\ell)$ is a stopping
criterion that  only depends on the current iterate
counter $\ell$ and $\textsc{GRADstop}(x)$ is a stopping criterion that
only depends on its argument $x$.
Lemma~\ref{ch1_lem:tseng_corollary} gives the iteration complexity of
\textbf{Algorithm~APG}.
\begin{lemma}
\label{ch1_lem:tseng_corollary}
Let $p$ and $f$ be a proper, lower semicontinuous, convex functions
such that $\dom p$ is closed and $\grad f$ is Lipschitz continuous on
$\cE$ with constant $L$. Fix $\epsilon>0$ and let
$\{x_1^{(\ell)},x_2^{(\ell)}\}_{\ell\in\integers_+}$ denote the
\textbf{Algorithm~APG} iterates % displayed in
                                % Figure~\ref{ch1_alg:pga} when the
                                % algorithm is not terminated
                                % according to
when both
$\textsc{ITERstop}$ and $\textsc{GRADstop}$ are disabled. Then
$p(x_1^{(\ell)})+f(x_1^{(\ell)})\leq
\min_{X\in\cE}\{p(x)+f(x)\}+\epsilon$ whenever
$\ell\geq\sqrt{\frac{2L}{\epsilon}}~\norm{x^\ast-x^{(0)}}-1$, where $x^\ast
\in \argmin_{x \in \cE}\{p(x) + f(x)\}$.
\end{lemma}
\begin{proof}
See Corollary~3 in \cite{Tseng08} and Theorem 4.4 in \cite{Beck09_1J} for
the details of proof.\end{proof}

\section{FALC Algorithm}
\label{ch3_sec:theory}
\begin{figure}[h!]
    \rule[0in]{6.5in}{1pt}\\
    \textbf{Algorithm FALC}$\big(\big\{\big(\lk,\epsk,\tk,\xik\big)\big\}_{k\in \integers_+}, X^{(0)}\big)$\\
    \rule[0.125in]{6.5in}{0.1mm}
    \vspace{-0.25in}
    {\footnotesize
    \begin{algorithmic}[1]
    %\STATE \textbf{input:} $\big\{\big(\lk,\epsk,\tk,\xik\big)\big\}_{k\in \integers_+}$, $X^{(0)}\in\reals^{m\times n}$ \textbf{such that} $\mathcal{A}(X^{(0)})=b$
    \STATE $y^{(0)}\gets \mathcal{A}(X^{(0)})-b, s^{(0)}\gets \mathcal{C}(X^{(0)})-d$
    \STATE $\eta\gets \mu_1\norm{\sigma(X^{(0)})}_\alpha + \mu_2\norm{\mathcal{C}(X^{(0)})-d}_\beta$
    \STATE $\theta_1^{(1)} \gets 0$, $\theta_2^{(1)} \gets 0$, $k \gets 0$
    \WHILE{($\textsc{FALCstop}$ is \FALSE)}
    \STATE $k \gets k + 1$
    \STATE $p^{(k)}(X,s,y) := \lambda^{(k)}(\mu_1\norm{\sigma(X)}_\alpha+\mu_2\norm{s}_\beta)$
    \STATE $f^{(k)}(X,s,y) := \frac{1}{2} \norm{\mathcal{A}(X)-y-b-\lambda^{(k)}\theta_1^{(k)}}_2^2+\frac{1}{2} \norm{\mathcal{C}(X)-s-d-\lambda^{(k)}\theta_2^{(k)}}_2^2$
    \STATE $\etak\gets \eta+\frac{\lambda^{(k)}}{2}\left(\norm{\tetki}_2^2+\norm{\tetkii}_2^2\right)$ \label{ch3_li:etak}
    \STATE $\cS^{(k)} := \{(X,s,y)\in\reals^{m\times n}\times\reals^p\times\reals^q:\ \mu_1\norm{\sigma(X)}_\alpha\leq\etak,\ \mu_2\norm{s}_\beta\leq\etak\}$ \label{ch3_alg:Fk}
    \STATE $\textsc{ITERstop}(\ell):=\{\ell\geq\ell_{\max}^{(k)}\}$, where $\ell_{\max}^{(k)}$ is defined in
    \eqref{ch3_eq:apg_iter_cond} \label{ch3_li:inner-stopping-condition1}
    \STATE $\textsc{GRADstop1}(X,s,y):=\left\{\exists (G,g)\in\partial_{X,s} P^{(k)}(.,.,.)|_{(X,s,y)} \mbox{ s.t. } \sqrt{\norm{G}_F^2 + \norm{g}_2^2} \leq \tau^{(k)}\right\}$ \label{ch3_li:gradstop1}
    \STATE $\textsc{GRADstop2}(X,s,y):=\left\{\norm{y-\Pi_{\cQ}\left(y-\frac{1}{L}\grad_y P^{(k)}(X,s,y)\right)}_2\leq\xik\right\}$
    \STATE $\textsc{GRADstop}:= \textsc{GRADstop1} \textbf{ and } \textsc{GRADstop2}$\label{ch3_li:inner-stopping-condition2}
    \STATE $(X^{(k)},s^{(k)},\yk)\gets\textbf{Algorithm APG}(p^{(k)}, f^{(k)}, \Sk, \left(X^{(k-1)},s^{(k-1)},y^{(k-1)}\right), \textsc{ITERstop}, \textsc{GRADstop})$ \label{algeq:apg_call}
    \STATE  $\theta_1^{(k+1)} \gets \theta_1^{(k)} - \frac{\mathcal{A}(X^{(k)})-\yk-b}{\lambda^{(k)}}$ \label{ch3_li:lagrangian_update_1}
    \STATE  $\theta_2^{(k+1)} \gets \theta_2^{(k)} - \frac{\mathcal{C}(X^{(k)})-\sk-d}{\lambda^{(k)}}$ \label{ch3_li:lagrangian_update_2}
    \ENDWHILE
    \RETURN $(\Xk,\sk,\yk)$
    \end{algorithmic}
    \rule[0.25in]{6.5in}{0.1mm}
    }
    \vspace{-0.5in}
    \caption{Outline of First-Order Augmented Lagrangian Algorithm~(FALC)}\label{ch3_alg:falc}
    \vspace{-0.25in}
\end{figure}
For the sake of notational simplicity,  we focus on the following
simpler problem in this section.
\begin{align}
\min_{X\in\reals^{m\times n}}\mu_1\norm{\sigma(X)}_\alpha +\mu_2\norm{\mathcal{C}(X)-d}_\beta
\hbox{ subject to }
\mathcal{A}(X)-b\in\cQ. \label{ch3_eq:simple_composite_norm_minimization}
\end{align}
We give convergence results for FALC when $\mu_1>0$ and $\mu_2>
0$. In
Section~\ref{ch3_sec:extensions} we briefly describe how to modify
the algorithm to solve \eqref{ch3_eq:general_problem}. %and also discuss the necessary
% modifications when $\mu_1=0$ or
% $\mu_2=0$.

The linear maps $\mathcal{A}$ and $\cC$ in
\eqref{ch3_eq:simple_composite_norm_minimization} can be represented as
$\mathcal{A}(X)=A\vect(X)$ and $\cC(X) = C \vect(X)$, where
$A\in\reals^{q\times mn}$ and  $C\in\reals^{p\times mn}$. Let
$\sigma_{\min}(A)$ and $\sigma_{\max}(A)$ denote the smallest and the
largest singular values of $A$, respectively. Since we assume that
$\cA$ is surjective~(see Assumption~\ref{ass:injective}), $A$ has full
row rank; consequently, $A^T$ has full column rank. %We also assume that
% $X^{(0)}\in\reals^{m\times n}$ is given,
% such that $\mathcal{A}(X^{(0)})=b$.
We set
\begin{math}
M := \left(
     \begin{array}{ccc}
       -I & 0 & C \\
       0 & -I & A \\
     \end{array}
    \right)
\end{math}
and $L=\sigma^2_{\max}(M)$. Let $\Xopt$ denote an optimal solution of
\eqref{ch3_eq:simple_composite_norm_minimization} and
$\Pi_{\cQ}:\reals^q\rightarrow \cQ$ denote the Euclidean projection onto
$\cQ\subset\reals^q$.

To obtain separable and efficiently solvable subproblems, we introduce
slack variables $s\in\reals^p$ and $y\in\reals^q$, and reformulate
\eqref{ch3_eq:simple_composite_norm_minimization} as
\begin{equation}
  \label{ch3_eq:CN-slack}
  \begin{array}{rlclcl}
    \min_{X, s, y} & \mu_1\norm{\sigma(X)}_\alpha& + & \mu_2\norm{s}_\beta,\\
    \mbox{s.t.}
    &\cC(X) & - &~s &  =& d,\\
    &\cA(X) & - &~y &  =& b,\\
    &       &   &~y &\in&\cQ.
  \end{array}
\end{equation}
We solve \eqref{ch3_eq:CN-slack} by inexactly solving a sequence of
optimization problems  of the form
\begin{equation}
  \label{ch3_eq:augmented_lagrangian_subproblem}
  \min_{X\in\reals^{m\times n},~s\in\reals^p,~
    y\in\cQ\subset\reals^q}
  \left\{
    \begin{array}[c]{l}
      \lambda^{(k)}(\mu_1\norm{\sigma(X)}_\alpha+\mu_2\norm{s}_\beta)\\
      \mbox{} -
      \lambda^{(k)}(\theta_1^{(k)})^T(\mathcal{A}(X)-y-b) +
      \frac{1}{2}\norm{\mathcal{A}(X)-y-b}_2^2 \\
      \mbox{} -
      \lambda^{(k)}(\theta_2^{(k)})^T(\mathcal{C}(X)-s-d) +
      \frac{1}{2}\norm{\mathcal{C}(X)-s-d}_2^2
    \end{array}
  \right\},
\end{equation}
for an appropriately chosen sequence %of parameters
$\{(\lambda^{(k)},\tetki, \tetkii)\}_{k \in \integers_+}$. By completing
squares, it is easy to see that
\eqref{ch3_eq:augmented_lagrangian_subproblem} is equivalent to
\begin{equation}
  \label{ch3_eq:augmented_lagrangian_subproblem_2}
  \min_{X\in\reals^{m\times n},~s\in\reals^p,~y\in\cQ\subset\reals^q}P^{(k)}(X,s,y),
\end{equation}
where
\begin{equation}
  \begin{array}{rcl}
    P^{(k)}(X,s,y) & := &
    \lambda^{(k)}(\mu_1\norm{\sigma(X)}_\alpha+\mu_2\norm{s}_\beta) +
    f^{(k)}(X,s,y), \\
    f^{(k)}(X,s,y) & := &
    \frac{1}{2}\norm{\mathcal{A}(X)-y-b-\lambda^{(k)}\theta_1^{(k)}}_2^2
    +\frac{1}{2}\norm{\mathcal{C}(X)-s-d-\lambda^{(k)}\theta_2^{(k)}}_2^2.
  \end{array}
  \label{ch3_eq:smooth_part}
\end{equation}

\textbf{Algorithm~FALC} displayed in Figure~\ref{ch3_alg:falc} is the
outline of Algorithm FALC. The algorithm takes as inputs the sequence
$\big\{\big(\lk,\epsk,\tk,\xik\big)\big\}_{k\in
  \integers_+}$ and a starting point
$\left(X^{(0)},s^{(0)},y^{(0)}\right)$ such that
$\mathcal{A}(X^{(0)})-b\in \cQ$, $y^{(0)}:=\mathcal{A}(X^{(0)})-b\in\cQ$,
and $s^{(0)}:=\mathcal{C}(X^{(0)})-d$.
In Section~\ref{ch3_sec:multiplier_selection} we describe how we set the input %multiplier
sequence. Let $(\Xkopt,\skopt,\ykopt)$ denote an optimal solution of
\eqref{ch3_eq:augmented_lagrangian_subproblem_2}.
%, i.e. $(\Xkopt,\skopt,\ykopt)\in \argmin_{X\in\reals^{m\times
%n},s\in\reals^p,y\in\reals^q}\{P^{(k)}(X,s,y):~\norm{y}_\gamma\leq\delta\}$.

\subsection{Algorithm APG for subproblems}
In each iteration of FALC, we call \textbf{Algorithm~APG} displayed in
Figure~\ref{ch1_alg:pga} to inexactly solve
\eqref{ch3_eq:augmented_lagrangian_subproblem_2}, which we call the
``$k$-th subproblem".
Let
$\mathbf{1}_\cQ$ denote the indicator function of the closed convex
set $\cQ\subset\reals^q$, i.e., if $y\in\cQ$, then
$\mathbf{1}_\cQ(y)=0$; otherwise,
$\mathbf{1}_\cQ(y)=\infty$. $\mu_1\norm{\sigma(X)}_\alpha
+\mu_2\norm{s}_\beta+\mathbf{1}_\cQ(y)$
is a proper, lower semicontinuous~(lsc), convex function of
$(X,s,y)$. Moreover, $f^{(k)}(X,s,y)$ is a proper, lsc, convex
function that has a Lipschitz continuous gradient, $\grad f^{(k)}$,
defined on $\reals^{m\times n}\times\reals^p\times\reals^q$ with
Lipschitz constant equal to $L$ for all $k\geq 1$. Thus,
\eqref{ch3_eq:augmented_lagrangian_subproblem_2} is of the form
described in  \eqref{ch1_eq:tseng_problem}.
In each update step of \textbf{Algorithm~APG}
i.e. line~\ref{ch1_algeq:tseng_z} in Figure~\ref{ch1_alg:pga}, we need
to solve one problem of the form
%\begin{equation}
%\label{ch3_eq:sub-problem}
%\min_{X\in\reals^{m\times n}, s\in\reals^p,
%  y:\norm{y}_\gamma\leq\delta}\left\{
%  \begin{array}[c]{l}
%    \lambda^{(k)}\mu_1\norm{\sigma(X)}_\alpha + \grad_X
%    f^{(k)}(\hat{X},\hat{s},\hat{y})^T(X-\hat{X})+\frac{L}{2}\norm{X-\hat{X}}_F^2 \\ \smskip
%    \mbox{} + \lambda^{(k)}\mu_2\norm{s}_\beta + \grad_s
%    f^{(k)}(\hat{X},\hat{s},\hat{y})^T(s-\hat{s})
%    +\frac{L}{2}\norm{s-\hat{s}}_2^2,\\
%    \mbox{} + \grad_y f^{(k)}(\hat{X},\hat{s},\hat{y})^T( y-\hat{y})
%    +\frac{L}{2}\norm{y-\hat{y}}_2^2
%  \end{array}
%  \right\},
%\end{equation}
\begin{equation}
  \label{ch3_eq:sub-problem}
  \min_{(X,s,y)\in\Sk~:~y\in\cQ}
  \left\{
    \begin{array}[c]{l}
      \lambda^{(k)}(\mu_1\norm{\sigma(X)}_\alpha+\mu_2\norm{s}_\beta)+
      \left[
        \begin{array}{c}
            \grad_X f^{(k)}(\tilde{X},\tilde{s},\tilde{y}) \\
            \grad_s f^{(k)}(\tilde{X},\tilde{s},\tilde{y}) \\
            \grad_y f^{(k)}(\tilde{X},\tilde{s},\tilde{y}) \\
        \end{array}
    \right]^T
    \left[
        \begin{array}{c}
            X-\tilde{X} \\
            s-\tilde{s} \\
            y-\tilde{y} \\
        \end{array}
    \right] \\
    +\frac{L}{2}\norm{X-\tilde{X}}_F^2+\frac{L}{2}\norm{s-\tilde{s}}_2^2
    +\frac{L}{2}\norm{y-\tilde{y}}_2^2
    \end{array}
  \right\},
\end{equation}
for a given $(\tilde{X},\tilde{s},\tilde{y})$.
%where
%\eq
%f^{(k)}(X,s,y) :=
%\begin{array}[t]{l}
%  \mbox{} - \lambda^{(k)}(\theta_1^{(k)})^T(\mathcal{A}(X)+y-b) +
%  \frac{1}{2}\norm{\mathcal{A}(X)+y-b}_2^2 \\
%  \mbox{} -
%  \lambda^{(k)}(\theta_2^{(k)})^T(\mathcal{C}(X)+s-d) +
%  \frac{1}{2}\norm{\mathcal{C}(X)+s-d}_2^2
%\end{array}
%\en
%denotes the ``smooth'' part of the objective
%function in \eqref{ch3_eq:augmented_lagrangian_subproblem}.
Note that \eqref{ch3_eq:sub-problem} is \emph{separable} in $X$, $s$
and $y$ variables. Solving \eqref{ch3_eq:sub-problem}
reduces to one ``constrained shrinkage'' in $X\in\reals^{m\times n}$,
see~\eqref{ch3_eq:matrix_constrained_shrinkage}; %~\cite{Ma08_1J}
one ``constrained shrinkage" in $s\in\reals^p$, see
\eqref{ch3_eq:vector_constrained_shrinkage}; %~\cite{Dau04_2J};
and one Euclidean projection onto $\cQ$ in $y\in\reals^q$.

\subsection{Convex set $\cS^{(k)}$ and the initial iterate for k-th
  subproblem}
In the $k$-th FALC iteration, we solve
\eqref{ch3_eq:augmented_lagrangian_subproblem_2} over the convex set
$\cS^{(k)}$ defined in Figure~\ref{ch3_alg:falc}, using
\textbf{Algorithm~APG} starting from the initial iterate $(X^{(k-1)},s^{(k-1)},y^{(k-1)})$.

Let $\eta:=\mu_1\norm{\sigma(X^{(0)})}_\alpha+
\mu_2\norm{\mathcal{C}(X^{(0)})-d}_\beta$ and
$\etak:=\eta
+\frac{\lambda^{(k)}}{2}\left(\norm{\tetki}_2^2+\norm{\tetkii}_2^2\right)$. For
all $k\geq 1$, since
$(\Xkopt,\skopt,\ykopt)\in\argmin_{X\in\reals^{m\times
    n},~s\in\reals^p,~y\in\reals^q}\{P^{(k)}(X,s,y):~y\in\cQ\}$, we
have
\begin{align}
\mu_1\norm{\sigma(\Xkopt)}_\alpha+\mu_2\norm{\skopt}_\beta \leq
\frac{P^{(k)}(\Xkopt,\skopt,\ykopt)}{\lk}\leq
\frac{P^{(k)}(X^{(0)},s^{(0)},y^{(0)})}{\lk}=\etak. \label{ch3_eq:etak1}
\end{align}
Above inequality ensures that $\cS^{(k)} \cap \argmin
\{P^{(k)}(X,s,y):X\in\reals^{m\times n},~s\in\reals^p,~y\in\cQ\} \neq
\emptyset$.

\subsection{\textsc{ITERstop} and \textsc{GRADstop}: Stopping criteria for
  Algorithm APG}
Next, we discuss the stopping criteria set in Line
\ref{ch3_li:inner-stopping-condition1} and Line
\ref{ch3_li:inner-stopping-condition2} of Figure~\ref{ch3_alg:falc}.
\subsubsection{\textsc{ITERstop}}
% Clearly,
% $\etakii$ depends on the properties of $\cQ$ and this is the reason we
% specify the value of $\etakii$ in Lemma~\ref{lem:y-bound} below
% depending on  $\cQ$.
% In the rest of this section, we show that if \textbf{Algorithm~APG}
% stops according to \textsc{ITERstop}, the iterate $(\Xk,\sk,\yk)$
% satisfies
% \begin{equation}
% \label{ch3_eq:modified-inner-a}
%     P^{(k)}(\Xk,\sk,\yk) \leq \min_{X\in\reals^{m\times n},
%     s\in\reals^p, y\in\reals^q}\{P^{(k)}(X,s,y):~y\in\cQ\} +
%     \epsilon^{(k)}.
% \end{equation}

Let
$\{(X_1^{(k,\ell)},s_1^{(k,\ell)},y_1^{(k,\ell)})\}_{\ell\in\integers_+}$
denote the sequence of $x_1$-iterates
% Figure~\ref{ch1_alg:pga}
when \textbf{Algorithm\\~APG} is called %  in
% Line~\ref{algeq:apg_call} of Figure~\ref{ch3_alg:falc}
to solve the $k$-th sub-problem. For the sake of notational simplicity, let
$h^{(k)}(X,s,y):=\norm{X-X^{(k-1)}}_F^2+\norm{s-s^{(k-1)}}_2^2+\norm{y-y^{(k-1)}}_2^2$. Hence,
Lemma~\ref{ch1_lem:tseng_corollary} establishes that
\begin{equation}
P^{(k)}(X_1^{(k,\ell)}, s_1^{(k,\ell)}, y_1^{(k,\ell)})\leq
\inf_{X,s,y}\{P^{(k)}(X,s,y):~y\in\cQ\} + \epsilon^{(k)}\ \mbox{  for
} \ \ell \geq
\sqrt{\frac{2Lh^{(k)}\left(\Xkopt,\skopt,\ykopt\right)}{\epsilon^{(k)}}}-1
\label{ch3_eq:eps_stop_1}
\end{equation}
where $L=\sigma^2_{\max}(M)$ is the Lipschitz constant of $\grad
f^{(k)}$ for all $k\geq 1$. Triangular inequality implies that
\begin{equation}
\sqrt{h^{(k)}\left(\Xkopt,\skopt,\ykopt\right)}\leq\norm{\Xkopt}_F+\norm{X^{(k-1)}}_F
+\norm{\skopt}_2+\norm{s^{(k-1)}}_2
+\norm{\ykopt}_2+\norm{y^{(k-1)}}_2. \label{ch3_eq:eps_stop_2}
\end{equation}
It is easy to show that
\begin{equation}
\label{ch3_eq:IJ_beta}
\frac{1}{I(\alphas)} \norm{\sigma(X)}_{\alpha} \leq \norm{X}_F \leq
I(\alpha) \norm{\sigma(X)}_{\alpha},
\qquad
\frac{1}{J(\betas)} \norm{x}_{\beta} \leq \norm{x}_2 \leq J(\beta)
\norm{x}_{\beta},
\end{equation}
where
\begin{equation}
    \label{ch3_eq:IJ_def}
I(\alpha) =
\left\{
  \begin{array}{ll}
    \sqrt{\min\{m,n\}}, & \alpha = \infty,\\
    1, & \mbox{otherwise},
  \end{array}
\right.
\qquad
J(\beta) =
\left\{
  \begin{array}{ll}
    \sqrt{p}, & \beta = \infty,\\
    1, & \mbox{otherwise},
  \end{array}
\right.
\end{equation}
and $\alphas$ %= \frac{\alpha}{\alpha-1}$
(resp. $\betas$ %=\frac{\beta}{\beta-1}$
) denotes the H{\"o}lder conjugate of the
$\alpha$ (resp. $\beta$), i.e., $\frac{1}{\alpha}+\frac{1}{\alphas}=1$.
In Lemma~\ref{lem:y-bound} in Appendix~\ref{app:proofs} we show that
\eq
\etakii
:=\sigma_{\max}(A)\frac{I(\alpha)}{\mu_1}~\etak+\norm{b+\lk\tetki}_F
+2\norm{\cA(X^{(0)})-b}_2
\en
is an upper bound on $\norm{\ykopt}_2$. Note that when
 $\cQ$ is a bounded set such that $\cQ \subseteq\{y:
\norm{y}_2 \leq \eta_2\}$. Then, one can set $\etakii:=\eta_2$ for all
$k\geq 1$.  Let
\begin{equation}
\label{ch3_eq:apg_iter_cond}
\ell_{\max}^{(k)}:=\sqrt{\frac{2L}{\epsk}}\left[\left(\frac{I(\alpha)}{\mu_1}+
    \frac{J(\beta)}{\mu_2}\right)\etak+\etakii+\norm{X^{(k-1)}}_F+\norm{s^{(k-1)}}_2+
  \norm{y^{(k-1)}}_2\right].
\end{equation}
Then, clearly $\ell_{\max}^{(k)}$ satisfies the following inequality
\begin{equation*}
    \sqrt{\frac{2Lh^{(k)}\left(\Xkopt,\skopt,\ykopt\right)}{\epsilon^{(k)}}}
    \leq \ell_{\max}^{(k)}.
\end{equation*}
Thus, \eqref{ch3_eq:eps_stop_1} implies that when \textbf{Algorithm~APG} terminates due to
$\textsc{ITERstop}(\ell)$, the iterate
$(X_1^{(k,\ell)},s_1^{(k,\ell)},$ $y_1^{(k,\ell)})$ is
$\epsk$-optimal. % and it is assigned to the $k$-th FALC iterate
% $(\Xk,\sk,\yk)$ -see Line~\ref{ch1_algeq:stop} of
% \textbf{Algorithm~APG}.

\subsubsection{\textsc{GRADstop}} % Other than limiting the number of
                                % iterations for the k-th subproblem
                                % \eqref{ch3_eq:augmented_lagrangian_subproblem_2} via \textsc{ITERstop},
The stopping condition \textsc{GRADstop} in Line
\ref{ch3_li:inner-stopping-condition2} of Figure~\ref{ch3_alg:falc} is used to
terminate the \textbf{Algorithm~APG} when a certain set of perturbed first-order
optimality conditions hold at the current iterate. Specifically,
\textbf{Algorithm~APG} stops
according to \textsc{GRADstop}, when we have
\begin{equation}
\label{ch3_eq:modified-inner-b}
\begin{array}{ll}
    (1) & \sqrt{\norm{G}_F^2 + \norm{g}_2^2} \leq \tau^{(k)}, \mbox{ for some } (G,g)\in\partial_{X,s} P^{(k)}(.,.,.)|_{(\Xk,\sk,\yk)}\\
    (2) & \norm{\yk-\Pi_{\cQ}\left(\yk-\frac{1}{L}\grad_y P^{(k)}(\Xk,\sk,\yk)\right)}_2\leq\xik,
\end{array}
\end{equation}
where
% Here $\partial_{X,s} P^{(k)}(.,.,.)|_{(\Xk,\sk,\yk)}$, which we call the set of partial subgradients, denotes the projection of $\partial P^{(k)}$ at $(\Xk,\sk,\yk)$ onto coordinates corresponding to $X$ and $s$. Hence,
\begin{eqnarray}
\label{ch3_eq:partial_subgradient}
\lefteqn{\partial_{X,s} P^{(k)}(.,.,.)|_{(\Xk,\sk,\yk)}}&& \nonumber\\
&=&\left\{
(G,g)\in\reals^{m\times n}\times\reals^p:~
\begin{array}{rl}
G&\in\lk\mu_1~\partial\norm{\sigma(.)}_\alpha|_{\Xk}+\grad_X f^{(k)}(\Xk,\sk,\yk),\\
g&\in\lk\mu_2~\partial\norm{.}_\beta|_{\sk}+\grad_s f^{(k)}(\Xk,\sk,\yk),
\end{array}
\right\}.
\end{eqnarray}
denotes the projection of $\partial P^{(k)}$ at $(\Xk,\sk,\yk)$ onto the
$X$ and $s$ co-ordinates. Note that \eqref{ch3_eq:modified-inner-b}
  would indeed be the first-order optimality conditions if $\tk$ and
  $\xik$ were both set to $0$.

In our numerical experiments, we found that
the calls to \textbf{Algorithm APG} were almost always terminated by
the gradient-based stopping condition \textsc{GRADstop}. This suggests
that relying on only \textsc{ITERstop} to terminate calls to {\bf
  Algorithm APG} is likely to be very inefficient.

% We show in
% Theorem~\ref{ch3_thm:epsilon_convergence} we can guarantee convergence
% to an $\epsilon$-optimal solution even when the solutions to the
% subproblems

\subsection{Convergence Results}
% In the rest of this section, we establish theoretical properties of
% FALC.
Given $\epsilon>0$, let $N_{\rm FALC}(\epsilon)$ be the number of
times \textbf{Algorithm~APG} is called within \textbf{Algorithm~FALC}
until an $\epsilon$-feasible and $\epsilon$-optimal solution to
\eqref{ch3_eq:simple_composite_norm_minimization} is found. During the
$k$-th call, \textbf{Algorithm~APG} inexactly solves the $k$-th
subproblem \eqref{ch3_eq:augmented_lagrangian_subproblem_2}. Let $\Nk$
denote the number of iterations \textbf{Algorithm~APG} needs until one
of the stopping criteria \textsc{ITERstop}, \textsc{GRADstop} or \textsc{FALCstop} is
met, where \textsc{FALCstop} is the stopping condition for \textbf{Algorithm~FALC}. Finally, let $N_{\rm inner}$ be the total number
\textbf{Algorithm~APG} iterations until an $\epsilon$-optimal and
$\epsilon$-feasible solution to
\eqref{ch3_eq:simple_composite_norm_minimization} is computed,
i.e. $N_{\rm inner}=\sum_{k=1}^{N_{\rm FALC}(\epsilon)}\Nk$.
% In Theorem~\ref{ch3_thm:limit-point}, we establish that every limit
% point of the FALC iterate sequence $\{(\Xk,\sk,$ $\yk)\}_{k \in
%   \integers_+}$, is an optimal solution for the composite norm
% minimization problem. However, before going into the details of the
% theorem, we first need to provide a couple of technical results
% below. These results help us prove Lemma~\ref{ch3_lem:bdd-theta},
% which will later be used in Lemma~\ref{lem:bounded} to show that with
% a proper choice of parameter sequence
% $\{\lk,\epsk,\tau^{(k)},\xik\}_{k \in \integers_+}$, we can ensure
% that $\{\tetki\}_{k\in\integers_+}$, $\{\tetkii\}_{k\in\integers_+}$
% and $\{(\Xk,\sk,\yk)\}_{k\in\integers_+}$ sequences are bounded, which
% is an extremely important step to prove
% Theorem~\ref{ch3_thm:limit-point}. First, we give a useful property of
% Euclidean projection that will be used in this section, where we state
% and prove the convergence properties of FALC.

We begin by establishing bounds on the sequence of dual iterates
$\{\tetki\}_{k\in\integers_+}$ and $\{\tetkii\}_{k\in\integers_+}$. In
order to establish this result, we need to bound the infeasibility of
an $\epsk$-optimal solution to the $k$-th sub-problem. In each
iteration of FALC we solve a sub-problem of the form: $\min_{X,s,y}\{P(X,s,y):~y\in\cQ\}$, where
\eq
P(X,s,y)  =  \lambda
(\mu_1\norm{\sigma(X)}_\alpha+\mu_2\norm{s}_\beta)  + \frac{1}{2}
\norm{\mathcal{A}(X)-y-b-\lambda\theta_1}_2^2+\frac{1}{2}
\norm{\mathcal{C}(X)-s-d-\lambda\theta_2}_2^2.
\en
Suppose
$(\bar{X},\bar{s},\bar{y})$ is $\epsilon$-optimal, %for the problem $\min_{X,s,y}\{P(X,s,y):~y\in\cQ\}$,
i.e., $ 0\leq P(\bar{X},\bar{s},\bar{y})- \min_{X,s,y}\{P(X,s,y):~y\in\cQ\}%\min_{X \in \reals^{m\times n},~s\in\reals^p,~y\in\cQ\subset\reals^q}P(X,s,y)
\leq \epsilon$.
In Lemma~\ref{ch3_cor:grad_norm_bound} in Appendix~\ref{app:proofs} we
establish that
\begin{equation}
  \label{eq:infeas-bnd}
  \begin{array}{rcl}
  \norm{\mathcal{C}(\bar{X})-\bar{s}-d-\lambda\theta_2}_2 & \leq &
J(\betas)\mu_2\lambda+\sigma_{\max}(M)\sqrt{2\epsilon},\\
\norm{\mathcal{A}^*\left(\mathcal{A}(\bar{X})-\bar{y}-b-\lambda\theta_1\right)+\mathcal{C}^*\left(\mathcal{C}(\bar{X})-\bar{s}-d-\lambda\theta_2\right)}_F
&\leq & I(\alphas)\mu_1\lambda+ \sigma_{\max}(M)\sqrt{2\epsilon},
\end{array}
\end{equation}
With this bound, we are now ready to show that the dual iterates are
bounded.
\begin{lemma}
    \label{ch3_lem:bdd-theta}
    For all  $k>1$, the elements of $\{\tetki\}_{k\in\integers_+}$ and $\{\tetkii\}_{k\in\integers_+}$ satisfy the following relation
\begin{subequations}
\label{ch3_eq:theta_induction}
\begin{align}
&\norm{\tetkii}_2\leq
\max\left\{\sigma_{\max}(M)\sqrt{\frac{2\epsilon^{(k-1)}}{(\lambda^{(k-1)})^2}},
  \
  \frac{\tau^{(k-1)}}{\lambda^{(k-1)}}\right\}+J(\betas)\mu_2, \label{ch3_eq:theta2_induction}\\
&\norm{\tetki}_2 \leq \frac{1}{\sigma_{\min}(A)}\left[\sigma_{\max}(C)\norm{\tetkii}_2 +
\max\left\{\sigma_{\max}(M)\sqrt{\frac{2\epsilon^{(k-1)}}{(\lambda^{(k-1)})^2}},
  \
  \frac{\tau^{(k-1)}}{\lambda^{(k-1)}}\right\}+I(\alphas)\mu_1\right].\label{ch3_eq:theta1_induction}
%&\norm{\mathcal{A}^*(\tetki)}_F \leq \norm{\mathcal{C}^*(\tetkii)}_F +
%\max\left\{\sigma_{\max}(M)\sqrt{\frac{2\epsilon^{(k-1)}}{(\lambda^{(k-1)})^2}},
%  \
%  \frac{\tau^{(k-1)}}{\lambda^{(k-1)}}\right\}+I(\alphas)\mu_1.\label{ch3_eq:theta1_induction}
\end{align}
\end{subequations}
\end{lemma}
\begin{proof}
Consider the following two cases:
\begin{enumerate}[(a)]
\item \emph{The $k$-th call to \textbf{Algorithm~APG} terminates with \textsc{ITERstop}}:
%the iterate $(\Xk,\sk,\yk)$ satisfying
%\eqref{ch3_eq:modified-inner}(a).
Since the iterate is $\epsk$-optimal, the bound~\eqref{eq:infeas-bnd}
implies that
%Lemma~\ref{ch3_cor:grad_norm_bound} in Appendix~\ref{app:proofs} guarantees that
\begin{subequations}
  \label{ch3_eq:subgradient_feasibility_bound_a}
  \begin{align}
    & \norm{\mathcal{C}(\Xk)-\sk-d-\lk\tetkii}_2    \leq  J(\betas)\mu_2\lk+ \sigma_{\max}(M)\sqrt{2\epsilon^{(k)}}\
    ,\label{ch3_eq:subgradient_feasibility_bound_a2}\\
    & \norm{\mathcal{A}^*(\mathcal{A}(\Xk)-\yk-b-\lk\tetki)+\mathcal{C}^*(\mathcal{C}(\Xk)-\sk-d-\lk\tetkii)}_F   \nonumber\\
      &\mbox{} \leq  I(\alphas)\mu_1\lk+\sigma_{\max}(M)\sqrt{2\epsilon^{(k)}}. \label{ch3_eq:subgradient_feasibility_bound_a1}
  \end{align}
\end{subequations}
\item  \emph{The $k$-th call to \textbf{Algorithm~APG} terminates with \textsc{GRADstop}}:
%an iterate $(\Xk,\sk,\yk)$ that satisfies \eqref{ch3_eq:modified-inner}(b).
In this case, there exists
  $Q^{(k)}\in\partial \norm{\sigma(.)}_\alpha|_{\Xk}$ and
  $q^{(k)}\in\partial \norm{.}_\beta|_{\sk}$ such that
  \begin{align*}
    \sqrt{\norm{\lambda^{(k)}\mu_1 Q^{(k)} + \grad_X
        f^{(k)}(\Xk,\sk,\yk)}_F^2+\norm{\lambda^{(k)}\mu_2 q^{(k)} +
        \grad_s f^{(k)}(\Xk,\sk,\yk)}_2^2} ~\leq~ \tau^{(k)}.
  \end{align*}
Since $\norm{q^{(k)}}_{\betas} \leq 1$ and $\norm{\sigma(Q^{(k)})}_{\alphas} \leq 1$, from the definition of $I(.)$ and $J(.)$ in \eqref{ch3_eq:IJ_beta}, it follows that
$\norm{\sigma(Q^{(k)})}_F\leq I(\alphas)$ and $\norm{q^{(k)}}_2\leq J(\betas)$. Then we have
  \begin{subequations}
    \label{ch3_eq:subgradient_feasibility_bound_b}
    \begin{align}
      &\norm{\mathcal{C}(\Xk)-\sk-d-\lk\tetkii}_2 \leq J(\betas)\mu_2\lambda^{(k)}+\tau^{(k)},\label{ch3_eq:subgradient_feasibility_bound_b2}\\
      &\norm{\mathcal{A}^*(\mathcal{A}(\Xk)-\yk-b-\lk\tetki)
        +\mathcal{C}^*(\mathcal{C}(\Xk)-\sk-d-\lk\tetkii)}_F
      \nonumber\\
      &\leq~
      I(\alphas)\mu_1\lambda^{(k)}+\tau^{(k)}. \label{ch3_eq:subgradient_feasibility_bound_b1}
    \end{align}
  \end{subequations}
\end{enumerate}
Thus, combining \eqref{ch3_eq:subgradient_feasibility_bound_a} and
\eqref{ch3_eq:subgradient_feasibility_bound_b}, and using triangular
inequality it follows that for all $k\geq 1$
\begin{subequations}
  \label{ch3_eq:feasibility_bound}
  \begin{align}
    \norm{\mathcal{C}(\Xk)-\sk-d-\lk\tetkii}_2 \leq
    & ~J(\betas)\mu_2\lambda^{(k)}+\max\Big\{\sigma_{\max}(M)\sqrt{2\epsilon^{(k)}},\tau^{(k)}\Big\},\label{ch3_eq:feasibility_bound_1}\\
    \norm{\mathcal{A}^*(\mathcal{A}(\Xk)-\yk-b-\lk\tetki)}_F \leq
    &
    ~I(\alphas)\mu_1\lambda^{(k)}+\max\Big\{\sigma_{\max}(M)\sqrt{2\epsilon^{(k)}},\tau^{(k)}\Big\}
    \nonumber\\
    & \mbox{} + \norm{\mathcal{C}^*(\mathcal{C}(\Xk)-\sk-d-\lk\tetkii)}_F. \label{ch3_eq:feasibility_bound_2}
  \end{align}
\end{subequations}
Since
$\theta_1^{(k+1)}=\theta_1^{(k)}-\frac{\mathcal{A}(\Xk)-\yk-b}{\lambda^{(k)}}$
and
$\theta_2^{(k+1)}=\theta_2^{(k)}-\frac{\mathcal{C}(\Xk)-\sk-d}{\lambda^{(k)}}$,
\eqref{ch3_eq:feasibility_bound} is obtained by dividing
\eqref{ch3_eq:theta_induction} into $\lk$ and using the fact that
Assumption~\ref{ass:injective} implies that $A$ has full row rank,
i.e. $\sigma_{\min}(A)>0$. Thus, $\big\{(\tetki,\tetkii)\big\}_{k\in
  \integers}$ satisfies \eqref{ch3_eq:theta_induction}.
\end{proof}

\noindent Next, we establish that the FALC iterate sequence
$\{(\Xk,\sk,\yk)\}_{k\in\integers_+}$ is  bounded.
\begin{lemma}
\label{lem:bounded}
Let $(\Xkopt,\skopt,\ykopt)$ be an optimal solution to
\eqref{ch3_eq:augmented_lagrangian_subproblem_2} and let
$\{(\Xk,\sk,\yk)\}_{k\in\integers_+}$ denote the sequence of FALC iterates
% generated by \textbf{Algorithm~FALC} displayed in
% Figure~\ref{ch3_alg:falc} for
%   a fixed sequence of
corresponding to a parameter sequence $\{(\lk,\epsk,\tau^{(k)},\xik)\}_{k
  \in \integers_+}$ such that
  \begin{enumerate}[(i)]
    \item penalty multiplier, $\lambda^{(k)} \searrow 0$,
    \item approximate optimality parameter, $\epsilon^{(k)}\searrow 0$
      such that $\frac{\epsilon^{(k)}}{(\lambda^{(k)})^2}\leq B$ for all
      $k\geq1$ for some $B>0$,
    \item subgradient tolerance parameters, $\tk\searrow 0$ and $\xik
      \searrow 0$ such that $\frac{\tk}{\lk}\rightarrow 0$ and
      $\frac{\xik}{\lk}\rightarrow 0$ as $k\rightarrow\infty$.
  \end{enumerate}
  Then $\{(\Xkopt,\skopt,\ykopt)\}_{k\in\integers_+}$ and
  $\{(\Xk,\sk,\yk)\}_{k\in\integers_+}$ are bounded sequences.
\end{lemma}
\begin{proof}
In the $k$-th FALC iteration, the call to \textbf{Algorithm~APG}
% displayed in Figure~\ref{ch1_alg:pga}
terminates in at most  $\ell^{(k)}_{\max}$ iterations. Since
$\ell^{(k)}_{\max}$ is finite for all $k \geq 1$, the sequence
$\{(\Xk,\sk,\yk)\}_{k\in\integers_+}$ exists. In order to show
$\{(\Xk,\sk,\yk)\}_{k\in\integers_+}$ is a bounded sequence, we first
establish that $\{\tetki\}_{k\in\integers_+}$ and
$\{\tetkii\}_{k\in\integers_+}$ are bounded sequences.

Because %$A$ has full row-rank,
$\frac{\epsilon^{(k)}}{(\lambda^{(k)})^2}\leq B$ and
$\frac{\tau^{(k)}}{\lambda^{(k)}}\rightarrow 0$,
\eqref{ch3_eq:theta_induction} implies that there exist constants $B_{\theta_1}>0$ and
$B_{\theta_2}>0$ such that
\begin{equation}
\label{ch3_eq:theta_bound}
\max_{k \geq 1} \{\norm{\tetki}_2\}\leq B_{\theta_1},\ \mbox{and}\ \max_{k
  \geq 1}\{\norm{\tetkii}_2\}\leq B_{\theta_2}.
\end{equation}
From~\eqref{ch3_eq:theta_bound}, it follows that for $i = 1,2$,
\begin{equation}
  \lim_{k\rightarrow\infty}\lambda^{(k)}\theta_i^{(k)}=0,\ \mbox{ and }\
  \lim_{k\rightarrow\infty}\lambda^{(k)}\norm{\theta_i^{(k)}}_2^2=0.
  \label{ch3_eq:lambda_theta_limit}
\end{equation}
Also,  $\frac{\epsilon^{(k)}}{(\lambda^{(k)})^2}\leq B$ for all $k\geq 1$
implies that $\lim_{k\rightarrow\infty}\frac{\epsilon^{(k)}}{\lambda^{(k)}}=0$.
%\begin{align}
%\lim_{k\rightarrow\infty}\frac{\epsilon^{(k)}}{\lambda^{(k)}}=0.
%\label{ch3_eq:epsilon_lambda_division}
%\end{align}

Now we can prove that the iterate sequence is bounded. Trivially, %for all
% $k\geq 1$, $\norm{\yk}_\gamma\leq\delta$
% and
the choice of $\Sk$ ensures that
$\mu_1\max\{\norm{\sigma(\Xk)}_\alpha,~\norm{\sigma(\Xkopt)}_\alpha\}\leq\etak$
and $\mu_2\max\{\norm{\sk}_\beta, \norm{\skopt}_\beta\}\leq\etak$. From the definition of $\etak$ in
Line~\ref{ch3_li:etak} of Figure~\ref{ch3_alg:falc} and
\eqref{ch3_eq:theta_bound}, it follows that for all $k\geq 1$
\begin{equation}
\label{ch3_eq:B_eta}
\etak \leq \eta+\lk\left(\frac{B^2_{\theta_1}+B^2_{\theta_2}}{2}\right)
\leq
\eta+\lambda^{(1)}\left(\frac{B^2_{\theta_1}+B^2_{\theta_2}}{2}\right):=
B_{\eta_1}.
\end{equation}
Hence, for all $k\geq 1$,
\begin{equation}
    \label{ch3_eq:falc_iter_bound}
    \mu_1\max\{\norm{\sigma(\Xk)}_\alpha,~\norm{\sigma(\Xkopt)}_\alpha\}\leq
    B_{\eta_1} \mbox{ and }\ \mu_2\max\{\norm{\sk}_\beta, \norm{\skopt}_\beta\}\leq B_{\eta_1}.
\end{equation}
Next we show that $\{\yk\}_{k\in\integers_+}$ and
$\{\ykopt\}_{k\in\integers_+}$ are bounded. From the definition of
$\theta_1^{(k+1)}$ in Line~\ref{ch3_li:lagrangian_update_1} of
Figure~\ref{ch3_alg:falc}, we have
$\yk=\lk(\theta_1^{(k+1)}-\tetki)+\cA(\Xk)-b$ for all $k\geq
1$. Hence,
$\norm{\yk}_2\leq\lk\norm{(\theta_1^{(k+1)}-\tetki)}_2+\norm{\cA(\Xk)-b}_2$
for all $k\geq 1$. From \eqref{ch3_eq:theta_bound} and
\eqref{ch3_eq:falc_iter_bound}, it follows that there exists $B_y>0$
such that
\begin{align}
\label{eq:y-bound-uniform}
\norm{\yk}_2\leq B_y, \quad \forall~k\geq 1.
\end{align}
Moreover, Lemma~\ref{lem:y-bound} in Appendix~\ref{app:proofs} guarantees
that for all $k\geq 1$,
$\norm{\ykopt}_2\leq\etakii$, where $\etakii$ is given in
\eqref{eq:eta2} in Appendix~\ref{app:proofs}. Since $\etak\leq B_{\eta_1}$ for all $k\geq 1$,
\eqref{ch3_eq:lambda_theta_limit} and \eqref{eq:eta2} imply that there
exists a constant $B_{\eta_2}>0$ such that
\begin{align}
\label{eq:yopt-bound-uniform}
\norm{\ykopt}_2\leq B_{\eta_2}, \quad \forall~k\geq 1.
\end{align}
\end{proof}
\begin{theorem}
  \label{ch3_thm:limit-point}
  Let $\mathcal{X}=\{\Xk\}_{k \in \integers_+}$ denote the FALC
  iterate sequence corresponding to a parameter sequence
  $\{(\lk,\epsk,\tau^{(k)},\xik)\}_{k \in \integers}$ satisfying the
  conditions in Lemma~\ref{lem:bounded}.
  Then any limit point $\bar{X}$ of the sequence $\mathcal{X}$ is an
  optimal solution of the composite norm
  minimization
  problem~(\ref{ch3_eq:simple_composite_norm_minimization}).
\end{theorem}
\begin{proof}
Since Lemma~\ref{lem:bounded} guarantees that $\cX$ is a bounded
sequence, there exists a subsequence
$\cK\subset\integers_+$ such that
  $\lim_{k\in\cK}X^{(k)}=\bar{X}$ exists. We have previously shown
  that $\lim_{k\rightarrow\infty}\lambda^{(k)}\theta_i^{(k)}=0$ for
  $i\in\{1,2\}$. Hence, \eqref{ch3_eq:feasibility_bound_1} guarantees that
  $\lim_{k\in\cK}s^{(k)}=\bar{s}$ exists; similarly
  \eqref{ch3_eq:feasibility_bound_2} and the full row-rank assumption on
  $A$ together guarantee that
  $\lim_{k\in\cK}y^{(k)}=\bar{y}$ exists. Then, taking the limit of both sides of
  \eqref{ch3_eq:feasibility_bound_1} for $k\in\cK$, we have
  $\norm{\mathcal{C}(\bar{X})-\bar{s}-d}_2\leq 0$,
  i.e. $\bar{s}=\mathcal{C}(\bar{X})-d$. Moreover, taking the limit of
  both sides of \eqref{ch3_eq:feasibility_bound_2} for $k\in\cK$
  and using the fact that $\bar{s}=\mathcal{C}(\bar{X})-d$, we have
  $\norm{\mathcal{A}^*(\mathcal{A}(\bar{X})-\bar{y}-b)}_2\leq 0$. Since
  $A$ has full row rank and $\yk\in\cQ$ for all $k\geq 1$, we have
  \begin{align}
    \mathcal{A}(\bar{X})-b=\bar{y}, \qquad
    \bar{y}\in\cQ. \label{ch3_eq:feasibilty}
  \end{align}
Therefore, we can conclude that $\bar{X}$ is feasible, i.e. $\mathcal{A}(\bar{X})-b\in\cQ$.

In the rest of the proof, we will show that
$\bar{X}\in\argmin\{\mu_1\norm{\sigma(X)}_\alpha+\mu_2\norm{\mathcal{C}(X)-d}_\beta:\
\mathcal{A}(X)-b\in\cQ\}$. We consider the following two cases:
\begin{enumerate}[(a)]
\item There exists a subsequence $\cK_a\subset\cK$ such
that for all $k\in\cK_a$, \textbf{Algorithm~APG} terminates with
\textsc{ITERstop}; hence, % $(\Xk,\sk,\yk)$ satisfies
% \eqref{ch3_eq:modified-inner-a}, i.e. the sequence $\{
the iterate $(\Xk, \sk,
\yk)$ computed in
Step~\ref{algeq:apg_call} of FALC satisfies
\begin{align}
0\leq P^{(k)}(\Xk, \sk, \yk)-P^{(k)}(\Xkopt, \skopt, \ykopt)\leq
\epsilon^{(k)} \hspace{5mm} \forall
k \in \cK_a. \label{ch3_eq:inexact_minimizer_relation}
\end{align}
Fix $\Xopt\in\argmin_{X\in\reals^{m\times
    n}}\{\mu_1\norm{\sigma(X)}_\alpha
+\mu_2\norm{\mathcal{C}(X)-d}_\beta:~\mathcal{A}(X)-b\in\cQ\}$,
let $\sopt := \mathcal{C}(\Xopt)-d$ and $\yopt := \mathcal{A}(\Xopt)-b$. Since
$(\Xkopt,\skopt,\ykopt)\in\argmin_{X\in\reals^{m\times
    n},~s\in\reals^p,~y\in\reals^q}\{P^{(k)}(X,s,y):~y\in\cQ\}$, it
follows that $P^{(k)}(\Xkopt, \skopt,\ykopt)\leq P^{(k)}(\Xopt, \sopt,
\yopt)$ for $k\geq 1$. Thus, \eqref{ch3_eq:inexact_minimizer_relation} implies that
$P^{(k)}(\Xk,\sk, \yk)\leq
P^{(k)}(\Xopt,\sopt,\yopt)+\epsilon^{(k)}$. Hence, for all $k\in \cK_a$,
\begin{align}
\mu_1\norm{\sigma(\Xk)}_\alpha+\mu_2\norm{\sk}_\beta&\leq \frac{P^{(k)}(\Xk,\sk,\yk)}{\lambda^{(k)}}
%,\nonumber\\&\leq
~\leq\frac{P^{(k)}(\Xopt,\sopt,\yopt)+\epsilon^{(k)}}{\lambda^{(k)}},\nonumber\\
&= \mu_1\norm{\sigma(\Xopt)}_\alpha+\mu_2\norm{\mathcal{C}(\Xopt)-d}_\beta+
\frac{\lambda^{(k)}}{2}\left(\norm{\theta_1^{(k)}}_2^2 + \norm{\theta_2^{(k)}}_2^2\right)
+ \frac{\epsilon^{(k)}}{\lambda^{(k)}}. \label{ch3_eq:inexact_minimizer_bound}
\end{align}
Taking the limit of both sides of \eqref{ch3_eq:inexact_minimizer_bound} along
the subsequence $\cK_a$, and using the fact that
$\bar{s}=\mathcal{C}(\bar{X})-d$, we get
\begin{eqnarray}
\mu_1\norm{\sigma(\bar{X})}_\alpha+\mu_2\norm{\mathcal{C}(\bar{X})-d}_\beta
&= &
\lim_{k\in\cK_a}\mu_1\norm{\sigma(X^{(k)})}_\alpha+\mu_2\norm{\sk}_\beta,
\nonumber \\
&\leq &
\mu_1\norm{\sigma(\Xopt)}_\alpha+\mu_2\norm{\mathcal{C}(\Xopt)-d}_\beta\nonumber\\
&&+\lim_{k\in\cK_a}\left\{\frac{\lambda^{(k)}}{2}\left(\norm{\theta_1^{(k)}}_2^2
    + \norm{\theta_2^{(k)}}_2^2\right)+\frac{\epsilon^{(k)}}{\lambda^{(k)}}\right\}, \nonumber \\
&= &
\mu_1\norm{\sigma(\Xopt)}_\alpha+\mu_2\norm{\mathcal{C}(\Xopt)-d}_\beta,
\label{ch3_eq:inexact_minimizer_optimality}
\end{eqnarray}
where \eqref{ch3_eq:inexact_minimizer_optimality} follows from the fact that
$\{\theta^{(k)}_i\}$ is uniformly bounded for $i\in\{1,2\}$, $\lambda^{(k)}
\rightarrow 0$, and $\epsilon^{(k)}/\lambda^{(k)} \rightarrow 0$. Thus, from
\eqref{ch3_eq:feasibilty}, \eqref{ch3_eq:inexact_minimizer_optimality} and the fact that $\Xopt$ is optimal,
%$\Xopt \in \argmin\{\mu_1\norm{\sigma(X)}_\alpha$$+\mu_2\norm{\mathcal{C}(\Xopt)-d}_\beta:~\mathcal{A}(X)-b\in\cQ\}$,
it follows that $\bar{X}$ is also an optimal solution for the composite norm minimization problem
\eqref{ch3_eq:simple_composite_norm_minimization}.

\item There exists $K\in\cK$ such that, for all $k\in\cK_b := \cK \cap \{k \geq K\}$, \textbf{Algorithm~APG} terminates with \textsc{GRADstop}; hence, $(\Xk,\sk, \yk)$ satisfies \eqref{ch3_eq:modified-inner-b}.

  For all $k \in \cK_b$, there exist $Q^{(k)}\in\partial
  \norm{\sigma(.)}_\alpha|_{\Xk}$ and $q^{(k)}\in\partial
  \norm{.}_\beta|_{\sk}$ such that \eqref{ch3_eq:modified-inner-b}
  holds. Hence, we have
  \begin{subequations}
    \label{ch3_eq:detailed_subgradient_condition}
    \begin{align}
      &\norm{\lambda^{(k)}\mu_2 q^{(k)}+\grad_s f^{(k)}(\Xk,\sk,\yk)}_2
      \leq \tau^{(k)},  \label{ch3_eq:detailed_subgradient_condition_s}\\
      &\norm{\lambda^{(k)}\mu_1 Q^{(k)}+\grad_X f^{(k)}(\Xk,\sk,\yk)}_F
      \leq \tau^{(k)},  \label{ch3_eq:detailed_subgradient_condition_X}\\
      &\norm{\yk-\Pi_{\cQ}\left(\yk-\frac{1}{L}\grad_y f^{(k)}(\Xk,\sk,\yk)\right)}_2\leq\xik. \label{ch3_eq:detailed_subgradient_condition_y}
    \end{align}
  \end{subequations}
  For all $k\in \cK_b$,
  $Q^{(k)}\in\partial\norm{\sigma(.)}_\alpha|_{X^{(k)}}$ and
  $q^{(k)}\in\partial\norm{.}_\beta|_{s^{(k)}}$, therefore,
  $\norm{\sigma(Q^{(k)})}_{\alphas} \leq 1$ and $\norm{q^{(k)}}_{\betas}
  \leq 1$. Hence, there exists a subsequence $\cK'_b\subset\cK_b$ such that
  $\lim_{k\in\cK'_b}(Q^{(k)},q^{(k)})=(\bar{Q},\bar{q})$ exists. One can
  easily show that $\bar{Q}\in\partial\norm{\sigma(.)}_\alpha|_{\bar{X}}$
  and $\bar{q}\in\partial\norm{.}_\beta|_{\bar{s}}$. Dividing both sides
  of \eqref{ch3_eq:detailed_subgradient_condition_s} by $\lambda^{(k)}$, we
  get
  \begin{align}
    \norm{\mu_2 q^{(k)}+\theta_2^{(k+1)}}_2 \leq
    \frac{\tau^{(k)}}{\lambda^{(k)}}, \label{ch3_eq:subgradient_theta2_convergence}
  \end{align}
  for all $k\in\cK_b\supset\cK'_b$. Since $\lim_{k\in\cK'_b}q^{(k)}=\bar{q}$
  and $\lim_{k\in\integers_+}\frac{\tau^{(k)}}{\lambda^{(k)}}=0$, it follows
  that $\lim_{k\in\cK'_b}\theta_2^{(k+1)}=\bar{\theta}_2$ exists and taking
  the limit of both sides of \eqref{ch3_eq:subgradient_theta2_convergence} along $k\in\cK'_b$, we
  have
  \begin{align*}
    \bar{\theta}_2=-\mu_2\bar{q}.
  \end{align*}
  Dividing both sides of \eqref{ch3_eq:detailed_subgradient_condition_X}
  by $\lambda^{(k)}$, we get

  \begin{align}
    \norm{\mu_1
      Q^{(k)}-\mathcal{A}^*(\theta_1^{(k+1)})-\mathcal{C}^*(\theta_2^{(k+1)})}_F
    \leq\frac{\tau^{(k)}}{\lambda^{(k)}}, \label{ch3_eq:subgradient_theta1_convergence}
  \end{align}
  for all $k\in\cK_b\supset\cK'_b$. Since $\lim_{k\in\cK'_b}Q^{(k)}=\bar{Q}$,
  $\lim_{k\in\integers_+}\frac{\tau^{(k)}}{\lambda^{(k)}}=0$ and $A$ has
  full row rank, it follows that
  $\lim_{k\in\cK'_b}\theta_1^{(k+1)}=\bar{\theta}_1$ exists and taking the
  limit of both sides of \eqref{ch3_eq:subgradient_theta1_convergence} along $k\in\cK'_b$, we have
  $\mu_1\bar{Q}+\mu_2\mathcal{C}^*(\bar{q})=\mathcal{A}^*(\bar{\theta}_1)$. Note
  that $\bar{q}\in\partial \norm{.}_\beta|_{\bar{s}}$ and
  $\bar{s}=\mathcal{C}(\bar{X})-d$. Hence,
  $\mathcal{C}^*(\bar{q})\in \partial
  \norm{\mathcal{C}(.)-d}_\beta|_{\bar{X}}$ and we have
  \begin{align}
    \label{ch3_eq:subgradient_optimality_1}
    \mathcal{A}^*(\bar{\theta}_1) = G_* \quad \hbox{and} \quad
    G_*\in\partial\mu_1\norm{\sigma(.)}_\alpha +\mu_2\norm{\mathcal{C}(.)-d}_\beta|_{\bar{X}},
  \end{align}
  where $G_*:=\mu_1\bar{Q}+\mu_2\mathcal{C}^*(\bar{q})$.

  Let $\yk_p:=\yk-\frac{1}{L}\grad_y f^{(k)}(\Xk,\sk,\yk)$ denote the
  gradient step. Since $\xik\searrow 0$, taking the limit of both
  sides of \eqref{ch3_eq:detailed_subgradient_condition_y} along
  $k\in\cK'_b$, we get
  \begin{align}
  \label{eq:y_limit}
  \bar{y}=\lim_{k\in\cK'_b}\yk=\lim_{k\in\cK'_b}\Pi_\cQ(\yk_p)=\Pi_\cQ(\lim_{k\in\cK'_b}\yk_p),
  \end{align}
  where the third equality follows from the fact that Euclidean
  projection $\Pi_\cQ(.)$ is continuous when $\cQ$ is a
  nonempty, closed, convex set; and $\lim_{k\in\cK'_b}\yk_p$ exists
  since $\grad f^{(k)}$ is Lipschitz continuous.

  Dividing both sides of \eqref{ch3_eq:detailed_subgradient_condition_y}
  by $\lk$ and taking the limit along $k\in\cK'_b$, we get
  \begin{align}
    \lim_{k\in\cK'_b}\norm{\yk/\lk-\Pi_{\cQ/\lk}(\yk_p/\lk)}_2=\lim_{k\in\cK'_b}
    \norm{\yk/\lk-\Pi_\cQ(\yk_p)/\lk}_2=0,
    \label{ch3_eq:subgradient_theta1_convergence_2}
  \end{align}
  where the first equality follows from Lemma~\ref{lem:euclidean_proj}.

  For all $k\in\cK_b$, \eqref{ch3_eq:subgradient_theta1_convergence_3}
  follows from the definition of Euclidean projection:
  \begin{align}
  \fprod{\Pi_{\cQ/\lk}(\yk_p/\lk)-\yk_p/\lk,~y/\lk-\Pi_{\cQ/\lk}(\yk_p/\lk)}\geq
  0, \quad \forall y\in\cQ \label{ch3_eq:subgradient_theta1_convergence_3}
  \end{align}
  Since $\yk_p/\lk=\yk/\lk-\theta_1^{(k+1)}/L$, multiplying the second
  term of the inner product in
  \eqref{ch3_eq:subgradient_theta1_convergence_3} by $\lk$ and using
  Lemma~\ref{lem:euclidean_proj}, it follows that for all $k\in\cK_b$
  \begin{align}
  \fprod{\Pi_{\cQ/\lk}(\yk_p/\lk)-\yk/\lk+\theta_1^{(k+1)}/L,~y-\Pi_\cQ(\yk_p)}\geq
  0, \quad \forall y\in\cQ \label{ch3_eq:subgradient_theta1_convergence_4}
  \end{align}
  Since
  $\lim_{k\in\cK'_b}\theta_1^{(k+1)}=\bar{\theta}_1$, taking the limit of
  both sides of \eqref{ch3_eq:subgradient_theta1_convergence_4} along
  $k\in\cK'_b\subset\cK_b$ and
  using \eqref{ch3_eq:subgradient_theta1_convergence_2}, we have
  \begin{align*}
    \fprod{\bar{\theta}_1,~y-\Pi_\cQ(\lim_{k\in\cK'_b}\yk_p)}\geq 0, \quad \forall y\in\cQ.
  \end{align*}
  Thus, it follows from \eqref{eq:y_limit} and above inequality that
  \begin{align}
    \fprod{\bar{\theta}_1,~y-\bar{y}}\geq 0 \qquad \forall
    y\in\cQ. \label{ch3_eq:subgradient_optimality_2}
  \end{align}
  Consequently, \eqref{ch3_eq:subgradient_optimality_1} and
  \eqref{ch3_eq:subgradient_optimality_2} together imply that
  $(\bar{X},\bar{y})$ satisfies the first order optimality conditions of
  the relaxed problem~\eqref{ch3_eq:relaxed_problem}.
  \begin{align}
    \label{ch3_eq:relaxed_problem}
    \min_{X\in\reals^{m\times n},y\in\reals^q}
    \Big\{\mu_1\norm{\sigma(X)}_\alpha
    +\mu_2\norm{\mathcal{C}(X)-d}_\beta-
    (\bar{\theta}_1)^T(\mathcal{A}(X)-y-b):\
    y\in\cQ\Big\}.
  \end{align}
  Since \eqref{ch3_eq:relaxed_problem} is convex, it follows that
  $(\bar{X},\bar{y})$ is an optimal solution to the relaxed
  problem~\eqref{ch3_eq:relaxed_problem}. Moreover, from
  \eqref{ch3_eq:feasibilty}, $(\bar{X},\bar{y})$ is feasible to the
  composite norm minimization problem,
  i.e. $\min\{\mu_1\norm{X}_\alpha+\mu_2\norm{\mathcal{C}(X)-d}_\beta:\
  \mathcal{A}(X)-y=b, y\in\cQ\}$. Therefore,
  $\bar{X}\in\argmin\{\mu_1\norm{X}_\alpha+\mu_2\norm{\mathcal{C}(X)-d}_\beta:\ \mathcal{A}(X)-b\in\cQ\}$.
\end{enumerate}
\end{proof}
Clearly, when \eqref{ch3_eq:simple_composite_norm_minimization} has a
unique solution, the FALC iterates converge to this unique solution.
\begin{corollary}
  \label{ch3_cor:unique}
  Suppose the composite norm minimization problem
  \eqref{ch3_eq:simple_composite_norm_minimization} has a {\em unique} optimal
  solution $X_*$. Let $\{X^{(k)}\}_{ k \in \integers_+}$ denote the sequence
  of FALC iterates % generated by \textbf{Algorithm~FALC} displayed in
                   % Figure~\ref{ch3_alg:falc} when
  corresponding to a parameter sequence $\{(\lk,\epsk,\tk,\xik)\}_{k \in \integers}$ satisfying
  the conditions in Lemma~\ref{lem:bounded}. %Theorem~\ref{ch3_thm:limit-point}.
  Then $\lim_{k\rightarrow\infty}X^{(k)}=X_*$.
  %where $X_*=\argmin_{X\in\reals^{m\times n}}\{\mu_1\norm{\sigma(X)}_\alpha+\mu_2\norm{\mathcal{C}(X)-d}_\beta: \norm{\mathcal{A}(X)-b}_\gamma\leq\delta\}$.
\end{corollary}

For most sparse optimization problems such as basis pursuit and the
matrix completion problems, the unknown signal can be
recovered only if the corresponding convex relaxation has a unique
solution.
% For compressed sensing, exact recovery occurs only if
% $\min_{x\in\reals^n}\{\norm{x}_1: Ax=b\}$ has a {\em unique}
% solution. Similarly, for matrix completion problems, one can recover
% the data matrix exactly only if $\min_{X\in\reals^{m\times
%     n}}\{\|X\|_* : X_{ij}=M_{ij}\ \ (i,j)\in\Omega\}$ has a {\em
%   unique} solution. That is to say uniqueness is a \emph{necessary}
% condition for the exact recovery.
Additionally,  when the set of
constraints for the basis pursuit or affine rank minimization problems
are defined by randomly generated Gaussian matrices, the unknown
target signal is, with very
high probability,  the unique solution to these
optimization problems.
%\red{In the rest of this section we will assume that the composite norm minimization
%problem~\eqref{ch3_eq:simple_composite_norm_minimization} has a unique solution $X_*$.}

We next establish a bound on the iteration complexity of FALC. % computing the
% iterate $(\Xk,\sk,\yk)$ when $\{(\lk,\epsk,$
% $\tk,\xik)\}_{k\in\integers_+}$ satisfies all the conditions in
% Theorem~\ref{ch3_thm:limit-point}. This result is presented in
In Lemma~\ref{ch3_lem:tseng} we prove a uniform bound on the
number of \textbf{Algorithm~APG} iterations required to inexactly solve
any subproblem encountered in FALC. % to obtain the next FALC iterate
\begin{lemma}
  \label{ch3_lem:tseng}
 Suppose  the parameter sequence
  $\{(\lk,\epsk,\tk,\xik)\}_{k\in\integers_+}$ satisfies all the
  conditions in Lemma~\ref{lem:bounded}. %Theorem~\ref{ch3_thm:limit-point}.
  Then there exists constant $\mathcal{N}$ such that, for all $k\geq 1$, the number of
  iterations $\Nk$ required by \textbf{Algorithm~APG}
  to compute $(\Xk,\sk,\yk)$ satisfies
  \begin{align}
    \label{ch3_eq:apg_iter_bound}
    \Nk\leq \frac{\cN}{\sqrt{\epsilon^{(k)}}}.
  \end{align}
  % The computational
  % complexity of each iteration of \textbf{Algorithm~APG}, displayed in
  % Figure~\ref{ch1_alg:pga}, is bounded by the complexity of computing
  % an SVD of $X\in\reals^{m\times n}$, which is
  % $\cO(\min\{nm^2,n^2m\})$.
\end{lemma}
\begin{proof}
%In the $k$-th FALC iteration, we apply \textbf{Algorithm~APG} to
%inexactly solve the optimization problem $\min_{X,s,y}\{\Pk(X,s,y):
%\norm{y}_\gamma \leq \delta\}$, i.e. the $k$-th subproblem
%\eqref{ch3_eq:augmented_lagrangian_subproblem_2}, by using the initial
%iterate $(X^{(k-1)},s^{(k-1)},y^{(k-1)})$ and setting the prox function
%to
%$h^{(k)}(X,s,y)=\frac{1}{2}\norm{X-X^{(k-1)}}_F^2+\frac{1}{2}\norm{s-s^{(k-1)}}_2^2
%+\frac{1}{2}\norm{y-y^{(k-1)}}_2^2$.
% n this proof we will bound $\ell_{\max}^{(k)}$ given in
 % \eqref{ch3_eq:apg_iter_cond}, i.e. the number of
 % \textbf{Algorithm~APG}
The number of interations $N^{(k)} \leq \ell_{\max}^{(k)}$, where $\ell_{\max}^{(k)}$
denotes the number of
iterations required to satisfy \textsc{ITERstop}.
%We terminate and set $(\Xk,\sk,\yk) =
%(X^{(k,\ell)},s^{(k,\ell)},y^{(k,\ell)})$ whenever
%\eqref{ch3_eq:modified-inner}(a) holds.
% The bound we compute is clearly an upper bound on the number of
% \textbf{Algorithm~APG} iterations required to compute the iterate
% $(\Xk,\sk,\yk)$ that satisfies one of the stopping conditions:
% \textsc{ITERstop} or \textsc{GRADstop}.
Since $(\Xk,\sk,\yk)\in\Sk$, \eqref{ch3_eq:apg_iter_cond} and
\eqref{ch3_eq:IJ_beta} imply that
\begin{align*}
\ell_{\max}^{(k)}& \leq \sqrt{\frac{2L}{\epsk}}
\left[\left(\frac{I(\alpha)}{\mu_1}+\frac{J(\beta)}{\mu_2}\right)
(\etak+\eta_1^{(k-1)}) +\etakii+\norm{y^{(k-1)}}_2\right],\\
&\leq
\sqrt{\frac{8L}{\epsk}}\left[\left(\frac{I(\alpha)}{\mu_1}+\frac{J(\beta)}
{\mu_2}\right)B_{\eta_1}
+ \frac{B_{\eta_2}+B_y}{2}\right]:=\mathcal{N}~\frac{1}{\sqrt{\epsilon^{(k)}}},
\end{align*}
where the second inequality follows from \eqref{ch3_eq:B_eta},
\eqref{eq:y-bound-uniform} and \eqref{eq:yopt-bound-uniform}. % Since
% $N^{(k)}\leq \ell^{(k)}_{\max}$ for
% all $k\geq 1$, we have
%\eqref{ch3_eq:apg_iter_bound}.
\end{proof}

In each iteration of \textbf{Algorithm~APG} %with $x$-update rule
                                %\eqref{ch1_eq:fast_x}
we need to %compute the solution of
solve one instance of each of the following
problems. % In the appendix, we prove
% Lemma~\ref{ch3_lem:unified_constrained_shrinkage}, which establishes
% the worst-case complexity of computing a solution to each of the
% following optimization problems.
\begin{enumerate}[(a)]
\item One \emph{constrained matrix shrinkage problem} of the form
 \begin{equation}
   \label{ch3_eq:matrix_constrained_shrinkage}
   \min_{X \in \reals^{m\times n}} \Big\{\lambda \norm{\sigma(X)}_{\alpha} +
   \frac{1}{2}\norm{X-\tilde{X}}_F^2: \norm{\sigma(X)}_{\alpha} \leq
   \tilde{\eta}\Big\}
 \end{equation}
 for a given $\tilde{X}\in\reals^{m\times n}$ and $\tilde{\eta}>0$. When
 $\alpha\in\{1,\infty\}$ the worst-case complexity of computing a
 solution to \eqref{ch3_eq:matrix_constrained_shrinkage} is the same as that
 of computing a full SVD, i.e. $\cO(\min\{nm^2,n^2m\})$, and when
 $\alpha = 2$, the worst-case complexity is $\cO(mn)$. See
 Lemma~\ref{ch3_lem:unified_constrained_shrinkage} in
 Appendix~\ref{app:aux-results} for details. %The basic operations'
                                %description per
 % \textbf{Algorithm~APG} iteration above implicitly assumed
 % that we need to compute this SVD exactly. This is not
 % necessary -- inexactly computing the SVD adds a small
 % additional error term to \eqref{ch3_eq:eps_stop_1}.
 Exact SVD computation is not necessary -- inexactly computing the SVD
 only adds a small additional error to \eqref{ch3_eq:eps_stop_1}.
\item  One \emph{constrained vector shrinkage problem} of the form
 \begin{equation}
   \label{ch3_eq:vector_constrained_shrinkage}
   \min_{s \in \reals^p} \Big\{\lambda \norm{s}_{\beta} +
   \frac{1}{2}\norm{s-\tilde{s}}_2^2: \norm{s}_{\beta} \leq
   \tilde{\eta}\Big\}
 \end{equation}
for a given $\tilde{s}\in\reals^p$ and $\tilde{\eta}>0$. The complexity of
solving the vector shrinkage problem
is $\cO(p\log(p))$ when $\beta \in \{1,\infty\}$ and $\cO(p)$ when
$\beta = 2$. See Lemma~\ref{ch3_lem:unified_constrained_shrinkage} in
Appendix~\ref{app:aux-results}.
\item  One \emph{Euclidean projection problem} of the form
\begin{equation}
 \label{ch3_eq:vector_constrained_shrinkage_y}
 \min_{y \in \reals^q} \big\{\frac{1}{2}\norm{y-\tilde{y}}_2^2: y\in\cQ\big\}
\end{equation}
for a given
$\tilde{y}\in\reals^q$. The complexity of solving the Euclidean projection
problem depends on $\cQ$.
%is $\cO(q\log(q))$ when $\gamma \in \{1,\infty\}$ and $\cO(q)$ when
%$\gamma = 2$ -see Lemma~\ref{ch3_lem:unified_constrained_shrinkage} in the appendix.
\end{enumerate}
In Theorem~\ref{ch3_thm:finite} we establish bounds on the
infeasibility and sub-optimality of the FALC iterate. This result
leads to a convergence rate result in
Theorem~\ref{ch3_thm:epsilon_convergence}.
\begin{theorem}
\label{ch3_thm:finite}
Let $\{(\Xk,\sk,\yk)\}_{k \in \integers_+}$ denote the sequence of
FALC iterates %generated by \ \textbf{Algorithm\ ~FALC}, displayed in
              %Figure~\ref{ch3_alg:falc}, whe
corresponding to a parameter sequence
$\{(\lk,\epsk,\tk,\xik)\}_{k\in\integers_+}$ satisfying the conditions in
Lemma~\ref{lem:bounded}.
 %Theorem~\ref{ch3_thm:limit-point}.
 In addition, suppose that, for all $k \geq 1$,  $\tk=\kappa_1 \epsk$ and
$\xik=\kappa_2 \epsk$ for some $\kappa_i\in(0,1)$ $i=1,2$. Then there
exist positive constants  $c_j$, $ j = 1, \ldots, 3$,  such
that for all $k \geq 1$,\vspace{2mm}
\begin{enumerate}[(i)]
    \item $\yk\in\cQ$ such that
      $\norm{\mathcal{A}(\Xk)-\yk-b}_2 \leq c_1
      \lambda^{(k)}$, \label{ch3_eq:infeasibility_bound}
    \item
      $\left|\left(\mu_1\norm{\sigma(\Xk)}_\alpha+\mu_2\norm{\mathcal{C}(\Xk)-d}_\beta\right)
        % \left(\mu_1\norm{\sigma(\Xopt)}_\alpha+\mu_2\norm{\mathcal{C}(\Xopt)-d}_\beta\right)
        - P^\ast
      \right|\leq
      c_2 \lambda^{(k)}+c_3\sqrt{\epsk}$, \label{ch3_eq:suboptimality_bound}
\end{enumerate}
\vspace{2mm}
where $P^\ast$ denotes the optimal value
of~\eqref{ch3_eq:simple_composite_norm_minimization}.
\end{theorem}
\begin{proof}
For all parameter sequences $\{(\lk,\epsk,\tk,\xik)\}_{k\in\integers_+}$
satisfying the conditions in Lemma~\ref{lem:bounded}, we show in
\eqref{ch3_eq:theta_bound} that $\norm{\theta_i^{(k)}}_2\leq B_{\theta_i}$
for $i\in\{1,2\}$. Therefore,
\eq
\begin{array}{rcl}
  \norm{\mathcal{A}(\Xk)-\yk-b}_2 & \leq &
  \norm{\mathcal{A}(\Xk)-\yk-b-\lambda^{(k)}\theta_1^{(k)}}_2 +
  \lambda^{(k)}\norm{\theta_1^{(k)}}_2,\\
  & = & \lambda^{(k)} \norm{\theta_1^{(k+1)}}_2 +
  \lambda^{(k)}\norm{\theta_1^{(k)}}_2,\\
  & \leq &  2B_{\theta_1}
  \lambda^{(k)}.
\end{array}
\en
This establishes (\ref{ch3_eq:infeasibility_bound}).
% where the last inequality follows the fact $\norm{\theta_1^{(k)}}_2
% \leq B_{\theta_1}$ for all $k\geq 1$; see \eqref{ch3_eq:theta_bound} for
%   details.

In the rest of the proof, we establish
\eqref{ch3_eq:suboptimality_bound}. Let $(\Xopt,\sopt,\yopt)$ denote any optimal
solution of \eqref{ch3_eq:simple_composite_norm_minimization}, i.e.
$P^\ast =
\mu_1\norm{\sigma(\Xopt)}_\alpha +\mu_2\norm{\mathcal{C}(\Xopt)-d}_\beta$. In \eqref{ch3_eq:falc_iter_bound} and
\eqref{eq:y-bound-uniform} we establish that  $\{(\Xk,\sk,\yk)\}_{k\in\integers_+}$
is a bounded sequence.
Therefore, there exists $\Gamma<\infty$ such
that for all $k\geq 1$
\begin{equation}
\label{ch3_eq:iter_bound}
\max\Big\{\norm{X_*}_F, \norm{s_*}_2, \norm{\yopt}_2, \norm{\Xk}_F,
\norm{\sk}_2, \norm{\yk}_2\Big\}\leq \Gamma.
\end{equation}
Consider the following two cases:
\begin{enumerate}[(a)]
\item  \emph{The $k$-th call to \textbf{Algorithm~APG} terminates
    with \textsc{ITERstop}}:  %$(\Xk,\sk,\yk)$ satisfies
                              %\eqref{ch3_eq:modified-inner}(a). Then,
  From \eqref{ch3_eq:inexact_minimizer_bound} it follows that
  \begin{equation}
  \label{ch3_eq:L1_norm_bound_a}
    \mu_1\norm{\sigma(X^{(k)})}_\alpha + \mu_2\norm{\sk}_\beta \leq
    \mu_1\norm{\sigma(X_*)}_\alpha +\mu_2\norm{\mathcal{C}(\Xopt)-d}_\beta
    +\frac{\lambda^{(k)}}{2}\left(\norm{\theta_1^{(k)}}_2^2+\norm{\theta_2^{(k)}}_2^2\right)
    +\frac{\epsilon^{(k)}}{\lambda^{(k)}}.
 \end{equation}
 \item \emph{The $k$-th call to \textbf{Algorithm~APG} terminates with
     \textsc{GRADstop}}: %$(\Xk,\sk,\yk)$ satisfies
                        %\eqref{ch3_eq:modified-inner}(b).
 Let $(G,g)$ belong to the set of partial subgradients $\partial_{X,s}
 P^{(k)}(.,.,.)|_{(X^{(k)},\sk,\yk)}$, defined in
 \eqref{ch3_eq:partial_subgradient}, and satisfy the stopping condition
 \textsc{GRADstop}. Then from the convexity of $P^{(k)}$ and
 Lemma~\ref{lem:grady_bound}, it follows that
   \begin{eqnarray}
     \lefteqn{P^{(k)}(\Xk,\sk,\yk)}\nonumber\\
    &\leq~& P^{(k)}(\Xopt,\sopt,\yopt)-\fprod{G, \Xopt-\Xk}-
    g^T(\sopt-\sk)-\grad_y P^{(k)}(\Xk,\sk,\yk)^T(\yopt-\yk),\nonumber\\
    &\leq~& P^{(k)}(\Xopt,\sopt,\yopt)+\norm{G}_F \norm{\Xopt-\Xk}_F+
    \norm{g}_2\norm{\sopt-\sk}_2+L\xik\norm{\yopt-\yk}_2\nonumber\\
    && +\xik\norm{\grad_y f^{(k)}(\Xk,\sk,\yk)}_2,\nonumber\\
    &\leq~&
    P^{(k)}(\Xopt,\sopt,\yopt)+\tk\left(\norm{\Xopt-\Xk}_F
      +\norm{\sopt-\sk}_2\right)+L\xik\norm{\yopt-\yk}_2\nonumber\\
    &&+\xik\norm{\grad_y
      f^{(k)}(\Xk,\sk,\yk)}_2.\label{ch3_eq:convexity_bound}
  \end{eqnarray}
  %where  and $\partial_{X,s} P^{(k)}(.,.,.)|_{(\Xk,\sk,\yk)}$ denotes the
  %set of partial subgradients of the function $P^{(k)}$ at
  %$(\Xk,\sk,\yk)$.
  Dividing \eqref{ch3_eq:convexity_bound} by $\lk$ and using the
  fact that $\theta_1^{(k+1)}=\grad_y f^{(k)}(\Xk,\sk,\yk)/\lk$, it
  follows that
  \begin{eqnarray}
  \label{ch3_eq:L1_norm_bound_b}
    \mu_1\norm{\sigma(X^{(k)})}_\alpha + \mu_2\norm{\sk}_\beta &\leq
    &\mu_1\norm{\sigma(X_*)}_\alpha +
    \mu_2\norm{\mathcal{C}(\Xopt)-d}_\beta+\frac{\lambda^{(k)}}{2}
    \left(\norm{\theta_1^{(k)}}_2^2+\norm{\theta_2^{(k)}}_2^2\right)
    +\xik\norm{\theta_1^{(k+1)}}_2
    \nonumber \\
    && \mbox{} +
    \frac{\tk}{\lk}\left(\norm{\Xopt-\Xk}_F+\norm{\sopt-\sk}_2\right)
    +\frac{\xik}{\lk}L\norm{\yopt-\yk}_2.
 \end{eqnarray}
\end{enumerate}
From \eqref{ch3_eq:iter_bound}, the fact that
$\norm{\theta_i^{(k)}}_2\leq B_{\theta_i}$ for $i\in\{1,2\}$,
\eqref{ch3_eq:L1_norm_bound_a} and \eqref{ch3_eq:L1_norm_bound_b}, it
follows that
 \begin{eqnarray}
    \label{ch3_eq:L1_norm_bound}
    \hspace{-6mm}\mu_1\norm{\sigma(X^{(k)})}_\alpha + \mu_2\norm{\sk}_\beta &\leq
    &\mu_1\norm{\sigma(X_*)}_\alpha +
    \mu_2\norm{\mathcal{C}(\Xopt)-d}_\beta+\left(\frac{B_{\theta_1}^2
        +B_{\theta_2}^2}{2}\right)\lk \nonumber\\
    &&+\max\left\{\frac{\epsilon^{(k)}}{\lambda^{(k)}},
      \frac{4\Gamma\tk}{\lk}+\frac{\xik}{\lk}(2L\Gamma+\lk
      B_{\theta_1})\right\}.
 \end{eqnarray}
 %From \eqref{ch3_eq:IJ_beta} and Corollary \ref{ch3_cor:grad_norm_bound} in
% Appendix~A, it follows that for all $k\geq 1$,
% \begin{eqnarray}
%   \norm{\mathcal{C}(\Xk)+\sk-d-\lk\tetkii}_\beta & \leq &
%   J(\betas)~\norm{\mathcal{C}(\Xk)+\sk-d-\lk\tetkii}_2 \nonumber\\
%   &\leq&
%   J(\betas)\left(\sqrt{2\epsk}\sigma_{\max}(M)+J(\betas)\mu_2\lk\right). \label{ch3_eq:feasibilty_A_C}
% \end{eqnarray}
  The bound \eqref{ch3_eq:IJ_beta} relating $\norm{\sigma(.)}_\alpha$
  %singular value norms
  to the Frobenius norm,  \eqref{ch3_eq:feasibility_bound_1},
   bound $\norm{\theta_2^{(k)}}_2
  \leq B_{\theta_2}$, together with triangle inequality imply that
  \begin{eqnarray}
    \hspace{-6mm}\norm{\mathcal{C}(\Xk)-d}_\beta &\leq& \norm{\sk}_\beta +
    \norm{\lk\tetkii}_\beta + J(\betas)\left(\max\left\{\sqrt{2\epsk}
      \sigma_{\max}(M),\ \tk\right\}+J(\betas)\mu_2\lk\right),\nonumber\\
    &\leq& \norm{\sk}_\beta +
    J(\betas)\left(B_{\theta_2}+\mu_2~J(\betas)\right)\lk+J(\betas)
    \max\left\{\sqrt{2}~\sigma_{\max}(M),\
      \kappa_1\sqrt{\epsk}\right\}~\sqrt{\epsk},
    \nonumber %\label{ch3_eq:s_bound}
  \end{eqnarray}
  where the second inequality uses the relation $\tk=\kappa_1\epsk$.
  The above inequality and \eqref{ch3_eq:L1_norm_bound} imply that
  \begin{eqnarray}
    \mu_1\norm{\sigma(X^{(k)})}_\alpha+\mu_2\norm{\mathcal{C}(\Xk)-d}_\beta
    &\leq&\mu_1\norm{\sigma(X_*)}_\alpha
    +\mu_2\norm{\mathcal{C}(\Xopt)-d}_\beta \nonumber\\
    && + \left(\frac{B_{\theta_1}^2 +
        B_{\theta_2}^2}{2}+\mu_2J(\betas)\left(B_{\theta_2} +
        \mu_2~J(\betas)\right)\right)\lambda^{(k)} \nonumber\\
    &&   +\max\left\{\frac{\epsilon^{(k)}}{\lambda^{(k)}},
      \frac{4\Gamma\tk}{\lk}+\frac{\xik}{\lk}(2L\Gamma+\lk
      B_{\theta_1})\right\} \nonumber\\
    &&+\mu_2 J(\betas)\max\left\{\sqrt{2}~\sigma_{\max}(M),\
      \kappa_1\sqrt{\epsk}\right\}~\sqrt{\epsk}.
    \label{ch3_eq:approx_suboptimality_upper_bound}
  \end{eqnarray}
  Since %$\{(\Xk,\sk,\yk)\}_{k\in\integers_+}$ is a bounded sequence,
  $\frac{\epsilon^{(k)}}{(\lambda^{(k)})^2}\leq B$, $\tk=\kappa_1 \epsk$
  and $\xik=\kappa_2 \epsk$ for all $k\geq 1$,
  \eqref{ch3_eq:approx_suboptimality_upper_bound} implies one side of the
  bound in \eqref{ch3_eq:suboptimality_bound}.

  Next, we establish a lower bound for $P^{(k)}(\Xkopt,\skopt,\ykopt)$ by
  comparing the following pairs of Lagrangian duals
  \begin{subequations}
    \label{ch3_eq:primal_dual_pair}
    \begin{equation}
      \label{ch3_eq:primal_problem}
      \begin{array}[t]{rl}
        \min_{X\in\reals^{m\times n}} &
        \mu_1\norm{\sigma(X)}_\alpha+\mu_2\norm{\mathcal{C}(X)-d}_\beta,\\
        \mbox{s.t.} & \mathcal{A}(X)- b\in\cQ.
      \end{array}
    \end{equation}
    \begin{equation}
      \label{ch3_eq:dual_problem}
      \begin{array}[t]{rl}
        \max_{w\in\reals^q, v\in\reals^p} & -b^Tw-d^Tv-\gamma_\cQ(w),\\
        \mbox{s.t.} & \norm{\sigma(\mathcal{A}^*(w)+\mathcal{C}^*(v))}_{\alphas} \leq \mu_1,\\
        &\norm{v}_{\betas} \leq \mu_2,
      \end{array}
    \end{equation}
  \end{subequations}
  where $\gamma_\cQ$ is the support function of $\cQ$, i.e.,
  $\gamma_\cQ(w):=\sup_{y\in\cQ}w^Ty$, and
  \begin{subequations}
    \label{ch3_eq:penalty_primal_dual}
    \begin{equation}
      \label{ch3_eq:penalty_primal_problem}
      \begin{array}[t]{rl}
        \min_{X\in\reals^{m\times n}, s\in\reals^p, y\in\reals^q} &
        \lambda (\mu_1\norm{\sigma(X)}_\alpha+\mu_2\norm{s}_\beta) +
        f(X,s,y),\\
        \mbox{s.t.} & y\in\cQ.
      \end{array}
    \end{equation}
    \begin{equation}
      \label{ch3_eq:penalty_dual_problem}
      \begin{array}[t]{rl}
        \max_{w\in\reals^q, v\in\reals^p}&
        -\lambda (b+\lambda\theta_1)^Tw-\lambda (d+\lambda\theta_2)^Tv
        -\lambda\gamma_\cQ(w)-\frac{\lambda^2}{2}(\norm{w}_2^2+\norm{v}_2^2),\\
        \mbox{s.t.} &
        \norm{\sigma(\mathcal{A}^*(w)+\mathcal{C}^*(v))}_{\alphas} \leq
        \mu_1,\\
        &\norm{v}_{\betas} \leq \mu_2.
      \end{array}
    \end{equation}
  \end{subequations}
  Above $(w_*,v_*)$ denotes the optimal solution of the dual
  \eqref{ch3_eq:dual_problem} and
  $f(X,s,y):=\frac{1}{2}\norm{\mathcal{A}(X)+y-b-\lambda\theta_1}_2^2 +
  \frac{1}{2}\norm{\mathcal{C}(X)+s-d-\lambda\theta_2}_2^2$.
  %Moreover, \eqref{ch3_eq:penalty_primal_problem} and \eqref{ch3_eq:penalty_dual_problem} below are also a Lagrange primal-dual pair of problems.
  Note that $(w_*,v_*)$ is  feasible for %  the maximization problem
  \eqref{ch3_eq:penalty_dual_problem}. Therefore, by Lagrangian duality it follows that
  \begin{eqnarray}
    \lefteqn{P^{(k)}(\Xkopt,\skopt,\ykopt)} \nonumber\\
    %&= & \min_{X\in\reals^{m\times n},
%      s\in\reals^p,y\in\reals^q}\Big\{\lambda^{(k)}(\mu_1\norm{\sigma(X)}_\alpha
%    + \mu_2\norm{s}_\beta)+f^{(k)}(X,s,y):\ \norm{y}_\gamma\leq\delta)\Big\}
%    \nonumber\\
    &\geq & \lambda^{(k)} \left(-b^Tw_* - d^Tv_* -\gamma_\cQ(w_*)-
      \frac{\lambda^{(k)}}{2}\left(\norm{w_*}_2^2
    + \norm{v_*}_2^2 + 2(\tetki)^T w_* + 2 (\tetkii)^T
    v_*\right)\right),\nonumber\\
    &\geq & \lambda^{(k)} \left( \mu_1\norm{\sigma(X_*)}_\alpha +
    \mu_2\norm{\mathcal{C}(X_*)-d}_\beta\right) \nonumber\\
    &&-\frac{(\lambda^{(k)})^2}{2}\left(\norm{w_*}_2^2
      +\norm{v_*}_2^2+2\norm{\tetki}_2 \norm{w_*}_2+2
      \norm{\tetkii}_2\norm{v_*}_2\right),\label{ch3_eq:dual_cauchy}
  \end{eqnarray}
  where the first inequality follows from weak duality for primal-dual pair
  in \eqref{ch3_eq:penalty_primal_dual}, and \eqref{ch3_eq:dual_cauchy} follows
  from strong duality for primal-dual pair in
  \eqref{ch3_eq:primal_dual_pair} % ,
%   i.e. $b^Tw_*+d^Tv_*-\delta\norm{w_*}_{\gammas}=\mu_1\norm{X_*}_\alpha
%   +\mu_2\norm{\mathcal{C}(X_*)-d}_\beta$,
  and the Cauchy-Schwartz inequality.

  From the definition of $\{\theta_i^{(k)}\}_{k\in\integers_+}$ in
  Figure~\ref{ch3_alg:falc}, it is clear that the FALC iterates
  $\{\Xk\}_{k \in \integers}$ satisfy
  \eq
  \frac{P^{(k)}(\Xk,\sk,\yk)}{\lk}  =
  \big(\mu_1\norm{\sigma(\Xk)}_\alpha+\mu_2\norm{\sk}_\beta\big)
  +\frac{\lk}{2}\left(\norm{\theta_1^{(k+1)}}_2^2+\norm{\theta_2^{(k+1)}}_2^2\right),
  \en
  and it follows that
  \begin{equation}
    \mu_1\norm{\sigma(\Xk)}_\alpha + \mu_2\norm{\sk}_\beta \geq
    \frac{P^{(k)}(\Xkopt,\skopt,\ykopt)}{\lk}-\frac{\lk}{2}
    \left(\norm{\theta_1^{(k+1)}}_2^2+\norm{\theta_2^{(k+1)}}_2^2\right).
    \label{ch3_eq:lbd-1}
  \end{equation}

  Thus, the bound on $\norm{\theta_i^{(k)}}_2$, $i\in\{1,2\}$
  established in \eqref{ch3_eq:theta_bound}, and the inequalities
  \eqref{ch3_eq:dual_cauchy} and  \eqref{ch3_eq:lbd-1}, together imply that
  \begin{eqnarray}
    \mu_1\norm{\sigma(X^{(k)})}_\alpha + \mu_2\norm{\sk}_\beta
    & \geq &
    \mu_1\norm{\sigma(X_*)}_\alpha + \mu_2\norm{\mathcal{C}(X_*)-d}_\beta
    \nonumber\\
    && -\frac{\lk}{2}\left((B_{\theta_1} +
      \norm{w_*}_2)^2+(B_{\theta_2}+\norm{v_*}_2)^2\right).
    \label{ch3_eq:approx_suboptimality_lower_bound}
  \end{eqnarray}
  The bound
  $\norm{\mathcal{A}^*(w_*)+\mathcal{C}^*(v_*)}_{F} \leq I(\alphas)
  \norm{\sigma(\mathcal{A}^*(w_*)+\mathcal{C}^*(v_*))}_{\alphas}\leq
  I(\alphas)\mu_1$ implies that
  \eq
  \sigma_{\min}(A) \norm{w_{\ast}}_2 \leq \norm{\mathcal{A}^*(w_*)}_F
  \leq  I(\alphas)\mu_1 + \norm{\mathcal{C}^*(v_*)}_{F}
  \leq  I(\alphas)\mu_1 + \sigma_{\max}(C)\norm{v_{\ast}}_2,
  \en
  and the bound $\norm{v_*}_{\betas} \leq \mu_2$ implies that
  $\norm{v_*}_2\leq J(\betas)~\mu_2$. Hence, both $\norm{v_*}_2$ and $\norm{w_*}_2$ in \eqref{ch3_eq:approx_suboptimality_lower_bound} are bounded.

  The bound \eqref{ch3_eq:IJ_beta}, the uniform bound $\norm{\theta_2^{(k)}}_2
  \leq B_{\theta_2}$ %, for all $k\geq 1$, established in see \eqref{ch3_eq:theta_bound}, \eqref{ch3_eq:feasibility_bound_1}
  and triangle inequality together imply that
  \begin{equation}
    \norm{\sk}_\beta \leq \norm{\mathcal{C}(\Xk)-d}_\beta +
    J(\betas)\left(B_{\theta_2}+\mu_2~J(\betas)\right)\lk
    +J(\betas)\max\left\{\sqrt{2}~\sigma_{\max}(M),\
      \kappa_1\sqrt{\epsk}\right\}~\sqrt{\epsk}, \label{ch3_eq:C_bound}
  \end{equation}
  where the second inequality uses the relation $\tk=\kappa_1\epsk$.
  From \eqref{ch3_eq:approx_suboptimality_lower_bound} and
  \eqref{ch3_eq:C_bound}, it follows that
  \begin{eqnarray*}
    \lefteqn{\mu_1\norm{\sigma(X^{(k)})}_\alpha
      +\mu_2\norm{\mathcal{C}(\Xk)-d}_\beta}&&\\
    & \geq &
    \mu_1\norm{\sigma(X_*)}_\alpha+\mu_2\norm{\mathcal{C}(\Xopt)-d}_\beta\\
    &&\mbox{}   -\left(\frac{(B_{\theta_1}+\norm{w_*}_2)^2 +
        (B_{\theta_2}+\norm{v_*}_2)^2}{2}
    + \mu_2J(\betas)\left( B_{\theta_2}+
      \mu_2~J(\betas)\right)\right)\lambda^{(k)}\\
    && \mbox{} -\mu_2 J(\betas)\max\left\{\sqrt{2}~\sigma_{\max}(M),\
      \kappa_1\sqrt{\epsk}\right\}~\sqrt{\epsk}.
  \end{eqnarray*}
  This establishes the result.
\end{proof}

Now, we have all the estimates we need to prove the main convergence rate
result in this paper.
\begin{theorem}
  \label{ch3_thm:epsilon_convergence}
  Fix $\kappa_1,\kappa_2,\nu\in(0,1)$, and strictly positive parameters
  $(\lambda^{(1)},\epsilon^{(1)},\tau^{(1)},\xi^{(1)})$. For $k \geq 1$, set
  parameter sequence $\{(\lk,\epsk,\tk,\xik)\}_{k\in\integers_+}$ as follows:
  \begin{equation}
  \label{ch3_eq:param-update}
  \begin{array}{rclrcl}
    \lambda^{(k+1)} &= & \nu\; \lambda^{(k)}, \qquad \xik &= &
    \kappa_2\; \epsilon^{(k)},\\
    \epsilon^{(k+1)} &= & \nu^2\; \epsilon^{(k)},\qquad \tk &= &
    \kappa_1\; \epsilon^{(k)}.
  \end{array}
\end{equation}
% where $\Gamma$ is defined in \eqref{ch3_eq:iter_bound}.
Then, for all $\epsilon > 0$,   \textbf{Algorithm~FALC}
computes an $\epsilon$-feasible and $\epsilon$-optimal solution
$\bar{X}\in\reals^{m\times n}$ to the composite norm minimization
problem~\eqref{ch3_eq:simple_composite_norm_minimization}, %  i.e. for some
% $\bar{y}\in\reals^q$ such that $\bar{y}\in\cQ$, we have
%   \eq
%   \norm{\mathcal{A}(\bar{X})-\bar{y} - b}_2  \leq \epsilon, \qquad
%   \left|(\mu_1\norm{\sigma(\bar{X})}_\alpha + \mu_2\norm{\mathcal{C}(\bar{X})-d}_\beta) -
%       (\mu_1\norm{\sigma(\Xopt)}_\alpha +
%       \mu_2\norm{\mathcal{C}(\Xopt)-d}_\beta)\right| \leq \epsilon,
%   \en
in $N_{\rm inner}=\cO\Big( \frac{1}{\epsilon} \Big)$
\textbf{Algorithm~APG} iterations.
\end{theorem}
\begin{proof}
For the specific choice of the parameter sequence in
\eqref{ch3_eq:param-update} we have that
$\frac{\epsilon^{(k)}}{(\lambda^{(k)})^2}=\frac{\epsilon^{(1)}}{(\lambda^{(1)})^2}$,
for all $k\geq 1$. % Hence, setting
% $B=\frac{\epsilon^{(1)}}{(\lambda^{(1)})^2}$,
Therefore, Theorem~\ref{ch3_thm:finite}
guarantees that there exist $c_2>0$ and $c_3>0$ such that for all $k \geq 1$,
\begin{align*}
   \big|(\mu_1\norm{\sigma(\Xk)}_\alpha+\mu_2\norm{\mathcal{C}(\Xk)-d}_\beta)
   -(\mu_1\norm{\sigma(\Xopt)}_\alpha
   &
   +\mu_2\norm{\mathcal{C}(\Xopt)-d}_\beta)\big|
   \leq\big(c_2 \lambda^{(1)}+c_3\sqrt{\epsilon^{(1)}}~\big)\;\nu^{(k-1)}.
\end{align*}
Thus, $
\big|(\mu_1\norm{\sigma(\Xk)}_\alpha+\mu_2\norm{\mathcal{C}(\Xk)-d}_\beta)
  -
  (\mu_1\norm{\sigma(\Xopt)}_\alpha
  +\mu_2\norm{\mathcal{C}(\Xopt)-d}_\beta)\big|
\leq \epsilon$,
for all % $k\in\integers_+$ such that
\begin{equation}
  \label{ch3_eq:K-suboptimality}
  k > \log_{\frac{1}{\nu}}\left(\frac{c_2 \lambda^{(1)} +
        c_3\sqrt{\epsilon^{(1)}}}{\epsilon}\right)+1.
\end{equation}
Moreover, Theorem~\ref{ch3_thm:finite} also implies that there exists
$c_1>0$ such that
$
\norm{\mathcal{A}(\Xk)-\yk-b}_2 \leq c_1 \lambda^{(1)}\; \nu^{k-1},
$
for $k \geq 1$,
Thus, $\norm{\mathcal{A}(X^{(k)})-\yk-b}_2 \leq \epsilon$ for all % $k\in\integers_+$ such that
\begin{equation}
  \label{ch3_eq:K-infeasibility}
  k > \log_{\frac{1}{\nu}}\left(\frac{c_1
        \lambda^{(1)}}{\epsilon}\right)+
  1.
\end{equation}
Then \eqref{ch3_eq:K-suboptimality} and \eqref{ch3_eq:K-infeasibility}
imply that for all $\epsilon > 0$, the number of FALC iterations  required
to compute an $\epsilon$-feasible and
$\epsilon$-optimal solution
\begin{equation}
  \label{ch3_eq:nout-bnd}
  N_{\rm FALC}(\epsilon) \leq
  \left\lceil\log_{\frac{1}{\nu}}\left(\frac{U}{\epsilon}\right)\right\rceil+1,
\end{equation}
where $U = \max\left\{c_2 \lambda^{(1)}+c_3\sqrt{\epsilon^{(1)}},\; c_1 \lambda^{(1)}\right\}$.

From Lemma~\ref{ch3_lem:tseng} it follows that there exists a constant
$\cN$ such that APG iteration
to solve the $k$-th FALC subproblem $N^{(k)} \leq
\frac{\cN}{\sqrt{\epsilon^{(k)}}}$. Therefore,
\eq
N_{\rm inner} = \cN
\sum_{k=1}^{N_{\rm FALC}(\epsilon)}\frac{1}{\sqrt{\epsilon^{(k)}}} =
\frac{\cN}{\sqrt{\epsilon^{(1)}}} \sum_{k=0}^{N_{\rm FALC}(\epsilon)-1}\nu^{-k}
  =\frac{\cN}{\sqrt{\epsilon^{(1)}}}\cdot
  \frac{\nu}{(1-\nu)}\cdot\left(\frac{1}{\nu}\right)^{N_{\rm
      FALC}(\epsilon)} \leq \left(\frac{\cN
      U}{\nu(1-\nu)~\sqrt{\epsilon^{(1)}}}\right)~\frac{1}{\epsilon},
\en
where the last bound follows from (\ref{ch3_eq:nout-bnd}). % it follows
% that for all $\epsilon>0$ an
% $\epsilon$-feasible and
% $\epsilon$-optimal solution can be
% computed in at most
% \eq
% N_{\rm inner}\leq \left(\frac{\cN
% U}{\nu(1-\nu)~\sqrt{\epsilon^{(1)}}}\right)
% ~\frac{1}{\epsilon}=\cO\left(\frac{1}{\epsilon}\right)
% \en
% iterations of \textbf{Algorithm~APG},
% %with $x$-update rule \eqref{ch1_eq:fast_x}
% where each iteration requires $\cO(\min\{nm^2,$ $n^2m\})$ work for one SVD computation.
\end{proof}

Note that we do not explicitly specify the constant hidden in the
$\cO(1/\epsilon)$ iteration complexity result of
Theorem~\ref{ch3_thm:epsilon_convergence}. Moreover, the bound given in the
proof is very crude. The main reason is that the
composite minimization problem
\eqref{ch3_eq:composite_norm_minimization} is very general and the
constant strongly depends on the specific problem
structure. However, the proof technique used to establish
Theorem~\ref{ch3_thm:epsilon_convergence} can be applied as is to
any special case of composite minimization problem to obtain much
sharper constants. For instance, the
complexity result of FALC for the basis pursuit problem in
\eqref{ch3_eq:l1_minimization} is given by
$$N_{\rm inner}\leq n\kappa(A)^2 \left( \frac{16
    \norm{\xbp}_1}{\nu(1-\nu)} \cdot\frac{1}{\epsilon}+
  \frac{9}{\nu}\cdot
  \log_{\frac{1}{\nu}}\left(\frac{8n\kappa^2(A)}{\epsilon}\right)
\right) =\cO\left(\frac{1}{\epsilon}\right),$$ where
$\kappa(A):=\sigma_{\max}(A)/\sigma_{\min}(A)$ is the condition
number of $A$. This is the same bound that was obtained for this
special case in~\cite{Aybat12}.

The convergence rate result in
Theorem~\ref{ch3_thm:epsilon_convergence} relies on the uniform bound
established in Lemma~\ref{ch3_lem:tseng}. This uniform bound in turn
assumes that all calls to \textbf{Algorithm APG} are terminated by
\textsc{ITERstop}. On the other hand, in our numerical experiments almost all
calls to \textbf{Algorithm APG} were terminated by \textsc{GRADstop}. This suggests that the
$\cO(\frac{1}{\epsilon})$  rate result has a lot of slack. Indeed, we
find that early terminating \textbf{Algorithm APG} iterations when the stopping condition \textsc{GRADstop} is satisfied reduces total number of \textbf{Algorithm APG} iterations significantly: in our numerical tests, FALC required only $\cO(\log(\frac{1}{\epsilon}))$ inner iterations to compute an $\epsilon$-optimal, $\epsilon$-feasible iterate. The augmented Lagrangian algorithm~FAL introduced in~\cite{Aybat12} is
an implementation of FALC for the basis pursuit problem. In Figure 6.1 in \cite{Aybat12} one
can clearly observe the $\cO(\log(1/\epsilon))$ empirical performance as opposed to the $\cO(1/\epsilon)$ worst case complexity.
We observe a similar empirical performance with FALC on the numerical problems tested in this paper.

\section{Implementation details of Algorithm FALC}
\label{ch3_sec:implementation}
In this section we describe all the details of FALC. Let
$\{(X_i^{(k,\ell)},s_i^{(k,\ell)},y_i^{(k,\ell)})\}_{\ell\in\integers_+}$
denote the sequence of $x^{(\ell)}_i$-iterates of \textbf{Algorithm~APG}
in Figure~\ref{ch1_alg:pga} for $i\in\{1,2\}$ when \textbf{Algorithm~APG} is
called in Line~\ref{algeq:apg_call} of Figure~\ref{ch3_alg:falc} to solve
the $k$-th subproblem.

% When we implemented \textbf{Algorithm~FALC}, instead of
% \textsc{GRADstop1} given in Line~\ref{ch3_li:gradstop1} of
% Figure~\ref{ch3_alg:falc},

\subsection{Subgradient selection}
\label{ch3_sec:subgradient_selection}
%If the iterate $\left(\Xk,\sk,\yk\right)$ is required to satisfy only \eqref{ch3_eq:modified-inner}(a), i.e.
%$P^{(k)}(\Xk,\sk,\yk) \leq \inf\{P^{(k)}(X,s,y):$ $X\in\reals^{m\times
%n},s\in\reals^p,\norm{y}_\gamma\leq\delta\} + \epsilon^{(k)}$ which is
%implied by the first stopping criterion given in
%step~\ref{ch3_li:inner-stopping-condition} of \textbf{Algorithm~FALC}
%displayed in Figure~\ref{ch3_alg:falc},

We used the following  slightly modified version of \textsc{GRADstop} in
our implementation.
\begin{equation}
\label{ch3_eq:modified_gradstop1}
\textsc{GRADstop1}:=\left\{\exists (G,g)\in\partial_{X,s}
  P^{(k)}(.,.,.)|_{(\Xk,\sk,\yk)} \mbox{ s.t. }
  \norm{G}_F\leq\tau_X^{(k)}\mbox{ \textbf{and} }
  \norm{g}_2\leq\tau_s^{(k)}\right\}
\end{equation}
for some tolerance values $\tau_X^{(k)}$ and $\tau_s^{(k)}$ such that
$\{\tau_X^{(k)}\}_{k\in\integers_+}$ and
$\{\tau_s^{(k)}\}_{k\in\integers_+}$ are decreasing sequences. Clearly, if
\eqref{ch3_eq:modified_gradstop1} holds, the original \textsc{GRADstop1}
given in Line~\ref{ch3_li:gradstop1} of \textbf{Algorithm~FALC} holds for
$\tk=\tau_X^{(k)}+\tau_s^{(k)}$.

We check the stopping condition \textsc{GRADstop1} in each
\textbf{Algorithm~APG} iteration. Let $\breve{Z}^{(k,\ell)} =\left(\breve{X}_1^{(k,\ell)},\breve{s}_1^{(k,\ell)},
  \breve{y}_1^{(k,\ell)}\right)$
denotes the unconstrained solution to the optimization problem in
Line~\ref{ch1_algeq:tseng_z} of
Figure~\ref{ch1_alg:pga}, i.e. when the contraint $(X,s,y)\in\Sk$  is not
enforced. A subgradient $(G,g)\in\partial_{X,s}
P^{(k)}(.,.,.)|_{\breve{Z}^{(k,\ell)}}$ can be computed as follows %in
% line~\ref{ch3_li:G} and
% line~\ref{ch3_li:g}
\eq
G = \lk \mu_1 Q + \grad_X
f^{(k)}\left(\breve{Z}^{(k,\ell)}\right)\ \mbox{ and }\ g = \lk \mu_2 q +
\grad_s f^{(k)}\left(\breve{Z}^{(k,\ell)}\right),
\en
where
\begin{align*}
Q &= \frac{L}{\lambda^{(k)}\mu_1} \left(X_2^{(k,\ell)}-\frac{1}{L}\grad_X
  f^{(k)}(X_2^{(k,\ell)},s_2^{(k,\ell)},y_2^{(k,\ell)})-\breve{X}_1^{(k,\ell)}\right),\\
q &= \argmin\left\{\norm{\lk \mu_2 r+\grad_s
    f^{(k)}\left(\breve{Z}^{(k,\ell)}\right)}_2: r\in\partial
\norm{.}_\beta|_{\breve{s}_2^{(k,\ell)}}\right\}.
\end{align*}
From the first order optimality condition, it can be easily shown that $Q\in \partial \norm{\sigma(.)}_\alpha |_{\breve{X}_1^{(k,\ell)}}$. Moreover, given $\grad_s f^{(k)}$ at $\breve{Z}^{(k,\ell)}$ the complexity of computing $q\in\partial \norm{.}_\beta|_{\breve{s}_1^{(k,\ell)}}\subset \reals^p$ is  $\cO(p)$ when $\beta\in\{1,2\}$ and $\cO(p\log(p))$  when $\beta=\infty$.
\subsection{Stopping criterion for FALC}
In our numerical experiments, we terminate \textbf{Algorithm~FALC} either the distance between
successive inner iterates are below a threshold $\varrho$ for each
component, i.e.
$\norm{X_1^{(k,\ell)}-X_1^{(k,\ell-1)}}_F\leq\varrho$,
$\norm{s_1^{(k,\ell)}-s_1^{(k,\ell-1)}}_2\leq\varrho$ or there exist
partial subgradients with sufficiently small norm for each component,
i.e. $\norm{G}_F\leq\varsigma_X$, $\norm{g}_2\leq\varsigma_s$ for some
$(G,g)\in\partial_{X,s} P^{(k)}(.,.,.)|_{\breve{Z}^{(k,\ell)}}$ and
$$\norm{\breve{y}_1^{(k,\ell)}-\Pi_\cQ\left(\breve{y}_1^{(k,\ell)}-\frac{1}{L}\grad_y P^{(k)}(\breve{Z}^{(k,\ell)})\right)}_2 \leq\varsigma_y.$$
In our numerical experiments we set $\varrho$, $\varsigma_X$,
$\varsigma_s$ and $\varsigma_y$ by experimenting with a small instance of
the problem.
\subsection{Multiplier selection}
\label{ch3_sec:multiplier_selection}
Given $\bar{c}_\tau\in(0,1)$, $\bar{c}_\xi\in(0,1)$, $\bar{c}_\lambda>0$,
$c_\tau\in(0,1)$, $c_\xi\in(0,1)$,
$c_\lambda\in(0,1)$, for all $k\geq 1$ the approximate optimality
parameters $\tau_X^{(k)}$, $\tau_s^{(k)}$, $\xik$ and the penalty parameter $\lk$ are set
as follows:
\begin{equation*}
  \begin{array}{rcl}
    \breve{X}^{(1)} & = & \argmin_{X\in\reals^{m\times n}}
    \frac{1}{2}\norm{X-\left(X^{(0)}-\frac{1}{L}\grad_X
        f^{(1)}(X^{(0)},s^{(0)},y^{(0)})\right)}_F^2
    +\frac{\lambda^{(1)}\mu_1}{L}\norm{\sigma(X)}_\alpha,\\
    \breve{s}^{(1)} & = & \argmin_{s\in\reals^p}
    \frac{1}{2}\norm{s-\left(s^{(0)}-\frac{1}{L}\grad_s
        f^{(1)}(X^{(0)},s^{(0)},y^{(0)})\right)}_2^2+
    \frac{\lambda^{(1)}\mu_2}{L}\norm{s}_\beta,\\
    \breve{y}^{(1)} & = & \argmin_{y\in\cQ\subset\reals^q}
    \norm{y-\left(y^{(0)}-\frac{1}{L}\grad_y
        f^{(1)}(X^{(0)},s^{(0)},y^{(0)})\right)}_2^2,\\
    G^{(1)}& = & L\left(X^{(0)}-\frac{1}{L}\grad_X
      f^{(1)}(X^{(0)},s^{(0)},y^{(0)})-\breve{X}^{(1)}\right)+ \grad_X
    f^{(1)}(\breve{X}^{(1)},\breve{s}^{(1)},\breve{y}^{(1)}),\\
    g^{(1)} &=& \argmin_{g\in\reals^p}\{\norm{g}_2: g=\lambda^{(1)}\mu_2 p+\grad_s
    f^{(1)}(\breve{X}^{(1)},\breve{s}^{(1)},\breve{y}^{(1)}),\
    p\in\partial\norm{.}_\beta|_{\breve{s}^{(1)}}\},\\
    \tau_X^{(1)} & = & \bar{c}_{\tau} \norm{G^{(1)}}_F, \\
    \tau_s^{(1)} & = & \bar{c}_{\tau} \norm{g^{(1)}}_2, \\
    \xi^{(1)} & = & \bar{c}_{\xi} \norm{\breve{y}^{(1)}-\Pi_\cQ\left(\breve{y}^{(1)}-\frac{1}{L}\grad_y
        f^{(1)}(\breve{X}^{(1)},\breve{s}^{(1)},\breve{y}^{(1)})\right)}_2,\\
      \lambda^{(1)} & = & \bar{c}_\lambda \norm{X^{(0)}}_2,\\
%    \end{array}
%  \end{equation*}
%  \begin{equation*}
%    \begin{array}{rcl}
    \breve{X}^{(k)} & = & \argmin_{X\in\reals^{m\times n}}
    \frac{1}{2}\norm{X-\left(X^{(k-1)}-\frac{1}{L}\grad_X
        f^{(k)}(X^{(k-1)},s^{(k-1)},y^{(k-1)})\right)}_F^2
    +\frac{\lambda^{(k)}\mu_1}{L}\norm{\sigma(X)}_\alpha,\\
    \breve{s}^{(k)} & = & \argmin_{s\in\reals^p}
    \frac{1}{2}\norm{s-\left(s^{(k-1)}-\frac{1}{L}\grad_s
        f^{(k)}(X^{(k-1)},s^{(k-1)},y^{(k-1)})\right)}_2^2
    +\frac{\lambda^{(k)}\mu_2}{L}\norm{s}_\beta,\\
    \breve{y}^{(k)} & = & \argmin_{y\in\cQ\subset\reals^q}
    \norm{y-\left(y^{(k-1)}-\frac{1}{L}\grad_y
        f^{(k)}(X^{(k-1)},s^{(k-1)},y^{(k-1)})\right)}_2^2,\\
    G^{(k)}& = & L\left(X^{(k-1)}-\frac{1}{L}\grad_X
      f^{(k)}(X^{(k-1)},s^{(k-1)},y^{(k-1)})-\breve{X}^{(k)}\right)+ \grad_X
    f^{(k)}(\breve{X}^{(k)},\breve{s}^{(k)},\breve{y}^{(k)}),\\
    g^{(k)} &=& \argmin_{g\in\reals^p}\{\norm{g}_2: g=\lambda^{(k)}\mu_2 p+\grad_s
    f^{(k)}(\breve{X}^{(k)},\breve{s}^{(k)},\breve{y}^{(k)}),\
    p\in\partial\norm{.}_\beta|_{\breve{s}^{(k)}}\},\\
    \tau_X^{(k)}& = &\min\big\{\  c_{\tau}\,\tau_X^{(k-1)}, \bar{c}_{\tau}
    \norm{G^{(k)}}_F\ \big\},\\
    \tau_s^{(k)}& = &\min\big\{\  c_{\tau}\,\tau_s^{(k-1)}, \bar{c}_{\tau}
    \norm{g^{(k)}}_2\ \big\},\\
    \xik & =& \min \{ c_{\xi}\,\xi^{(k-1)},\ \norm{\breve{y}^{(k)}-\Pi_\cQ\left(\breve{y}^{(k)}-\frac{1}{L}\grad_y f^{(k)}(\breve{X}^{(k)},\breve{s}^{(k)},\breve{y}^{(k)})\right)}_2\}\\
    %\end{array}
%    \end{equation*}
%    \begin{equation*}
%    \begin{array}{rcl}
    \lambda^{(k)} & = & c_\lambda \, \lambda^{(k-1)},
  \end{array}
\end{equation*}
for all $k\geq 2$. In all our experiments, $\bar{c}_\tau=0.999$ and $\bar{c}_\xi=0.999$.

We initialize FALC with $\left(X^{(0)}, s^{(0)},y^{(0)}\right)$ such that
$\mathcal{A}(X^{(0)})- b\in\cQ$, $s^{(0)}=\mathcal{C}(X^{(0)})-d$ and $y^{(0)}=\mathcal{A}(X^{(0)})-b$. In first
iteration of FALC, we solve the problem
\begin{align*}
\min_{(X,s,y)\in \cS^{(1)},~y\in\cQ}P^{(1)}(X,s,y)
=\min_{(X,s,y)\in \cS^{(1)},~y\in\cQ}\lambda^{(1)}
(\mu_1\norm{\sigma(X)}_\alpha+\mu_2\norm{s}_\beta)
+ f^{(1)}(X,s,y),
\end{align*}
where $\cS^{(1)}=\{(X,s,y):\ \mu_1\norm{\sigma(X)}_\alpha+\mu_2\norm{s}_\beta\leq\eta^{(1)}\}$ and
$\eta^{(1)}=\mu_1\norm{\sigma(X^{(0)})}_\alpha+\mu_2\norm{s^{(0)}}_\beta$. Since
$X^{(0)}$ is feasible, $f^{(1)}(X^{(0)},s^{(0)},y^{(0)})=0$ and
$P^{(1)}(X^{(0)},s^{(0)},$ $y^{(0)})= \lambda^{(1)}\eta^{(1)}$. Then $P^{(1)}(X,s,y)\geq
0$ for all $X\in\reals^{m\times n}$ implies that the initial duality gap is less than
or equal to $\lambda^{(1)}\eta^{(1)}$. Hence, we initialize
$\epsilon^{(1)}=0.99 \lambda^{(1)}\eta^{(1)}$ and then set
$\epsilon^{(k+1)} = c_{\lambda}^2 \epsilon^{(k)}$ %  as the target
% optimality
% error for the next subproblem
for all $k\geq 1$.

\section{Numerical experiments}
\label{ch3_sec:computations}
In our numerical experiments, we focused on problems where both $\mu_1>0$ and $\mu_2>0$.
Two important problems of this form  are the principal component pursuit and stable
principal component pursuit problems, given in \eqref{ch3_eq:component_pursuit} and \eqref{ch3_eq:stable_component_pursuit}, respectively. In the first set
of experiments we solved a set of randomly generated instances of principal
component pursuit problems. In
this setting, we compare FALC with another augmented Lagrangian algorithm
I-ALM~\cite{Ma09_1J}, a proximal gradient algorithm APG~\cite{Ma09_1R} and a soft-thresholding algorithm
SVT~\cite{Can08_2J}. In the second set of experiments, we solved a set of
randomly generated instances of stable principal component pursuit
problem. Since I-ALM, APG and SVT are not able to solve this problem, we
only report statistics for FALC.  In
Section~\ref{ch3_sec:expt-setup}, we describe the methodology we have used in
both experimental settings for generating random problem
instances. All the numerical experiments were conducted
on an IBM Thinkpad laptop with a Intel Core 2 CPU T7200~@2.0 GHz processor, 3GB
SDRAM running MATLAB 7.2 on Windows XP Professional operating system.

The augmented Lagrangian algorithm~FAL introduced in~\cite{Aybat12} is
an implementation of FALC for the basis pursuit problem.
 The numerical results reported in~\cite{Aybat12} show that FAL was 2-7 times faster than the specialized
algorithms NESTA~\cite{Can09_4J},
FPC and FPC-BB~\cite{Yin07_1R,Yin08_1J}, FPC-AS~\cite{Wen09_1R},
YALL1~\cite{Yang09} and SPGL1~\cite{Ber08_1J}. See Tables 6.8, 6.9, 6.10, 6.11,
6.13 in \cite{Aybat12} for details.

The numerical results in this paper and those in~\cite{Aybat12}
clearly show that  FALC is very competitive
with the state-of-the-art algorithms for the special cases of the
composite norm minimization problem.

\subsection{Data generation}
\label{ch3_sec:expt-setup}
We tested FALC on randomly generated data matrices $D=X_0+S_0+\zeta_0$,
where
\begin{enumerate}[i.]
\item $X_0=UV^T$, such that $U\in\reals^{n\times r}$, $V\in\reals^{n\times
    r}$ for $r=0.05 n$ and $U_{ij}\sim \mathcal{N}(0,1)$, $V_{ij}\sim
  \mathcal{N}(0,1)$ for all $i,j$ are independent standard Gaussian
  variables,
\item $\Lambda\subset\{(i,j):\ 1 \leq i,j\leq n\}$ such that cardinality
  of $\Lambda$, $|\Lambda|=p$ for $p=0.05 n^2$,
\item $(S_0)_{ij}\sim\mathcal{U}[-1,1]$ for all $(i,j)\in\Lambda$ are
  independent uniform random variables between $-1$ and $1$,
\item $(\zeta_0)_{ij}\sim\delta\mathcal{U}[-1,1]$ for all $i,j$ are independent
  Gaussian variables.
\end{enumerate}

\subsection{Principal Component Pursuit Problem}
\label{ch3_sec:rPCA_results}
In this section we solve the problem
\begin{equation}
\label{ch3_eq:test_problem1}
  \begin{array}{rl}
    \mbox{min}_{X,S \in \reals^{m\times n}} &\norm{X}_* + \mu_2\norm{\vect(S)}_1,\\
    \mbox{subject to} &  X+S=D,
  \end{array}
\end{equation}
and  report the results of our numerical experiments
comparing FALC with I-ALM~\cite{Ma09_1J}, APG~\cite{Ma09_1R} and
SVT~\cite{Can08_2J}. All the codes for I-ALM, APG and SVT, can be found at
~\url{http://perception.csl.uiuc.edu/matrix-rank/home.html}. Note
that SVT~\cite{Can08_2J} algorithm was originally proposed for solving the
matrix completion problem. The algorithm we used in our numerical study is
an adaptation of the SVT algorithm by % John
Wright and % Shankar
Rao at the
Perception and Decision Laboratory in University of Illinois,
Urbana-Champaign to solve robust PCA problem.

We created 10 random problems of size $n=500$, i.e. $D\in\reals^{500
  \times 500}$ using the procedure
described in Section~\ref{ch3_sec:expt-setup}, where $\delta$ is set to $0$,
i.e. $\zeta_0=\mathbf{0}$. We chose parameter values for
each of the four algorithms so that they produce a  solution $X_{sol}$ and $S_{sol}$
with relative-infeasibility  approximately equal to $5\times 10^{-9}$,
i.e. $\frac{\norm{X_{sol}+S_{sol}-D}_F}{\norm{D}_F} \approx 5\times 10^{-9}$. For each
algorithm   we set the parameters  by solving a set of small size problems
and these parameter values were fixed throughout the experiments, all
other parameters are set to their
default values. The termination criteria are not directly comparable due
to different formulations of the problem solved by different
solvers. For FALC we attempted  to set the stopping parameter
$\varrho$ such that on  average the stopping criterion for FALC is
tighter than the stopping criteria of all the other algorithms we tested.

\begin{enumerate}[1.]
\item\textbf{FALC}: Problem~\eqref{ch3_eq:test_problem1} is a special case of
  problem~\eqref{ch3_eq:composite_norm_minimization} with $\delta=0$, $\alpha=1$ and $\beta=1$. Therefore,
  $f^{(k)}(X,s,y)$ defined in \eqref{ch3_eq:smooth_part} simplifies to
  $f^{(k)}(X,S)=\frac{1}{2}\norm{\vect(X+S)-\vect(D)-\lk\tetki}_2^2$ (note
  that for all $k\geq 1$, $\tetkii=0$). We set $c_\tau=0.4$, $c_\epsilon=0.4$,
  $c_\lambda=0.4$, $\bar{c}_\tau=0.999$, $\bar{c}_\epsilon=0.999$,
  $\bar{c}_\lambda=2$ and initialize $\theta_1^{(1)}$ as in
  \cite{Ma09_1J}, i.e.
    \begin{align}
    \label{ch3_eq:theta_initial}
    \theta_1^{(1)}=\frac{1}{\max\{\norm{{\rm sign}(D)}_2,~\sqrt{n}\norm{\vect({\rm sign}(D))}_\infty\}}
    \vect({\rm sign}(D)).
    \end{align}
    Finally, we set $\varrho=1\times10^{-5}$ and terminate FALC when the
    distance between successive inner iterates are below the threshold
    $\varrho$ for each component, i.e.
    $\norm{X_1^{(k,\ell)}-X_1^{(k,\ell-1)}}_F\leq\varrho$ and
    $\norm{s_1^{(k,\ell)}-s_1^{(k,\ell-1)}}_2\leq\varrho$ for any $k\geq
    1$.  We  used PROPACK~\cite{propack} for computing partial singular
    value decompositions. In order to estimate the rank of $X_0$, we
    followed the scheme proposed in Equation~(17) in \cite{Ma09_1J}. The
    code for PROPACK is available at
    [\url{http://soi.stanford.edu/~rmunk/PROPACK/}].
\item \textbf{I-ALM}: I-ALM solves
  $\min\{\norm{X}_*+\frac{1}{\sqrt{n}}\norm{\vect(S)}_1:\ X+S=D\}$. Let
  $(\Xk,\Sk)$ be the $k$-th iterate. I-ALM terminates when
  $\frac{\norm{\Xk+\Sk-D}_F}{\norm{D}_F}\leq 1\times 10^{-8}$.
\item \textbf{APG}: For some $\bar{\lambda}>0$, APG solves
  $\min\left\{\bar{\lambda}\left(\norm{X}_*+\frac{1}{\sqrt{n}}\norm{\vect(S)}_1\right)
    +\frac{1}{2}\norm{X+S-D}_F^2\right\}$. Stopping
  tolerance is set to $5\times 10^{-11}$ (the definition of stopping
  criteria is complicated, for details see partial APG code at
  [\url{http://perception.csl.uiuc.edu/matrix-rank/home.html}]. In the
  code, by default $\bar{\lambda}$ is set to $\sigma_{\max}(D)\times
  10^{-9}$.
\item \textbf{SVT}: SVT solves a relaxation of the robust PCA problem,
    $$\min\left\{\bar{\lambda}\left(\norm{X}_*+\frac{1}{\sqrt{n}}\norm{\vect(S)}_1\right)+\frac{1}{2}(\norm{X}^2_F+\norm{S}^2_F):\ X+S=D\right\}.$$
Let $(\Xk,\Sk)$ be the $k$-th iterate when $\bar{\lambda}$ is set to $1\times 10^3$. SVT terminates
  $\frac{\norm{\Xk+\Sk-D}_F}{\norm{D}_F}\leq 5\times 10^{-4}$. Note
  that we have chosen a weaker stopping criterion for SVT.
\end{enumerate}
The results of the experiments are displayed in Table~\ref{ch3_tab:comperative_test_results}.
In Table~\ref{ch3_tab:comperative_test_results},
the row labeled $\mathbf{CPU}$ lists the running time of each algorithm in
\emph{seconds} and all other rows are self-explanatory. The column
labeled {\tt average} lists the average taken over the $10$ random
instances, the columns labeled {\tt min} (resp. {\tt max}) list the
minimum (resp. maximum) over the $10$ instances. The experimental results in
Table~\ref{ch3_tab:comperative_test_results},
show that FALC is competitive with the state of the art algorithms,
e.g. I-ALM, APG and SVT, specialized for solving robust PCA problem. Even
though FALC is not special purpose algorithm for robust PCA, % \eqref{ch3_eq:general_problem},
in our numerical experiments, FALC required fewer singular value decompositions
when compared to APG and SVT. In addition, for all 10 randomly created problems in the
test set, only FALC and I-ALM accurately identified
the zero-set of the sparse component $S_0$, i.e. $I_0=\{(ij) \in\{1,2,...,n\}\times \{1,2,...,n\}:
(S_0)_{ij}=0\}$ without any thresholding. This feature of FALC is very
appealing in practice. For signals
with a large dynamic range, almost all of the state-of-the-art efficient
algorithms produce a solution with many small non zeros terms, and it is
often
hard to determine the threshold.

\begin{table}[!htb]
    \centering
    {\footnotesize
    \begin{tabular}{cc|c|c|c|c|c|}
    \cline{2-7}
    &\multicolumn{3}{|c|}{\textbf{FALC}} &\multicolumn{3}{|c|}{\textbf{I-ALM}}\\ \cline{2-7}
    &\multicolumn{1}{|c|}{\textbf{Average}} &\textbf{Min}&\textbf{Max}&\textbf{Average}&\textbf{Min}&\textbf{Max}\\ \hline
    \multicolumn{1}{|c|}{\textbf{svd \#}}
    &40&39&44&31.6&30&33 \\ \hline
    \multicolumn{1}{|c|}{$\mathbf{\norm{X_{sol}-X_0}_F/\norm{X_0}_F}$}
    &3.5E-09&2.7E-09&4.5E-09&1.9E-09&5.9E-10&3.4E-09 \\ \hline
    \multicolumn{1}{|c|}{$\mathbf{\norm{S_{sol}-S_0}_F/\norm{S_0}_F}$}
    &1.3E-07&1.0E-07&1.8E-07&1.9E-07&4.8E-08&3.8E-07 \\ \hline
    \multicolumn{1}{|c|}{$\mathbf{|~\norm{X_{sol}}_*-\norm{X_0}_*|/\norm{X_0}_*}$}
    &1.6E-10&2.4E-11&3.6E-10&1.1E-11&3.7E-12&2.1E-11 \\ \hline
    \multicolumn{1}{|c|}{$\mathbf{\max\{|\sigma_i-\sigma^0_i|: \sigma^0_i>0\}}$}
    &2.1E-07&1.0E-07&4.0E-07&8.7E-08&2.3E-08&2.5E-07 \\ \hline
    \multicolumn{1}{|c|}{$\mathbf{\max\{|\sigma_i|: \sigma^0_i=0\}}$}
    &1.2E-13&7.2E-14&1.9E-13&1.5E-13&5.9E-14&3.7E-13 \\ \hline
    \multicolumn{1}{|c|}{$\mathbf{|~\norm{\vect(S_{sol})}_1-\norm{\vect(X_0)}_1|/\norm{\vect(X_0)}_1}$} &1.4E-08&4.1E-09&2.6E-08&2.2E-09&4.1E-10&5.1E-09 \\ \hline
    \multicolumn{1}{|c|}{$\mathbf{\max\{|(S_{sol})_{ij}-(S_0)_{ij}|: |(S_0)_{ij}|>0\}}$}
    &8.0E-07&5.0E-07&1.4E-06&1.1E-05&2.3E-06&2.5E-05 \\ \hline
    \multicolumn{1}{|c|}{$\mathbf{\max\{|(S_{sol})_{ij}|: (S_0)_{ij}=0\}}$}
    &0&0&0&0&0&0 \\ \hline
    \multicolumn{1}{|c|}{$\mathbf{rank}$}
    &25&25&25&25&25&25 \\ \hline
    \multicolumn{1}{|c|}{$\mathbf{\norm{X_{sol}+S_{sol}-D}_F/\norm{D}_F}$}
    &3.5E-09&2.6E-09&4.5E-09&4.7E-09&1.1E-09&9.6E-09 \\ \hline
    \multicolumn{1}{|c|}{$\mathbf{CPU}$}
    &22.9&19.6&27.8&16.8&13.5&24.3 \\ \hline\\
    \cline{2-7}
    &\multicolumn{3}{|c|}{\textbf{APG}} &\multicolumn{3}{|c|}{\textbf{SVT}}\\ \cline{2-7}
    &\multicolumn{1}{|c|}{\textbf{Average}}&\textbf{Min}&\textbf{Max}&\textbf{Average}&\textbf{Min}&\textbf{Max}\\ \hline
    \multicolumn{1}{|c|}{\textbf{svd \#}}
    &187.7&187&188&833.9&819&857 \\ \hline
    \multicolumn{1}{|c|}{$\mathbf{\norm{X_{sol}-X_0}_F/\norm{X_0}_F}$}
    &4.1E-09&4.0E-09&4.4E-09&1.8E-04&1.8E-04&1.8E-04 \\ \hline
    \multicolumn{1}{|c|}{$\mathbf{\norm{S_{sol}-S_0}_F/\norm{S_0}_F}$}
    &1.6E-07&1.6E-07&1.7E-07&2.0E-02&2.0E-02&2.1E-02 \\ \hline
    \multicolumn{1}{|c|}{$\mathbf{|~\norm{X_{sol}}_*-\norm{X_0}_*|/\norm{X_0}_*}$}
    &4.0E-09&3.8E-09&4.2E-09&1.7E-05&1.5E-05&1.9E-05 \\ \hline
    \multicolumn{1}{|c|}{$\mathbf{\max\{|\sigma_i-\sigma^0_i|: \sigma^0_i>0\}}$}
    &2.0E-06&1.9E-06&2.1E-06&1.5E-02&1.2E-02&1.7E-02 \\ \hline
    \multicolumn{1}{|c|}{$\mathbf{\max\{|\sigma_i|: \sigma^0_i=0\}}$}
    &1.3E-13&6.8E-14&1.9E-13&2.4E-13&7.6E-14&6.8E-13 \\ \hline
    \multicolumn{1}{|c|}{$\mathbf{|~\norm{\vect(S_{sol})}_1-\norm{\vect(X_0)}_1|/\norm{\vect(X_0)}_1}$} &1.8E-07&1.8E-07&1.9E-07&5.0E-03&4.9E-03&5.1E-03 \\ \hline
    \multicolumn{1}{|c|}{$\mathbf{\max\{|(S_{sol})_{ij}-(S_0)_{ij}|: |(S_0)_{ij}|>0\}}$} &2.0E-07&1.8E-07&2.3E-07&1.2E-01&1.1E-01&1.3E-01 \\ \hline
    \multicolumn{1}{|c|}{$\mathbf{\max\{|(S_{sol})_{ij}|: (S_0)_{ij}=0\}}$}
    &3.7E-08&2.1E-08&6.6E-08&5.5E-03&3.6E-03&8.5E-03  \\ \hline
    \multicolumn{1}{|c|}{$\mathbf{rank}$}
    &25&25&25&25&25&25 \\ \hline
    \multicolumn{1}{|c|}{$\mathbf{\norm{X_{sol}+S_{sol}-D}_F/\norm{D}_F}$}
    &5.4E-09&5.2E-09&5.8E-09&5.0E-04&5.0E-04&5.0E-04 \\ \hline
    \multicolumn{1}{|c|}{$\mathbf{CPU}$}
    &87.7&71.6&101.6&265.2&252.0&273.1 \\ \hline
    \end{tabular}
    }
    \caption{FALC vs I-ALM, APG, SVT: Numerical Test Results for PCP problem with $n=500$, $r=0.05n^2$, $p=0.05n$}
    \label{ch3_tab:comperative_test_results}
    \vspace{-.25in}
\end{table}
\begin{table}[!htb]
    \centering
    {\footnotesize
    \begin{tabular}{c|c|c|c|}
    \cline{2-4}
    &\multicolumn{3}{|c|}{\textbf{FALC}}\\ \cline{2-4}
    &\textbf{Average}  &\textbf{Min}&\textbf{Max}\\ \hline
    \multicolumn{1}{|c|}{\textbf{svd \#}}
    &59.3&55&64 \\ \hline
    \multicolumn{1}{|c|}{$\mathbf{\norm{X_{sol}-X_0}_F/\norm{X_0}_F}$}
    &1.7E-05&1.7E-05&1.7E-05 \\ \hline
    \multicolumn{1}{|c|}{$\mathbf{\norm{S_{sol}-S_0}_F/\norm{S_0}_F}$}
    &3.7E-04&3.0E-04&4.4E-04 \\ \hline
    \multicolumn{1}{|c|}{$\mathbf{|~\norm{X_{sol}}_*-\norm{X_0}_*|/\norm{X_0}_*}$}
    &1.6E-05&1.6E-05&1.6E-05 \\ \hline
    \multicolumn{1}{|c|}{$\mathbf{\max\{|\sigma_i-\sigma^0_i|: \sigma^0_i>0\}}$}
    &9.9E-03&9.7E-03&1.0E-02 \\ \hline
    \multicolumn{1}{|c|}{$\mathbf{\max\{|\sigma_i|: \sigma^0_i=0\}}$}
    &1.6E-13&3.6E-14&3.1E-13 \\ \hline
    \multicolumn{1}{|c|}{$\mathbf{|~\norm{\vect(S_{sol})}_1-\norm{\vect(X_0)}_1|/\norm{\vect(X_0)}_1}$}
    &2.3E-04&2.2E-04&2.4E-04 \\ \hline
    \multicolumn{1}{|c|}{$\mathbf{\max\{|(S_{sol})_{ij}-(S_0)_{ij}|: |(S_0)_{ij}|>0\}}$}
    &3.9E-03&3.0E-03&4.6E-03 \\ \hline
    \multicolumn{1}{|c|}{$\mathbf{\max\{|(S_{sol})_{ij}|: (S_0)_{ij}=0\}}$}
    &6.4E-05&0.0E+00&2.3E-04 \\ \hline
    \multicolumn{1}{|c|}{$\mathbf{rank}$}
    &25&25&25 \\ \hline
    \multicolumn{1}{|c|}{$\mathbf{\norm{X_{sol}+S_{sol}-D}_F/\norm{D}_F}$}
    &2.1E-05&2.0E-05&2.2E-05 \\ \hline
    \multicolumn{1}{|c|}{$\mathbf{CPU}$}
    &34.6&26.1&48.3 \\ \hline
    \end{tabular}
    }
    \caption{FALC: Numerical Test Results for SPCP problem with $n=500$, $r=0.05n^2$, $p=0.05n$, $\delta=1\times 10^{-4}$}
    \label{ch3_tab:sRPCA_test_results}
    \vspace{-.25in}
\end{table}
\subsection{Stable Principal Component Pursuit Problem}
\label{ch3_sec:test}
In this section, we solve the problem
\begin{equation}
\label{ch3_eq:test_problem2}
  \begin{array}{rl}
    \mbox{min}_{X,S \in \reals^{m\times n}} &\norm{X}_* + \mu_2\norm{\vect(S)}_1,\\
    \mbox{subject to} &  \norm{\vect(X+S-D)}_\infty\leq\delta,
  \end{array}
\end{equation}
and report the results of our numerical experiments using FALC. To best of
our knowledge, there are no publicly available code specialized for
solving problem in \eqref{ch3_eq:test_problem2}, other than general
purpose SDP solvers.

We created 10 random problems of size $n=500$, i.e. $D\in\reals^{500
  \times 500}$ using the procedure described in
Section~\ref{ch3_sec:expt-setup}, where $\delta$ is set to $1\times 10^4$,
i.e. each entry of the noise term $\zeta_0$ is coming from a uniform
distribution between $[-\delta, \delta]$. We chose the value of the
stopping parameter so that FALC produces a solution $X_{sol}$ and
$S_{sol}$ with $\frac{\norm{X_{sol}+S_{sol}-D}_F}{\norm{D}_F} \approx
1\times 10^{-5}$.

Problem in \eqref{ch3_eq:test_problem2} is a special case of
\eqref{ch3_eq:composite_norm_minimization} and $f^{(k)}(X,s,y)$ defined in
\eqref{ch3_eq:smooth_part} simplifies to $f^{(k)}(X,S,y)$
$=\frac{1}{2}\norm{\vect(X+S)-y-\vect(D)-\lk\tetki}_2^2$ (note that for
all $k\geq 1$, $\tetkii=0$). We set the parameter values for FALC by
solving a set of small size problems and these parameter values were fixed
throughout the experiments, all other parameters are set to their default
values, i.e. $c_\tau=0.4$, $c_\epsilon=0.4$, $c_\xi=0.4$, $c_\lambda=0.4$,
$\bar{c}_\tau=0.999$, $\bar{c}_\epsilon=0.999$, $\bar{c}_\xi=0.999$. We
set $\bar{c}_\lambda=1.5$ and initialize $\theta_1^{(1)}$ as in
\cite{Ma09_1J}, i.e. as in \eqref{ch3_eq:theta_initial}.

Finally, We set $\varrho=1\times10^{-5}$, $\varsigma=1\times 10^{-3}$ and
terminate FALC when either the distance between successive inner iterates
are below a threshold $\varrho$ for each component,
i.e. $\norm{\vect\left(X_1^{(k,\ell)}\right)-\vect\left(X_1^{(k,\ell-1)}\right)}_\infty\leq\varrho$,
$\norm{\vect\left(s_1^{(k,\ell)}\right)-\vect\left(s_1^{(k,\ell-1)}\right)}_\infty\leq\varrho$
for any $k\geq 1$ or there exist partial subgradients with sufficiently
small norm for each component, i.e.
$$\norm{G}_F\leq\varsigma/2,\ \norm{g}_2\leq\varsigma\ \hbox{ for some }\
(G,g)\in\partial_{X,s} P^{(k)}(.,.,.)|_{\breve{Z}^{(k,\ell)}}$$
and $$\norm{\breve{y}_1^{(k,\ell)}-\Pi_\cQ\left(\breve{y}_1^{(k,\ell)}-\frac{1}{L}\grad_y
    P^{(k)}(\breve{Z}^{(k,\ell)})\right)}_2\leq\varsigma,$$
where
$\breve{Z}^{(k,\ell)}=\left(\breve{X}_1^{(k,\ell)},
  \breve{s}_1^{(k,\ell)},\breve{y}_1^{(k,\ell)}\right)$
is defined at the beginning of Section~\ref{ch3_sec:implementation}.

We have used PROPACK~\cite{propack} for computing partial singular value
decompositions. In order to estimate the rank of $X_0$, we followed the
scheme proposed in Equation~(17) in \cite{Ma09_1J}. The results of the
experiments are displayed in Table~\ref{ch3_tab:sRPCA_test_results}.

\section{Extension to general composite norm problem}
\label{ch3_sec:extensions}

The algorithmic framework proposed in this paper extends to the
following much more general class of problems given in
\eqref{ch3_eq:general_problem}. % Let $\mu_1>0$ and $\mu_2>0$.
By introducing slack variables, \eqref{ch3_eq:general_problem} can be
reformulated as follows.
\begin{equation}
  \label{ch3_eq:CN-ext-reform}
  \begin{array}{rll}
    \min & \mu_1 \norm{\sigma(S)}_{\alpha}+\mu_2 \norm{s}_{\beta}+\mu_3 H(X),\\
    \mbox{subject to} &
    \begin{array}[t]{lll}
      \cF(X) - S = G, \\
      \cC(X) - s = d, \\
      \cA(X) - y = b, & y\in\cQ,
   \end{array}
 \end{array}
\end{equation}
where the decision variables $X\in\reals^{m\times n}$, $S\in\reals^{r_1\times
  r_2}$, $s\in\reals^p$, and  $y\in\reals^q$.
% Here the affine mapping $\cF(.)-G:\reals^{m\times
%   n}\rightarrow\reals^{r\times r}$ maps $X$ to a symmetric matrix and
$H(.)$ is a strongly convex function with convexity parameter $\varsigma$. We continue to assume that
$\cA$ is surjective; however,  when
$\mu_3 > 0$, we no longer require that at least one of
that at least one of $\cF$ and $\cC$ is an injective linear map

In this more general setting, the FALC inexactly solves subproblems of the form:
\begin{equation}
  \label{ch3_eq:CN-ext-dual}
  \min_{X,~S,~s,~y\in\cQ}\Pk(X,S,s,y),
\end{equation}
where
\begin{eqnarray*}
\Pk(X,S,s,y)&:=&\lk\big(\mu_1\norm{\sigma(S)}_\alpha+\mu_2\norm{s}_\beta+\mu_3 H(X)\big)+f^{(k)}(X,S,s,y),\\
f^{(k)}(X,S,s,y)&:=&
    \begin{array}[c]{l}
      \frac{1}{2}\norm{\mathcal{F}(X)-S-G-\lk\theta_3^{(k)}}_F^2
      +\frac{1}{2}\norm{\mathcal{C}(X)-s-d-\lk\theta_2^{(k)}}_2^2
      \\
      \mbox{} +\frac{1}{2}\norm{\mathcal{A}(X)-y-b-\lk\theta_1^{(k)}}_2^2.
    \end{array}
\end{eqnarray*}
%Note that we do {\em not} dualize neither the norm constraint $\norm{y}
%\leq \delta$ nor the cone constraints: $\norm{v}_2\leq t$ and $S\succeq
%0$.
Let $(\Xkopt,S_*^{(k)},\skopt,\ykopt)\in\argmin\Pk(X,S,s,y)$. Suppose the
initial iterate $X^{(0)}$ is feasible, i.e. $\cA(X^{(0)})-b\in\cQ$. Let
$S^{(0)}:=\cF(X^{(0)})-G$, $s^{(0)}:=\cC(X^{(0)})-d$,
$y^{(0)}:=\cA(X^{(0)})-b$, and
$\eta:=\mu_1\norm{\sigma(\cF(X^{(0)})-G)}_\alpha+\mu_2\norm{\cC(X^{(0)})-d}_\beta$.

The particular implementation of FALC depends on the nature of the
objective function. In all cases we
need to ensure that the iterate  sequence $\{\Xk\}_{k\in\integers_+}$  is bounded so
that it has a limit point.
%Suppose that we are solving the composite norm minimization problem with
%just $\cD(X)=c$ constraint. If the boundedness condition on
%$\{\Xk\}_{k\in\integers_+}$ is not satisfied, then
%$\{\Xk\}_{k\in\integers_+}$ sequence can escape to infinity in the affine
%space defined by the null space of $\cD$.
First consider the case where $\mu_3 > 0$. Strong convexity property of
$H(.)$ implies that
\begin{eqnarray*}
&\mu_1\norm{\sigma(S_*^{(k)})}_\alpha+\mu_2\norm{\skopt}_\beta
+\frac{\varsigma}{2}\norm{\Xkopt-\left(X^{(0)}-\frac{\grad
      H(X^{(0)})}{\varsigma}\right)}_F^2\\
&\mbox{}\leq\eta +\frac{1}{2\varsigma} \norm{\grad
  H(X^{(0)})}_F^2+\frac{\lk}{2}\sum_{i=1}^3\norm{\theta_i^{(k)}}_F^2.
\end{eqnarray*}
Therefore, we can define
$\Sk$ in line~\ref{ch3_alg:Fk} % of \textbf{Algorithm~FALC}
% displayed
in Figure~\ref{ch3_alg:falc} as follows:
\eq
\Sk:=\left\{(X,S,s,y):
  \mu_1\norm{\sigma(S)}_\alpha+\mu_2\norm{s}_\beta\leq\etak,~
  \left\|X-\left(X^{(0)}-\frac{\grad
        H(X^{(0)})}{\varsigma}\right)\right\|_F\leq
  \sqrt{\frac{2}{\varsigma}\etak}\right\},
\en
where
\begin{equation}
\label{ch3_eq:etak}
\etak:=\eta+\frac{1}{2\varsigma}\norm{\grad
  H(X^{(0)})}_F^2+\frac{\lk}{2}\sum_{i=1}^3\norm{\theta_i^{(k)}}_F^2.
\end{equation}
% The implementable algorithm almost does not change.
The only change in the algorithm is that we need to compute $\grad H$ at
every iteration of \textbf{Algorithm~APG}
additional to one $\grad f$ computation per iteration. %The case with $\mu_3$ was considered in the previous sections.

When $\mu_3 = 0$, we set
\eq
\Sk:=\{(X,S,s,y):~\mu_1\norm{\sigma(S)}_\alpha+\mu_2\norm{s}_\beta\leq\etak\}.
\en
ensuring that the iterates $\{(S^{(k)},s^{(k)},\yk)\}_{k\in\integers_+}$
are bounded, see Lemma~\ref{lem:bounded} for details. Since at
least one of $\cF$ and $\cC$ is injective, this implies that  $\{\Xk\}_{k\in\integers_+}$
is a bounded sequence. Note that without the injectivity assumption, the
sequence
$\{\Xk\}_{k\in\integers_+}$ may not have a limit point.

The general formulation \eqref{ch3_eq:general_problem} subsumes a number
of different problems as special cases -- see
Section~\ref{sec:special-cases} for details.  And,  our experience with
FALC leads us to believe that this algorithm
is likely to be very competitive for solving all these special
cases. However, the assumption that the operator $\cA$, defining the
constraints, be surjective can be restrictive. In some applications, the feasible region is the intersection of cones of the form: $\{X: \cA_1(X) -
b_1 \in \cQ_1, \cA_2(X) - b_2
\in \cQ_2\}$. While, it is often the case that each $\cA_i$, $i = 1, 2$,
is surjective, the product operator $\cA(X) = [\cA_1(X), \cA_2(X)]$ is
not. Thus, FALC cannot be used for these problems. The extension to
intersection of cones is  non-trivial and one would have to
design a completely new set of techniques.

% Define
% \begin{equation}
% \label{ch3_eq:etak}
% \etak:=\eta+\frac{1}{2\varsigma}\norm{\grad
% F(X^{(0)})}_F^2+\frac{\lk}{2}\sum_{i=1}^3\norm{\theta_i^{(k)}}_F^2.
% \end{equation}
% We redefine

% Next, suppose that $\mu_3=0$.

% FALC solves \eqref{ch3_eq:CN-ext-dual} by solving constrained shrinkage
% problems of the following form.
% \begin{enumerate}
% \item Matrix optimization problem over simple sets:
%   \begin{align}
%   &\min_S\Big\{\lambda\norm{\sigma(S)}_{\alpha}+\frac{1}{2}\norm{S-\tilde{S}}_F^2:\;\norm{\sigma(S)}_{\alpha}\leq\eta\},\\
%   &\min_{X}\Big\{\frac{1}{2}\norm{X-\tilde{X}}_F^2\Big\}. \label{ch3_eq:ineq_matrix_subproblem}
%   \end{align}
%   For given $\tilde{X}\in\reals^{m\times n}$ and $\tilde{S}\in\reals^{m\times n}$, $\lambda>0$ and $\eta>0$, these problems can be efficiently solved when $\norm{\sigma(.)}_\alpha$ is the either the the nuclear norm, Frobenius norm, or the $\ell_2$-norm.

% \item Vector optimization problem over simple sets:
%   \begin{align}
%   &\min_{s\in\reals^p}\Big\{\lambda\norm{s}_{\beta}+\frac{1}{2}\norm{s-\tilde{s}}_2^2:\;\norm{s}_{\beta}\leq\eta\Big\},\\
%   &\min_{y\in\reals^q}\Big\{\frac{1}{2}\norm{y-\tilde{y}}_2^2:y\in\cQ\Big\}, \label{ch3_eq: ineq_norm_vector_subproblem}.
%   \end{align}
%   For given $\tilde{s}$, $\tilde{y}$ and $\lambda>0$, these problems can be efficiently solved when $\beta$ is either $\ell_2$, $\ell_1$ or $\ell_{\infty}$ vector norms.
% \end{enumerate}
% The extension~\eqref{ch3_eq:general_problem} allows us to model a wider
% class of problems discussed in the introduction.

The main contribution of this paper is an efficient first-order augmented
Lagrangian algorithm~(FALC) for the composite norm minimization problem
\eqref{ch3_eq:composite_norm_minimization} and for its extension
\eqref{ch3_eq:general_problem}. FALC solves the composite norm
minimization problem
by solving a sequence of augmented Lagrangian subproblems, where
each subproblem is solved using \textbf{Algorithm~APG} in
Figure~\ref{ch1_alg:pga}. \textbf{Algorithm~APG} is essentially Algorithm~2
in~\cite{Tseng08} (see also FISTA~\cite{Beck09_1J}) with early
termination.
 %with $x$-update rule \eqref{ch1_eq:fast_x}.
We show that the continuation scheme on penalty parameter $\lambda$ used
in FALC guarantees that the iterate sequence provably converges to the
solution and we are also able to compute a convergence rate.
% We found that for a fixed measurement ratio  $p/(mn)$, sparsity ratio $q/(mn)$,
% and solution accuracy parameter $\gamma$,  the total number of inner
% iterations is effectively independent of the dimension $mn$ of the
% target signal; thus, one can tune the parameters on the smallest problem
% and use these parameters for all larger problems.
The performance of FALC in our numerical experiments has been very
promising. To best of our knowledge, for the stable PCA problem, FALC is
the first algorithm with a known complexity bound.

\bibliographystyle{siam}
\bibliography{thesis}
\appendix

\section{Proofs of technical results}
\label{app:proofs}

\subsection*{Lemma~\ref{lem:y-bound}  and proof}
$\mbox{}$
\begin{lemma}
\label{lem:y-bound}
Let $\cQ\subset\reals^q$ be nonempty closed convex set such that
$\{X\in\reals^{m\times n}: \cA(X)-b\in\cQ\}\neq\emptyset$, where $\cA$
is surjective; and let $(\Xkopt,\skopt,\ykopt)$ is an optimal solution
to \eqref{ch3_eq:augmented_lagrangian_subproblem_2}. Then, for all
$k\geq 1$,
\begin{align}
\norm{\ykopt}_2\leq \sigma_{\max}(A)\norm{\Xkopt}_F
+\norm{b+\lk\tetki}_2+ 2~\min_{\tilde{y} \in \cQ}\{\norm{\tilde{y}}_2\}.   \label{eq:y_bound}
\end{align}
% Hence the upper bound $\etakii$ on $\norm{\ykopt}_2$ can be chosen as follows:
% \begin{align}
% \etakii:=\sigma_{\max}(A)\frac{I(\alpha)}{\mu_1}~\etak
% +\norm{b+\lk\tetki}_F+2\norm{\cA(X^{(0)})-b}_2. \label{eq:eta2}
% \end{align}
% Suppose $0\in\cQ$. Then we can choose
% , then $X^{(0)}$ can be chosen such that
% $\norm{\cA(X^{(0)})-b}_F$ in \eqref{eq:y_bound} and in \eqref{eq:eta2}
% is $0$.
\end{lemma}
\begin{proof}
%The first statement is trivial since $\ykopt\in\cQ$. In the rest of
%the proof, $\cQ$ is not assumed to be bounded.
From the first order optimality conditions for
\eqref{ch3_eq:augmented_lagrangian_subproblem_2}, we have
$\ykopt=\Pi_\cQ(\cA(\Xkopt)-b-\lk\tetki)$. Since Euclidean
projection is nonexpansive, we have
\begin{align}
\norm{\ykopt-\tilde{y}}_2\leq \norm{\cA(\Xkopt)-b
  -\lk\tetki-\tilde{y}}_2 \quad \forall
\tilde{y}\in\cQ. \label{eq:nonexpansive}
\end{align}
% Since $\{X\in\reals^{m\times n}: \cA(X)-b\in\cQ\}\neq\emptyset$, there
% exists
The result now follows from the triangular inequality.
\end{proof}

This result implies several simple bounds on $\norm{\ykopt}_2$.
Since the initial iterate $X^{(0)}$ is feasible, i.e. $\cA(X^{(0)}) -
b \in \cQ$, it follows that
\begin{equation}
\label{eq:eta2}
\norm{\ykopt}_2\leq \etakii := \sigma_{\max}(A)\norm{\Xkopt}_F
+\norm{b+\lk\tetki}_2+ 2\norm{\cA(X^{(0)}) - b}_2.
\end{equation}
Suppose $0 \in \cQ$. Then
$
\norm{\ykopt}_2\leq \etakii := \sigma_{\max}(A)\norm{\Xkopt}_F
+\norm{b+\lk\tetki}_2.
$
When $\cQ$ is bounded with $\cQ \subseteq\{y:
\norm{y}_2 \leq \eta_2\}$. Then, one can set $\etakii:=\eta_2$ for all
$k\geq 1$.

%The above lemma is true for all no nempty, closed convex sets. However, in many applications $\cQ$ has further properties, e.g., $\cQ$ is a cone or a bounded set. Below we refine the results depending on the further properties of $\cQ$. Hence, $\etakii$ can be chosen accordingly.
% \end{remark}
\subsection*{Lemma~\ref{ch3_cor:grad_norm_bound} and proof} $\mbox{}$

\begin{lemma}
\label{ch3_cor:grad_norm_bound}
Fix $\alpha$, $\beta \in \{1,2,\infty\}$. Let
\eq
P(X,s,y)=\lambda (\mu_1\norm{\sigma(X)}_\alpha+\mu_2\norm{s}_\beta) +
f(X,s,y)
\en
where
\eq
f(X,s,y) = \frac{1}{2}
\norm{\mathcal{A}(X)-y-b-\lambda\theta_1}_2^2+\frac{1}{2}
\norm{\mathcal{C}(X)-s-d-\lambda\theta_2}_2^2.
\en
Suppose $(\bar{X},\bar{s},\bar{y})$ is $\epsilon$-optimal for the problem $\min_{X,s,y}\{P(X,s,y):~y\in\cQ\}$, i.e.
\eq
0\leq P(\bar{X},\bar{s},\bar{y})- \min_{X \in \reals^{m\times n},~s\in\reals^p,~y\in\cQ\subset\reals^q}P(X,s,y) \leq \epsilon.
\en
Then we have
\begin{eqnarray*}
\norm{\mathcal{C}(\bar{X})-\bar{s}-d-\lambda\theta_2}_2 & \leq &
J(\betas)\mu_2\lambda+\sigma_{max}(M)\sqrt{2\epsilon},\\
\norm{\mathcal{A}^*\left(\mathcal{A}(\bar{X})-\bar{y}-b-\lambda\theta_1\right)+\mathcal{C}^*\left(\mathcal{C}(\bar{X})-\bar{s}-d-\lambda\theta_2\right)}_F &\leq & I(\alphas)\mu_1\lambda+ \sigma_{max}(M)\sqrt{2\epsilon},
\end{eqnarray*}
where
$
M = \left(
     \begin{array}{ccc}
       -I & 0 & C \\
       0 & -I & A \\
     \end{array}
    \right),
$
$\frac{1}{\alphas}+\frac{1}{\alpha}=1$ (resp. $\frac{1}{\betas}+\frac{1}{\beta}=1$) is the H{\"o}lder conjugate of $\alpha$
(resp. $\beta$) and the functions $I(\cdot)$ and $J(\cdot)$ are defined in \eqref{ch3_eq:IJ_beta}.
\end{lemma}

\noindent In order to prove for Lemma~\ref{ch3_cor:grad_norm_bound}, we need the
following result.
\begin{theorem}
\label{ch3_thm:grad_convergence}
Let  $f:\reals^{m\times n}\times \reals^p \times \reals^q \rightarrow\reals$ denote a convex function with  a Lipschitz continuous gradient $\grad f$ with a
Lipschitz constant $L$ with respect to the norm %  $\norm{.}$ on $\reals^{m\times n}\times \reals^p \times \reals^q$ defined as follows:
% \begin{align}
$\norm{(X,s,y)}=\sqrt{\norm{X}_F^2+\norm{s}_2^2+\norm{y}_2^2}$. % \label{ch3_eq:composite_norm}
% \end{align}
Let $(X_*,s_*,y_*) \in \argmin_{X,s,y}\{\lambda(\mu_1\norm{\sigma(X)}_\alpha+\mu_2\norm{s}_\beta)+f(X,s,y): y\in\cQ\}$. Suppose $(\bar{X},\bar{s},\bar{y}) \in\reals^{m\times n} \times \reals^p\times \reals^q$ such that $\bar{y}\in\cQ$ satisfies
\eq
\lambda\big(\mu_1\norm{\sigma(\bar{X})}_\alpha+\mu_2\norm{\bar{s}}_\beta\big)+f(\bar{X},\bar{s},\bar{y})
\leq \lambda\big(\mu_1\norm{\sigma(X_*)}_\alpha+\mu_2\norm{s_*}_\beta\big)+f(X_*,s_*,y_*)+
\epsilon
\en
for some $\epsilon>0$. Then
\eq
    \norm{\grad_X f(\bar{X}, \bar{s}, \bar{y})}_F \leq \big(\sqrt{2L\epsilon}+I(\alphas)  \lambda \mu_1\big),
    \qquad
    \norm{\grad_s f(\bar{X}, \bar{s}, \bar{y})}_2 \leq \big(\sqrt{2L\epsilon}+ J(\betas)\lambda \mu_2\big)
\en
where $\frac{1}{\alphas}+\frac{1}{\alpha}=1$ (resp. $\frac{1}{\betas}+\frac{1}{\beta}=1$) is the H{\"o}lder conjugate of $\alpha$
(resp. $\beta$) and the functions $I(\cdot)$ and $J(\cdot)$ are defined in \eqref{ch3_eq:IJ_def}.
%\begin{enumerate}[(a)]
%\item $\norm{\grad_X f(\bar{X}, \bar{s})}_2 \leq \sqrt{2L\epsilon}+\lambda \mu_1$ when $\norm{.}_{\alpha}$ denotes the nuclear norm; and  $\norm{\grad_X f(\bar{X}, \bar{s})}_F \leq \sqrt{2L\epsilon}+\lambda \mu_1$, when $\norm{.}_{\alpha}$ is either the Frobenius or the $\ell_2$-operator norms.
%\item $\norm{\grad_s f(\bar{X}, \bar{s})}_\infty \leq \sqrt{2L\epsilon}+\lambda \mu_2$ when $\norm{.}_{\beta}$ denotes the $\ell_1$ norm; and $\norm{\grad_s f(\bar{X}, \bar{s})}_2 \leq \sqrt{2L\epsilon}+\lambda \mu_2$ when $\norm{.}_{\beta}$ is either the $\ell_2$ or $\ell_{\infty}$ norms.
%\end{enumerate}
\end{theorem}
\begin{proof}
% Given $\epsilon>0$, fix $\bar{X}\in\reals^{m\times n}$ and $\bar{s}\in\reals^p$ such that
% $\lambda(\mu_1\norm{\bar{X}}_\alpha+\mu_2\norm{\bar{s}}_\beta)+f(\bar{X},\bar{s})-[\lambda(\mu_1\norm{X^*}_\alpha+\mu_2\norm{s^*}_\beta)+f(X^*, s^*)]\leq\epsilon$.
Since
$\grad f$ is Lipschitz continuous with constant $L$, the triangular
inequality for $\norm{\sigma(.)}_\alpha$ and $\norm{.}_\beta$ implies
that for any $X\in\reals^{m\times n}$, $s\in\reals^p$ and
$y\in\reals^q$
\begin{eqnarray*}
\lefteqn{\lambda(\mu_1\norm{\sigma(X)}_\alpha+\mu_2\norm{s}_\beta)+f(X,s,y)}&&
\nonumber \\
& \leq &
\lambda(\mu_1\norm{\sigma(\bar{X})}_\alpha+\mu_2\norm{\bar{s}}_\beta)+
f(\bar{X},\bar{s},\bar{y})+\lambda(\mu_1\norm{\sigma(X-\bar{X})}_\alpha
+\mu_2\norm{s-\bar{s}}_\beta)\nonumber
\\
&& \mbox{} + \fprod{\grad_X
  f(\bar{X},\bar{s},\bar{y}),(X-\bar{X})}+\grad_s f(\bar{X},
\bar{s},\bar{y})^T (s-\bar{s})+\grad_y f(\bar{X}, \bar{s},\bar{y})^T
(y-\bar{y})\nonumber\\
&&\mbox{}+\frac{L}{2}\norm{X-\bar{X}}_F^2
+\frac{L}{2}\norm{s-\bar{s}}_2^2+\frac{L}{2}\norm{y-\bar{y}}_2^2,
\nonumber
\end{eqnarray*}
where $\fprod{X,Y} =\Tr(X^T Y)\in\reals$ denotes the usual Euclidean
inner product of $X\in\reals^{m\times n}$ and $Y\in\reals^{m\times
  n}$. Since $X$, $s$ and $y$ are arbitrary, it follows that
\begin{eqnarray}
\lefteqn{\lambda(\mu_1\norm{\sigma(X_*)}_\alpha+\mu_2\norm{s_*}_\beta)+f(X_*,s_*,y_*)}&& \nonumber\\
&\leq&\lambda(\mu_1\norm{\sigma(\bar{X})}_\alpha+\mu_2\norm{\bar{s}}_\beta)+f(\bar{X},\bar{s},\bar{y}) \nonumber \\
&&\mbox{}+\min_{X\in\reals^{m\times n}}\left\{\fprod{\grad_X f(\bar{X},\bar{s},\bar{y}),X-\bar{X}}+\frac{L}{2}\norm{X-\bar{X}}_F^2+\lambda\mu_1\norm{\sigma(X-\bar{X})}_\alpha\right\} \nonumber\\
&&\mbox{}+\min_{s\in\reals^p}\left\{\grad_s f(\bar{X},\bar{s},\bar{y})^T(s-\bar{s})+\frac{L}{2}\norm{s-\bar{s}}_2^2+\lambda\mu_2\norm{s-\bar{s}}_\beta\right\} \nonumber\\
&&\mbox{}+\min_{y\in\cQ\subset\reals^q}\left\{\grad_y f(\bar{X},\bar{s},\bar{y})^T(y-\bar{y})+\frac{L}{2}\norm{y-\bar{y}}_2^2\right\}.
\label{ch3_eq:combined_upper_bound}
\end{eqnarray}
The first minimization problem on the right hand side of
\eqref{ch3_eq:combined_upper_bound} can be simplified as follows:
\begin{align}
&\min_{X\in\reals^{m\times n}}\left\{\fprod{\grad_X
    f(\bar{X},\bar{s},\bar{y}),X-\bar{X}}+\frac{L}{2}\norm{X-\bar{X}}_F^2
  + \lambda\mu_1\norm{\sigma(X-\bar{X})}_\alpha\right\} \nonumber\\
=&\max_{W:\norm{\sigma(W)}_{\alphas}\leq
\lambda\mu_1}\min_{X\in\reals^{m\times
  n}}\left\{\frac{L}{2}\norm{X-\bar{X}}_F^2+\fprod{\grad_X
  f(\bar{X},\bar{s},\bar{y})+
  W,~X-\bar{X}}\right\}, \label{ch3_eq:minimax_problem_alpha}\\
=&\max_{W:\norm{\sigma(W)}_{\alphas}\leq
  \lambda\mu_1}\left\{\frac{L}{2}\norm{X^*(W)-\bar{X}}_F^2+\fprod{\grad_X
    f(\bar{X},\bar{s},\bar{y})+W,~X^*(W)-\bar{X}}\right\}, \nonumber\\
=&-\min_{W:\norm{\sigma(W)}_{\alphas}\leq \lambda\mu_1}\frac{\norm{\grad_X
    f(\bar{X},\bar{s},\bar{y})+W}_F^2}{2L}, \label{ch3_eq:simple_minimax_problem_alpha}
\end{align}
$X^*(W)=\bar{X}-\frac{\grad_X f(\bar{X},\bar{s},\bar{y})+W}{L}$ is the minimizer of the inner minimization problem in \eqref{ch3_eq:minimax_problem_alpha}.

The second minimization problem on the right hand side of \eqref{ch3_eq:combined_upper_bound} can be simplified as follows:
\begin{align}
&\min_{s\in\reals^p}\left\{\grad_s f(\bar{X},\bar{s},\bar{y})^T(s-\bar{s})+\frac{L}{2}\norm{s-\bar{s}}_2^2+\lambda\mu_2\norm{s-\bar{s}}_\beta\right\}\nonumber\\
=&\max_{u:\norm{u}_{\betas}\leq \lambda\mu_2}\min_{s\in\reals^p}\left\{\frac{L}{2}\norm{s-\bar{s}}_2^2 +(\grad_s f(\bar{X},\bar{s},\bar{y})+ u)^T(s-\bar{s})\right\}, \label{ch3_eq:minimax_problem}\\
=&\max_{u:\norm{u}_{\betas}\leq\lambda\mu_2}\left\{\frac{L}{2}\norm{s^*(u)-\bar{s}}_2^2 +(\grad_sf(\bar{X},\bar{s},\bar{y})+u)^T(s^*(u)-\bar{s})\right\}, \nonumber\\
=&-\min_{u:\norm{u}_{\betas}\leq \lambda\mu_2}\frac{\norm{\grad_s
    f(\bar{X},\bar{s},\bar{y})+u}_2^2}{2L}, \label{ch3_eq:simple_minimax_problem_beta}
\end{align}
$s^*(u)=\bar{s}-\frac{\grad_s f(\bar{X},\bar{s},\bar{y})+u}{L}$ is the
minimizer of the inner minimization problem in
\eqref{ch3_eq:minimax_problem}.

Since $\bar{y}\in\cQ$, the following is true for the third
minimization problem on the right hand side of
\eqref{ch3_eq:combined_upper_bound}.
\begin{align}
\label{ch3_eq:simple_bound_y_problem}
\min_{y\in\cQ\subset\reals^q}\left\{\grad_y f(\bar{X},\bar{s},\bar{y})^T(y-\bar{y})+\frac{L}{2}\norm{y-\bar{y}}_2^2\right\}\leq 0.
\end{align}
Thus, \eqref{ch3_eq:combined_upper_bound},
\eqref{ch3_eq:simple_minimax_problem_alpha},
\eqref{ch3_eq:simple_minimax_problem_beta} and
\eqref{ch3_eq:simple_bound_y_problem} together imply that
\begin{align}
\lambda(\mu_1\norm{\sigma(X_*)}_\alpha+\mu_2\norm{s_*}_\beta)+f(X_*,s_*,y_*)
\leq~&\lambda(\mu_1\norm{\sigma(\bar{X})}_\alpha+\mu_2\norm{\bar{s}}_\beta)+f(\bar{X},\bar{s},\bar{y}) \nonumber \\
&\mbox{}-\min_{W:\norm{\sigma(W)}_{\alphas}\leq \lambda\mu_1}\frac{\norm{\grad_X f(\bar{X},\bar{s},\bar{y})+W}_F^2}{2L}\nonumber\\
&\mbox{}-\min_{u:\norm{u}_{\betas}\leq \lambda\mu_2}\frac{\norm{\grad_s f(\bar{X},\bar{s},\bar{y})+u}_2^2}{2L}.\nonumber
\end{align}
Since
$\Big(\lambda(\mu_1\norm{\sigma(\bar{X})}_\alpha+\mu_2\norm{\bar{s}}_\beta)+f(\bar{X},\bar{s},\bar{y})\Big)
-\Big(\lambda(\mu_1\norm{\sigma(X_*)}_\alpha+\mu_2\norm{s_*}_\beta)+f(X_*,s_*,y_*)\Big)\leq\epsilon$,
we have that
\begin{align}
\min_{W:\norm{\sigma(W)}_{\alphas}\leq \lambda\mu_1}\norm{\grad_X f(\bar{X},\bar{s},\bar{y})+W}_F^2+
\min_{u:\norm{u}_{\betas}\leq \lambda\mu_2}\norm{\grad_s f(\bar{X},\bar{s},\bar{y})+u}_2^2\leq 2L \epsilon. \label{ch3_eq:L_bound}
\end{align}
From % the definition of $I(.)$ in
\eqref{ch3_eq:IJ_beta}, it follows that $\norm{W}_F\leq
I(\alphas)\norm{\sigma(W)}_{\alphas}$. Thus, \eqref{ch3_eq:L_bound}
implies that
\begin{align}
\min_{W:\norm{W}_F\leq I(\alphas)\lambda\mu_1}\norm{\grad_X f(\bar{X},\bar{s},\bar{y})+W}_F^2\leq 2L \epsilon. \label{ch3_eq:L_alpha_bound_2inf}
\end{align}
Suppose $\norm{\grad_X f(\bar{X},\bar{s},\bar{y})}_F> I(\alphas)\lambda\mu_1$. Then the optimal solution of the optimization problem in~\eqref{ch3_eq:L_alpha_bound_2inf} is
\eq
W^*=-I(\alphas)\lambda\mu_1 \cdot \frac{\grad_X f(\bar{X},\bar{s},\bar{y})}{\norm{\grad_X f(\bar{X},\bar{s},\bar{y})}_F}.
\en
Then \eqref{ch3_eq:L_bound} implies that $(\norm{\grad_X f(\bar{X},\bar{s},\bar{y})}_F-I(\alphas)\lambda\mu_1)^2\leq 2L\epsilon$, i.e. $\norm{\grad_X f(\bar{X},\bar{s},\bar{y})}_F\leq \sqrt{2L\epsilon}+I(\alphas)\lambda\mu_1$. This is trivially true when $\norm{\grad_X f(\bar{X},\bar{s},\bar{y})}_F\leq I(\alphas)\lambda\mu_1$. Therefore, we can conclude that always
\eq
\norm{\grad_X f(\bar{X},\bar{s},\bar{y})}_F\leq\sqrt{2L\epsilon}+I(\alphas)\lambda\mu_1.
\en
A similar analysis establishes that $\norm{\grad_s f(\bar{X},\bar{s},\bar{y})}_2\leq\sqrt{2L\epsilon}+J(\betas)\lambda\mu_2$.
\end{proof}

\noindent Now we are ready to prove Lemma~\ref{ch3_cor:grad_norm_bound}.
\begin{proof}
Let $f(X,s,y)=\frac{1}{2} \norm{\mathcal{A}(X)-y-b-\lambda
  \theta_1}_2^2+\frac{1}{2} \norm{\mathcal{C}(X)-s-d-\lambda
  \theta_2}_2^2$ and let $\norm{(X,s,y)}=\sqrt{\norm{X}_F^2+\norm{s}_2^2+\norm{y}_2^2}$, then for any $X_1, X_2
\in\reals^{m\times n}$, $s_1, s_2 \in\reals^p$ and $y_1, y_2
\in\reals^q$, we have
\begin{eqnarray*}
\lefteqn{\norm{\grad f(X_1, s_1, y_1)-\grad f(X_2, s_2, y_2)}^2}&&\\
&=&\left\|\left(
\begin{array}{c}
\grad_X f(X_1,s_1,y_1)-\grad_X f(X_2,s_2,y_2)\\
\grad_s f(X_1,s_1,y_1)-\grad_s f(X_2,s_2, y_2)\\
\grad_y f(X_1,s_1,y_1)-\grad_y f(X_2,s_2, y_2)
\end{array}
\right)\right\|^2, \\
&=&\norm{\grad_X f(X_1,s_1,y_1)-\grad_X f(X_2,s_2,y_2)}_F^2+\norm{\grad_s f(X_1,s_1,y_1)-\grad_s f(X_2,s_2,y_2)}_2^2\\
&&\mbox{}+\norm{\grad_y f(X_1,s_1,y_1)-\grad_y f(X_2,s_2,y_2)}_2^2,\\
&=&\norm{\mathcal{A}^*(\mathcal{A}(X_1-X_2)-y_1+y_2)+\mathcal{C}^*(\mathcal{C}(X_1-X_2)-s_1+s_2)}_F^2\\
&&\mbox{}+\norm{\mathcal{C}(X_1-X_2)-s_1+s_2}_2^2+\norm{\mathcal{A}(X_1-X_2)-y_1+y_2}_2^2,\\
&=&\norm{A^T(A\ \vect(X_1-X_2)-y_1+y_2)+C^T(C\ \vect(X_1-X_2)-s_1+s_2)}_2^2\\
&&\mbox{}+\norm{C\ \vect(X_1-X_2)-s_1+s_2}_2^2+\norm{A\ \vect(X_1-X_2)-y_1+y_2}_2^2,\\
&=&\left\|M^TM
    \left(
     \begin{array}{c}
       s_1-s_2 \\
       y_1-y_2 \\
       \vect(X_1-X_2)\\
     \end{array}
    \right)\right\|^2_2.
\end{eqnarray*}
Hence,
\begin{align}
\norm{\grad f(X_1, s_1, y_1)-\grad f(X_2, s_2, y_2)}\leq
&\ \sigma_{\max}^2(M)~\left\|
    \left(
     \begin{array}{c}
       s_1-s_2 \\
       y_1-y_2 \\
       \vect(X_1-X_2)\\
     \end{array}
    \right)\right\|_2, \nonumber\\
=&\ \sigma_{\max}^2(M)~\sqrt{\norm{X_1-X_2}_F^2+\norm{s_1-s_2}_2^2+\norm{y_1-y_2}_2^2}, \nonumber\\
=&\ \sigma_{\max}^2(M)~\norm{(X_1,s_1,y_1)-(X_2,s_2,y_2)}, \nonumber
\end{align}
where $\sigma_{\max}(M)$ is the maximum singular-value of $M$. Thus, $f:\reals^{m\times n}\times \reals^p\times \reals^q \rightarrow \reals$ is a convex function and $\grad f$ is Lipschitz continuous with respect to $\norm{.}$ with Lipschitz constant $L=\sigma_{\max}^2(M)$.

Since  $(\bar{X},\bar{s},\bar{y})$ is an $\epsilon$-optimal solution to the problem $\min\{P(X,s,y):X\in\reals^{m\times n}, s\in\reals^p, y\in\cQ\subset\reals^q\}$, Theorem~\ref{ch3_thm:grad_convergence} guarantees that
\begin{eqnarray}
\norm{\grad_X f(\bar{X},\bar{s},\bar{y})}_F &=& \norm{\mathcal{A}^*(\mathcal{A}(\bar{X})-\bar{y}-b-\lambda\theta_1) + \mathcal{C}^*(\mathcal{C}(\bar{X})-\bar{s}-d-\lambda\theta_2)}_F \nonumber\\
&\leq& \sqrt{2\epsilon}~\sigma_{\max}(M)+I(\alphas)\lambda\mu_1, \label{ch3_eq:grad_x_norm_bound}\\
\norm{\grad_s f(\bar{X},\bar{s},\bar{y})}_2 &=& \norm{\mathcal{C}(\bar{X})-\bar{s}-d-\lambda\theta_2}_2
\leq \sqrt{2\epsilon}~\sigma_{\max}(M)+J(\betas)\lambda\mu_2. \label{ch3_eq:grad_s_norm_bound}
\end{eqnarray}

\end{proof}

\subsection*{Lemma~\ref{lem:euclidean_proj} and proof} $\mbox{}$

\begin{lemma}
\label{lem:euclidean_proj}
Let $\cQ\subset\reals^q$ be a nonempty, closed, and convex set. Then for all $\tilde{y}\in\reals^q$ and $\lambda>0$, we have
$\Pi_\cQ(\lambda\tilde{y})=\lambda~\Pi_{\cQ /\lambda}(\tilde{y})$, or
equivalently, $\Pi_\cQ(\tilde{y})=\lambda~\Pi_{\cQ
  /\lambda}(\tilde{y}/\lambda)$, where $\cQ/\lambda = \{x: \lambda x
\in \cQ\}$.
\end{lemma}
\begin{proof}
Fix $\tilde{y}\in\reals^q$ and $\lambda>0$. Then
\begin{align}
\Pi_\cQ(\lambda\tilde{y})=\argmin_{x\in\cQ}\norm{x-\lambda\tilde{y}}_2% \argmin_{\lambda
  % y\in\cQ}\norm{\lambda y-\lambda\tilde{y}}_2=\argmin_{\lambda
  % y\in\cQ}\norm{ y-\tilde{y}}_2
=\lambda\argmin_{y\in\cQ/\lambda}\norm{
  y-\tilde{y}}_2 = \lambda~\Pi_{\cQ /\lambda}(\tilde{y}).
\end{align}
\end{proof}

\subsection*{Lemma~\ref{lem:grady_bound} and proof} $\mbox{}$

\begin{lemma}
\label{lem:grady_bound}
Let $(\Xopt,\sopt,\yopt)$ be an optimal solution to
\eqref{ch3_eq:CN-slack} and suppose that
$\norm{\Pi_{\cQ}\left(\yk_p\right)-\yk}_2\leq\xik$ for some $k\geq 1$,
where $\yk_p:=\yk-\frac{1}{L}\grad_y f^{(k)}(\Xk,\sk,\yk)$. Then we have
\begin{align}
-\fprod{\grad_y f^{(k)}(\Xk,\sk,\yk),~\yopt-\yk}\leq
L\xik\norm{\yopt-\yk}_2+\xik\norm{\grad_y
  f^{(k)}(\Xk,\sk,\yk)}_2. \label{eq:grady_bound}
\end{align}
\end{lemma}
\begin{proof}
From the definition of $\Pi_\cQ(.)$, we have
\begin{eqnarray}
&&\fprod{\Pi_\cQ(\yk_p)-\yk_p,~y-\Pi_\cQ(\yk_p)}\geq 0, \quad
\forall~y\in\cQ,\nonumber\\
&\Rightarrow&\fprod{\Pi_\cQ(\yk_p)-\yk,~y-\yk}
+\fprod{\Pi_\cQ(\yk_p)-\yk,~\yk-\Pi_\cQ(\yk_p)}\nonumber\\
&&+\fprod{\yk-\yk_p,~y-\yk}+\fprod{\yk-\yk_p,~\yk-\Pi_\cQ(\yk_p)}\geq 0,
\quad \forall~y\in\cQ. \label{eq:grady_bound_1}
\end{eqnarray}
Since $\yopt\in\cQ$, $\yk-\yk_p=\frac{1}{L}\grad_y f^{(k)}(\Xk,\sk,\yk)$
and $\norm{\Pi_{\cQ}\left(\yk_p\right)-\yk}_2\leq\xik$,
\eqref{eq:grady_bound} follows from \eqref{eq:grady_bound_1}.
\end{proof}

\section{Auxiliary results for simple optimization problems}
\label{app:aux-results}
\begin{lemma}
\label{ch3_lem:boundry_solution}
Let $(\cE,\norm{.})$ be a normed vector space,
$f:\cE\rightarrow \reals$ be a strictly convex function and
$\chi\subset\cE$ be a closed, convex set with a non-empty interior. Let
$\bar{x}=\argmin_{x\in\chi}f(x)$ and
$x^*=\argmin_{x\in\cE}f(x)$. If $x^*\not\in \chi$, then
$\bar{x}\in\bdry \chi$, where $\bdry \chi$ denotes the boundary of $\chi$.
\end{lemma}
\begin{proof}
We will establish the result by contradiction. Assume $\bar{x}$ is in
the interior of $\chi$, i.e. $\bar{x}\in \intr(\chi)$. Then
$\exists\;\epsilon>0$ such that
$B(\bar{x},\epsilon)=\{x\in\cE\;:\;\norm{x-\bar{x}}<\epsilon\}\subset
\chi$. Since $f$ is strictly convex and $x^*\neq \bar{x}$,
$f(x^*)<f(\bar{x})$. Choose
$0<\lambda<\frac{\epsilon}{\norm{\bar{x}-x^*}}<1$ so that $\lambda x^*
+(1-\lambda) \bar{x} \in B(\bar{x},\epsilon) \subset \chi$.
Since $f$ is strictly convex,
\begin{align}
f(\lambda x^* + (1-\lambda) \bar{x})<\lambda f(x^*) + (1-\lambda)
f(\bar{x}) < f(\bar{x}).
\end{align}
However, $\lambda x^* + (1-\lambda) \bar{x} \in
B(\bar{x},\epsilon)\subset \chi$ and $f(\lambda x^* + (1-\lambda)
\bar{x})<f(\bar{x})$
contradicts the fact that $f(\bar{x})<f(x)$ for all $x\in
\chi$. Therefore, $\bar{x} \not \in \intr(\chi)$. Since $\bar{x}\in
\chi$, it follows that $\bar{x}\in \bdry \chi$.
\end{proof}

Next, we collect together complexity results for optimization problems of
the form
\begin{eqnarray*}
&\min_{X\in\reals^{m\times n}}&\{\lambda \norm{\sigma(X)}_{\alpha}+ \frac{1}{2}\norm{X-\tilde{X}}_F^2: \norm{\sigma(X)}_{\alpha} \leq \eta\}\\
&\min_{s\in\reals^p}&\{\lambda\norm{s}_{\beta} +
\frac{1}{2}\norm{s-\tilde{s}}_2^2:\norm{s}_{\beta} \leq \eta\}
\end{eqnarray*}
that need to be solved in each \textbf{Algorithm~APG} update step, displayed in Figure~\ref{ch1_alg:pga}.
\begin{lemma}
  \label{ch3_lem:unified_constrained_shrinkage}
  Let $\bar{X} = \argmin_{X \in \reals^{m\times n}}\big\{\lambda \norm{\sigma(X)}_{\alpha} + \frac{1}{2}\norm{X-\tilde{X}}_F^2: \norm{\sigma(X)}_{\alpha} \leq \eta\big\}$ of the constrained matrix shrinkage problem. Then
  \eq
  \bar{X} = U \diag(\bar{s}) V^T,
  \en
  where $U\diag(\sigma)V^T$ denotes the SVD of $\tilde{X}$ such that $\sigma \in \reals_+^r$ and $r=\Rank(\tilde{X})$; and $\bar{s}$ denotes the optimal solution of the constrained vector shrinkage problem
  \eq
  \min_{s \in \reals^r}\big\{\lambda\norm{s}_{\alpha} +   \frac{1}{2}\norm{s-\sigma}_2^2:~\norm{s}_{\alpha} \leq \eta\big\}.
  \en
  Since the worst
  case complexity of computing the SVD of $\tilde{X}$ is
  $\cO(\min\{n^2m,m^2n\})$ the complexity of the computing $\bar{X}$ is
  $\cO(\min\{n^2m,m^2n\} + T_v(r,\alpha))$, where $T_v(r,\alpha)$ denotes the
  complexity of computing the solution of an $r$-dimensional constrained vector
  shrinkage problem with norm $\norm{.}_{\alpha}$.
  The function
  \begin{equation}
    T_v(p,\alpha) =
    \left\{
      \begin{array}{ll}
        \cO(p\log(p)) & \alpha = 1, \infty,\\
        % worst-case complexity and $\cO(n)$
        % expected complexity with randomized search
        \cO(p), & \alpha = 2,\\ % worst-case complexity.
       %  \cO(p\log(p)), & \beta =  % worst-case complexity and $\cO(n)$
% %     expected complexity with randomized search
      \end{array}
      \right.
    \end{equation}
  \end{lemma}
\begin{proof}
  The standard results in non-linear convex optimization over matrices implies
  that $\bar{X}$ is of the form $\bar{X} = U \diag(\bar{s}) V^T$~(see Corollary 2.5 in~\cite{Lewis95_1J}).

  Now, consider the  vector constrained shrinkage problem
  $$\min_{s\in\reals^p}\big\{\lambda\norm{s}_{\beta} + \frac{1}{2}\norm{s-\tilde{s}}_2^2:~\norm{s}_{\beta} \leq \eta\big\}.$$
  \begin{enumerate}[(i)]
  \item $\beta = 1$: First considered the unconstrained case, i.e. $\eta =\infty$. The unconstrained solution $s^*$ has a closed form $s^*=\mbox{sign}(\tilde{s})\odot\max\{|\tilde{s}|-\lambda\mathbf{1}, \mathbf{0}\}$ and can be computed with $\cO(p)$ complexity, where $\odot$ denotes componentwise multiplication and $\mathbf{1}$ is a vector of ones.

    When $\eta<\infty$, the constrained optimal solution, $\bar{s}$, can be computed with $\cO(p\log(p))$ complexity. See Lemma A.4 in~\cite{Ser09_1J}.
  \item $\beta = 2$: First considered the unconstrained case, i.e. $\eta =\infty$.
    Since $\ell_2$-norm is self dual, $\lambda\norm{s}_2 = \max\{u^Ts: \norm{u}_2
    \leq 1\}$. Thus,
    \begin{eqnarray}
      \min_{s\in\reals^p}\left\{\lambda\norm{s}_2+\frac{1}{2}\norm{s-\tilde{s}}_2^2\right\}
      & =& % \min_{s\in\reals^p}\left\{\lambda\norm{s}_\infty +
      % \frac{1}{2}\norm{s-q}_2^2\right\},  \nonumber\\
      % =&
      \min_{s\in\reals^p}\ \max_{u:\ \norm{u}_2\leq \lambda} \left\{u^T s  +
        \frac{1}{2}\norm{s-\tilde{s}}_2^2\right\}, \nonumber\\
      & =&\max_{u:\ \norm{u}_2\leq \lambda}\ \min_{s \in \reals^p}\left\{u^Ts
        +\frac{1}{2}\norm{s-\tilde{s}}_2^2\right\}, \nonumber\\
      & =&\max_{u:\ \norm{u}_2\leq \lambda} \left\{u^T
        (\tilde{s}-u)+\frac{1}{2}\norm{u}_2^2\right\}, \label{ch3_eq:u-2}\\
      & =&\frac{1}{2}\norm{\tilde{s}}_2^2-\min_{u:\ \norm{u}_2\leq \lambda}
      \frac{1}{2}\norm{u-\tilde{s}}_2^2, \nonumber
    \end{eqnarray}
    where \eqref{ch3_eq:u-2} follows from the fact that $s^*(u):=\argmin_{s\in\reals^p} \{u^T
    s+\frac{1}{2}\norm{s-\tilde{s}}_2^2\}=\tilde{s}-u$.

    Define
    \begin{eqnarray*}
      u^*:=\argmin_{u:\ \norm{u}_2\leq \lambda}
      \frac{1}{2}\norm{u-\tilde{s}}_2^2=\tilde{s}~\min\left\{\frac{\lambda}{\norm{\tilde{s}}_2},\
        \mathbf{1}\right\}.
    \end{eqnarray*}
    Then the unconstrained optimal solution $s^*=s^*(u^*)=\tilde{s}\max\left\{1-\frac{\lambda}{\norm{\tilde{s}}_2},\ 0\right\}$ and
    the complexity of computing $\bar{s}$ is $\cO(p)$.

    Next, consider the constrained optimization problem, i.e. $\eta<\infty$.
    The constrained optimum $\bar{s}=s^\ast$,
    whenever $s^*$ is feasible, i.e. $\norm{s^*}_2 \leq \eta$.
    Since $f(s):=\lambda\norm{s}_2+\frac{1}{2}\norm{s-\tilde{s}}_2^2$ is
    strongly convex, Lemma~\ref{ch3_lem:boundry_solution} implies that $\norm{\bar{s}}_{2}
    = \eta$ whenever $\norm{s^*}_2 > \eta$. Thus,
    \eq
    \min \Big\{\lambda\norm{s}_2+\frac{1}{2}\norm{s-\tilde{s}}_2^2:\
    \norm{s}_2\leq \eta\Big\}  = \lambda \eta + \min\Big\{
    \frac{1}{2}\norm{s - \tilde{s}}_2^2 : \norm{s}_2^2 = \eta^2\Big\}.
    \en
    The unique KKT point for the optimization problem $\min\big\{
    \frac{1}{2}\norm{s - \tilde{s}}_2^2 : \frac{1}{2}\norm{s}_2^2 = \frac{\eta^2}{2}\big\}$, is given by
    $\bar{s} = \eta\frac{\tilde{s}}{\norm{\tilde{s}}}$ and KKT multiplier for
    the constraint $\frac{1}{2}\norm{s}_2^2 = \frac{\eta^2}{2}$ is
    $\vartheta =  \frac{\norm{\tilde{s}}_2}{\eta} -1$. It is easy to check that $\vartheta
    > 0$ whenever $\norm{s^*}_2 > \eta$. Thus,
    $\bar{s}$ is optimal for the convex optimization problem $\min\Big\{
    \frac{1}{2}\norm{s - \tilde{s}}_2^2 : \norm{s}_2^2 \leq \eta^2\Big\}$;
    consequently, optimal for equality constrained optimization problem $\min\big\{
    \frac{1}{2}\norm{s - \tilde{s}}_2^2 : \norm{s}_2 = \eta\big\}$. Hence, the
    complexity of computing $\bar{s}$ is $\cO(p)$
  \item $\beta = \infty$: First consider the unconstrained problem.  Since
    $\ell_1$-norm is the dual norm of the $\ell_{\infty}$-norm, we
    have that
    \begin{eqnarray}
      \min_{s\in\reals^p}\left\{\lambda\norm{s}_\infty+\frac{1}{2}\norm{s-\tilde{s}}_2^2\right\}
      & =& % \min_{s\in\reals^p}\left\{\lambda\norm{s}_\infty +
           % \frac{1}{2}\norm{s-q}_2^2\right\}, \nonumber\\
      % =&
      \min_{s\in\reals^p}\ \max_{u:\ \norm{u}_1\leq \lambda} \left\{u^T
        s+\frac{1}{2}\norm{s-\tilde{s}}_2^2\right\}, \nonumber\\
      & =&\max_{u:\ \norm{u}_1\leq \lambda}\ \min_{s \in
        \reals^p}\left\{u^Ts+\frac{1}{2}\norm{s-\tilde{s}}_2^2\right\}, \nonumber\\
      & =&\max_{u:\ \norm{u}_1\leq \lambda} \left\{u^T (\tilde{s}-u) +
        \frac{1}{2}\norm{u}_2^2\right\}, \label{ch3_eq:u-1}\\
      & =&\frac{1}{2}\norm{\tilde{s}}_2^2-\min_{u:\ \norm{u}_1\leq \lambda} \frac{1}{2}\norm{u-\tilde{s}}_2^2, \nonumber
    \end{eqnarray}
    where \eqref{ch3_eq:u-1} follows from the fact that $s^*(u):=\argmin_{s\in\reals^p}\{u^T
    s+\frac{1}{2}\norm{s-\tilde{s}}_2^2\}=\tilde{s}-u$.
    The result in (i) implies that
    complexity of computing $u^*=\min_{u:\ \norm{u}_1\leq \lambda}
    \frac{1}{2}\norm{u-\tilde{s}}_2^2$ is $\cO(p\log(p))$. Thus, the
    unconstrained optimal solution  $s^* = s^*(u^*) =
    \tilde{s} - u^\ast$ can be computed in $\cO(p\log(p))$ operations.

    Next, consider the constrained optimization problem.
    The constrained optimum, $\bar{s}=s^\ast$
    whenever $s^*$ is feasible, i.e. $\norm{s^*}_{\infty} \leq \eta$.
     Since $f(s) = \lambda \norm{s}_{\infty} + \frac{1}{2}\norm{s-\tilde{s}}_2^2$ is strictly convex,
    Lemma~\ref{ch3_lem:boundry_solution} implies that $\norm{\bar{s}}_{\infty}
    = \eta$, whenever $\norm{s^*}_{\infty} > \eta$. Therefore,
    \eq
    \min\Big\{\lambda\norm{s}_\infty+\frac{1}{2}\norm{s-\tilde{s}}_2^2:\norm{s}_\infty
    \leq \eta\Big\}  = \lambda \eta +
    \min\Big\{\frac{1}{2}\norm{s-\tilde{s}}_2^2:\norm{s}_\infty =\eta\Big\}.
    \en
    Then, it is easy to check $\mbox{sign}(\bar{s}_i) = \mbox{sign}(\tilde{s}_i)$ for all $i = 1, \ldots, p$. Moreover, $\norm{s^*}_\infty>\eta$ implies that $\norm{\tilde{s}}_\infty>\eta$. These two facts imply that
    \eq
    \min\Big\{\frac{1}{2}\norm{s-\tilde{s}}_2^2:\norm{s}_\infty =\eta\Big\} =
    \min\Big\{\frac{1}{2}\norm{s-\abs{\tilde{s}}}_2^2: 0 \leq s_i \leq \eta\Big\}.
    \en
    For $1\leq i\leq p$, we have $\min\{\abs{\tilde{s}_i},\eta\}=\argmin_{s_i\in\reals}\big\{\frac{1}{2}(s_i - \abs{\tilde{s}_i})^2: 0
    \leq s_i \leq \eta\big\}$. Thus, it follows that $\bar{s}= \mbox{sign}(\tilde{s}) \odot\min\{|\tilde{s}|, \eta \mathbf{1}\}$. Hence the complexity of computing $\bar{s}$ is $\cO(p\log(p))$.
  \end{enumerate}
\end{proof}

\end{document}